\documentclass[11pt,a4paper]{article}
\usepackage[english]{babel}
\usepackage[utf8]{inputenc}
\usepackage{amsbsy}
\usepackage{amsmath}
\usepackage{amssymb}
\usepackage{amsfonts}
\usepackage{amsthm}
\usepackage{mathtools}
\mathtoolsset{showonlyrefs}  
\usepackage{bm}
\usepackage{esint}
\usepackage[percent]{overpic}
\usepackage{graphicx}
\usepackage{hyperref}
\usepackage[active]{srcltx}
\usepackage{upgreek}
\usepackage{xcolor}
\usepackage{srcltx}
\usepackage{psfrag}
\usepackage[percent]{overpic}
\usepackage{subcaption} 
\usepackage[nottoc,numbib]{tocbibind}
\usepackage{dsfont}
\usepackage[skip=5pt]{caption}
\usepackage[shortlabels]{enumitem}
\setlist[enumerate]{topsep=1pt,itemsep=0pt,parsep=0pt}
\setenumerate[1]{label=\textup{(\roman*)}}
\setenumerate[2]{label=(\alph*)}
\setenumerate[3]{label=(\arabic*)}

\newcommand{\beps}{\bm{\varepsilon}}

\newcommand{\C}{\mathbb{C}}
\newcommand{\N}{\mathbb{N}}

\newcommand{\R}{\mathbb{R}}
\renewcommand{\SS}{\mathbb{S}}

\newcommand{\Z}{\mathbb{Z}}

\newcommand{\boA}{\mathcal{A}}

\newcommand{\boC}{\mathcal{C}}

\newcommand{\boE}{\mathcal{E}}

\newcommand{\boH}{\mathcal{H}}

\newcommand{\boK}{\mathcal{K}}

\newcommand{\boN}{\mathcal{N}}
\newcommand{\boO}{\mathcal{O}}
\newcommand{\boP}{\mathcal{P}}

\newcommand{\boR}{\mathcal{R}}
\newcommand{\boS}{\mathcal{S}}

\newcommand{\boV}{\mathcal{V}}

\newcommand{\gE}{\mathfrak{E}}

\newcommand{\gS}{\mathfrak{S}}

\newcommand{\ga}{\mathfrak{a}}
\newcommand{\gb}{\mathfrak{b}}
\newcommand{\gc}{\mathfrak{c}}

\newcommand{\gm}{\mathfrak{m}}
\newcommand{\gn}{\mathfrak{n}}

\newcommand{\gs}{\mathfrak{s}}
\newcommand{\gu}{\mathfrak{u}}
\newcommand{\gv}{\mathfrak{v}}

\newcommand{\ve}{\varepsilon}
\newcommand{\grad}{\nabla}
\newcommand{\ptl}{\partial}

\renewcommand{\P}{\mathcal{P}}
\newcommand{\M}{M} 
\newcommand{\SSS}{\bm{S}}
\newcommand{\RRR}{\bm{R}}
\newcommand{\A}{\bm{A}}
\newcommand{\B}{\bm{B}}
\renewcommand{\u}{\bm{u}}
\renewcommand{\v}{\bm{v}}

\newcommand{\m}{\bm{m}}

\newcommand{\mm}{\bm{m}}

 \newcommand{\QQ}{\boldsymbol{Q}}

\newcommand{\ba}{\bm{a}}

\newcommand{\be}{\bm{e}}
\newcommand{\bmm}{\bm{m}}

\newcommand{\bv}{\bm{v}}
\newcommand{\bu}{\bm{u}}

\newcommand{\f}{\bm{f}}
\newcommand{\g}{\bm{g}}
\newcommand{\h}{\bm{h}}

\renewcommand{\div}{\mathop{{\rm div}}\nolimits}
\DeclareMathOperator{\diag}{{\rm diag}}

\renewcommand{\Im}{\mathop{{\rm Im}}\nolimits}

\renewcommand{\Re}{\mathop{{\rm Re}}\nolimits}

\newcommand{\loc}{\mathrm{loc}}
 
\providecommand{\abs}[1]{|#1 |}
\providecommand{\norm}[1]{\|#1 \|}
\renewcommand{\div}{\operatorname{div}}
\newcommand{\sech}{\operatorname{sech}}
\newcommand{\dist}{\operatorname{dist}}

\newcommand{\Erf}{\operatorname{Erf}}
 \newcommand{\BMO}{\textup{BMO}}
\newcommand{\bqq}{\begin{equation*}}
\newcommand{\eqq}{\end{equation*}}
\newcommand{\bq}{\begin{equation}}
\newcommand{\eq}{\end{equation}}


\theoremstyle{plain}
\newtheorem{thm}{Theorem}[section]
\newtheorem{cor}[thm]{Corollary}
\newtheorem{corollary}[thm]{Corollary}

\newtheorem{prop}[thm]{Proposition}
\newtheorem*{thm*}{Theorem}
\theoremstyle{definition}

\newtheorem{remark}[thm]{Remark}

 \usepackage{microtype}
\begin{document}
	\title{Recent results for the Landau--Lifshitz equation}
	\author{
		\renewcommand{\thefootnote}{\arabic{footnote}}
Andr\'e de Laire\footnotemark[1]}
	\footnotetext[1]{ 
		Univ.\ Lille, CNRS, UMR 8524, Inria - Laboratoire Paul Painlev\'e, F-59000 Lille, France.\\
		E-mail: {\tt andre.de-laire@univ-lille.fr}}
	\date{}
	\maketitle
\begin{abstract}
We give a survey on some recent results concerning the Landau--Lifshitz equation, a
fundamental nonlinear PDE with a strong geometric content, describing the dynamics of  the magnetization  in ferromagnetic materials.
We revisit the Cauchy problem for the anisotropic LL equation,
without dissipation, for smooth solutions, and also in the energy space 
in  dimension one.
We also  examine two approximations of the LL equation given by 
of the Sine--Gordon equation and cubic  Schr\"odinger equations, arising in certain singular limits of strong easy-plane and easy-axis anisotropy, respectively.

Concerning localized solutions, we review the orbital and asymptotic stability problems for a
sum of solitons in dimension one, exploiting  the variational nature of the solitons in the  hydrodynamical framework. 

Finally, we survey results concerning the existence, uniqueness  and stability
of self-similar solutions (expanders and shrinkers) for the isotropic  LL equation with Gilbert term. Since expanders are associated with a singular initial condition
with a jump discontinuity, we also review their well-posedness in spaces linked to the BMO space.
\end{abstract}

\section{Introduction}

The Landau--Lifshitz (LL) equation has been  introduced in 1935 by L.~Landau and E.~Lifshitz in~\cite{LandLif1}   and it  constitutes nowadays a fundamental tool in the magnetic recording industry, 
due to its applications to ferromagnets \cite{wei2012}.
This PDE  describes the dynamics of the orientation of the 
magnetization (or spin) in ferromagnetic materials, 
and it is given by 
\begin{equation}\label{LL-general}
	\partial_t \bmm+  \bmm \times H_{\text{eff}}(\bmm)=0, 
\end{equation}
where  $\bmm=(m_1, m_2,m_3):\R^N \times I\longrightarrow \SS^2$
is the spin vector, $I\subset \R$ is a time interval,  $\times$ denotes
the usual cross-product in $\R^3$, and $\mathbb{S}^2$ is the unit
sphere in $\mathbb{R}^3$. Here  $H_{\text{eff}}(\bmm)$ is
the effective magnetic field, corresponding to (minus) the $L^2$-derivative of the  magnetic energy of 
the material. We will focus on energies of the form 
$	E_{\rm LL}(\bmm) =E_{\rm ex}(\bmm)+E_{\rm ani}(\bmm), $
where the {\em exchange} energy 
\begin{equation}
	E_{\rm ex}(\bmm)= \frac{1}{2} \int_{\R^N}  |\nabla \bmm|^2=  \frac{1}{2} \int_{\R^N}  |\nabla m_1|^2+|\nabla m_2|^2+|\nabla m_3|^2,
\end{equation}
accounts for the  local tendency of $\m$ to align the magnetization field, and the {\em anisotropy} energy
\begin{equation}
	E_{\rm ani}(\bmm)=\frac12 \int_{\R^N}\langle \bmm , J\bmm\rangle_{\R^3}, \ \quad J\in \textup{Sym}_3(\R), 
\end{equation}
accounts for the likelihood of $\bmm$ to attain one or more directions of magnetization, 
which determines the {\em easy} directions.  Due to the invariance of \eqref{LL-general} 
under rotations,
 we can assume that $J$ is a diagonal matrix  $J = \diag(J_1, J_2, J_3)$,
and thus the anisotropy energy reads
\begin{equation}
	\label{Eani}
	E_{\rm ani}(\bmm)=\frac12 \int_{\R^N} ( \lambda_1 m_1^2 + \lambda_3 m_3^2),
\end{equation}
with  $\lambda_1 = J_2 - J_1$ and $\lambda_3 = J_2 - J_3$. Therefore \eqref{LL-general} can be recast as
\begin{equation}
	\label{LL-ani}
	\partial_t \bmm+  \bmm \times (\Delta \bmm-\lambda_1m_1 \be_1-\lambda_3m_3 \be_3)=0,
\end{equation}
where $(\be_1,\be_2,\be_3)$ is the canonical basis of $\R^3$.
Notice that for finite energy solutions, \eqref{Eani} formally implies that 
$m_1(x)\to0$ and $m_3(x)\to0$, as $\abs{x}\to \infty$, and hence 
$\abs{m_2(x)}\to 1$, as $\abs{x}\to \infty.$

For biaxial ferromagnets, all the numbers $J_1$, $J_2$ and $J_3$ are different, so that $\lambda_1\neq \lambda_3$ and $\lambda_1 \lambda_3 \neq 0$. Uniaxial ferromagnets are characterized by the property that only two of the numbers $J_1$, $J_2$ and $J_3$ are equal. For instance, the case  $J_1 = J_2$ corresponds to $\lambda_1 = 0$ and $\lambda_3 \neq 0$, so that the material has a uniaxial anisotropy in the direction  $\be_3$. Hence, the ferromagnet owns an {\em easy-axis} anisotropy along the vector $\be_3$ if $\lambda_3 < 0$, while the anisotropy is {\em easy-plane} along the plane $x_3 = 0$ if $\lambda_3 > 0$. 
Finally, in the isotropic case $\lambda_1 = \lambda_3 = 0$,  equation \eqref{LL-ani} reduces to the well-known Schr\"odinger map equation
\begin{equation}
	\label{SM}
	\partial_t \bmm+  \bmm \times \Delta \bmm=0. 
\end{equation}

The LL equation \eqref{LL-ani} is a nonlinear dispersive PDE, with dispersion relation 
\begin{equation}
	\label{dispersion}
	\omega(k)=\pm\sqrt{\abs{k}^4+(\lambda_1+\lambda_3)\abs{k}^2+\lambda_1\lambda_3},
\end{equation}
for linear sinusoidal waves of
frequency $\omega$ and wavenumber $k$, i.e.~solutions of the form $e^{i(k \cdot x- \omega t)}$.
From \eqref{dispersion}, we can recognize similarities with some classical dispersive equations. For instance, for the Schr\"odinger equation 
$	i\partial_t \psi+\Delta \psi= 0,$
the dispersion relation is $\omega(k)=\abs{k}^2$, corresponding to $\lambda_1=\lambda_3=0$ in \eqref{dispersion}, i.e.\ the Schr\"odinger map equation \eqref{SM}.

When considering  Schr\"odinger equations with nonvanishing conditions at infinity, the typical
example is the Gross--Pitaesvkii equation \cite{deLaire2}
\begin{equation}
	\label{GP}
	i\partial_t \psi+\Delta \psi+\sigma \psi (1-\abs{\psi}^2)= 0, 
\end{equation}
$\sigma >0$, and the dispersion relation for the linearized equation at the constant solution equal to $1$ is
	$\omega(k)=\pm\sqrt{\abs{k}^4+2\sigma\abs{k}^2}.$
This  corresponds to taking $\lambda_1=0$ or $\lambda_3=0$, with $\lambda_1+\lambda_3=2\sigma$, in \eqref{dispersion}.

Finally, let us consider the Sine--Gordon equation
$
\partial_{tt} \psi-\Delta \psi+\sigma \sin(\psi)= 0,
$
$\sigma>0$, whose linearized equation at 0 is given by the Klein--Gordon equation, with dispersion relation
	$\omega(k)=\pm\sqrt{\abs{k}^2+\sigma},$
that behaves like \eqref{dispersion} for $\lambda_1\lambda_3=\sigma$ 
and $\lambda_1+\lambda_3=1$, at least for $k$ small. 

In this context, the Landau--Lifshitz equation is considered as a universal model from which it is possible to derive other completely integrable equations \cite{FaddTak0}. We review some recent  rigorous results in this context in Section~\ref{sec:regime}.

%
%
%

\subsection{The dissipative model}
\label{preamble:dissipative}

In 1955, T.~Gilbert  proposed in \cite{gilbert} a modification of  equation \eqref{LL-general} to incorporate a damping term. The so-called Landau--Lifshit--Gilbert (LLG) equation then reads
\begin{equation}
	\label{LLG-Heff}
	\ptl_t \bmm=-\beta \bmm\times H_{\text{eff}}(\bmm) -\alpha \bmm \times (\bmm\times H_{\text{eff}}(\bmm)),
\end{equation}
where  $\beta\geq 0$ and  $\alpha\geq 0$, so that there is dissipation when $\alpha>0$, and in that case we refer to $\alpha$ as the Gilbert damping coefficient.
Note that, by performing a time scaling, we assume w.l.o.g.\ that
$$
\alpha\in[0,1]\quad \text{and}\quad  \beta=\sqrt{1-\alpha^2}.
$$ 

Let us remark that the identity $a\times (b\times  c)=  b(a\cdot  c)- c( a\cdot  b)$, for all $ a, b, c\in \R^3$, implies that 
for any smooth function $\bv$, valued in $\SS^2$, satisfies 
\begin{equation}
	\label{iden:cross-prod}
-	\bv\times (\bv\times \Delta \bv)=\Delta \bv+\abs{\grad \bv}^2\bv.
	\end{equation}
Then, we see  that  in the limit case  $\beta = 0$ (and so $\alpha=1$), 
the LLG equation reduces to the heat-flow equation for
harmonic maps
\begin{equation}\label{HFHM}
	\ptl_t \bmm-\Delta\bmm=\abs{\grad{\bmm}}^2\bmm.
\end{equation}
This classical equation is an important model in several areas such as differential geometry
and calculus of variations. It is also related with other problems 
such as the theory of liquid crystals and the Ginzburg--Landau equation. 
For more details, we refer to the surveys \cite{eells2,lin,struwe96}. 

As before, one way to start the study of the LLG equation
is noticing the link with other PDEs. Let us 
illustrate this point in the isotropic case $H_{\text{eff}}(\bmm)=\Delta \bmm$.
To simplify our notation, we consider the equation for the opposite vector  $\bmm\to -\bmm$, 
which yields the equation
\begin{equation}
	\label{LLG-intro}
	\ptl_t \bmm=\beta \bmm\times \Delta\bmm -\alpha \bmm \times (\bmm\times \Delta \bmm).
\end{equation}
For a smooth solution $\bmm$ with \mbox{$m_3>-1$}, 
we can use the stereographic projection 
\bq 
\label{def:stereo}
u=\boP(\m)=\frac{m_1+im_2}{1+m_3},
\eq
that satisfies the quasilinear Schr\"odinger equation
	\begin{equation}\label{DNLS}\tag{DNLS}
	iu_t +(\beta-i\alpha)\Delta u=2(\beta-i\alpha)\frac{\bar u (\grad u)^2}{1+\abs{u}^2},
\end{equation}
where we used the notation $(\grad u)^2=\grad u \cdot \grad u =\sum_{j=1}^N (\partial_{x_j}u)^2$ (see e.g.\ \cite{lak-nak} for details).
When $\alpha>0$, one can use the properties of the  semigroup $e^{(\alpha+i\beta)t\Delta}$ to establish a Cauchy theory for rough initial data, as we will see in Section~\ref{chap:self}.

When $N=1$, the LLG equation is also related to the
Localized Induction Approximation (LIA), also called 
{\it{binormal flow}}, a geometric curve flow modeling the
self-induced motion of a vortex filament within an inviscid fluid in
$\R^3$ \cite{lakshmanan,daniel-lak}. As we will in Section~\ref{chap:self}, this is related with the geometric representation of the LLG equation in a Serret--Frenet system.


There are several variants of previous equation  considering more complex models 
including for instance a demagnetization field and the 
effects of the boundary in bounded domains. We refer to \cite{lakshmanan} for
an  overview of different models, to \cite{cimrak} for 
recent developments on the approximation of  solutions,
to the survey \cite{GuoDing0}
for more details of the  derivation and results on the initial value problem, 
and to  \cite{otto-kohn2006} for a review of methods for  pattern formation based on asymptotic analysis. 

\subsection{The hydrodynamical formulation}
\label{sec:hidro}
We end this introduction by explaining  another useful transformation for the analysis of the LL equation.  For simplicity, we assume that there is no dissipation. In the seminal work \cite{Madelun1}, Madelung showed that 
the nonlinear Schr\"odinger equation (NLS) can be recast into the form of a hydrodynamical system.
For instance, for the   NLS equation 
\begin{equation}
	i\partial_t \Psi+\Delta \Psi +\Psi f(\abs{\Psi}^2)=0,
\end{equation}
assuming that $\rho=\abs{\Psi}^2$  does not vanish, the  {Madelung transform}
$\psi =\sqrt{\rho}e^{i\phi}$
leads to the system 
\begin{equation*}
\partial_t \rho+2\div(\rho  \grad \phi )=0,\quad 
\partial_t \phi+\abs{\grad\phi}^2+f(\rho)=
		\frac{\Delta(\sqrt{\rho})}{\sqrt{\rho}}.
\end{equation*}
Therefore, setting $\bv=2\grad \phi$, we get 
the  Euler--Korteweg system
\begin{equation}\label{NLS-hydro2}
\partial_t \rho+\div(\rho \bv)=0,\quad 
\partial_t \bv+(\bv\cdot \grad) \bv+2\grad (f(\rho))=2 \grad\Big(
		\frac{\Delta(\sqrt{\rho})}{\sqrt{\rho}}
		\Big),
\end{equation}
which 
is a dispersive perturbation of the classical  Euler equation for compressible fluids, 
with the additional term 
$2\grad (
\Delta(\sqrt{\rho})/\sqrt{\rho}$, 
which is interpreted as quantum pressure in the quantum fluids models
\cite{CarDaSa1,ingrid-fluids}.

The Madelung transform is  useful to study properties of NLS equations
with nonvanishing conditions at infinity  (see \cite{BetGrSm1,ChirRou2}).
Coming back to the LL equation \eqref{LL-ani}, let $\bmm$ be a solution 
of this equation such that 
the map $\check{m} = m_1 + i m_2$  does not vanish.  
In the spirit of the Madelung transform,  we set 
$$\check{m} = (1 - m_3^2)^\frac{1}{2} \big( \sin(\phi) + i\cos(\phi) \big).$$
Thus, setting the  hydrodynamical variables $u = m_3$ and $\phi$, we get  the system
\begin{equation}
	\tag{H}
	\label{HLL}
	\left\{
	\begin{aligned}
		\partial_t u &= \div \big( (1 - u^2) \nabla \phi \big) - \frac{\lambda_1}{2} (1 - u^2) \sin(2 \phi),\\
		\partial_t \phi &= - \div \Big( \frac{\nabla u}{1 - u^2} \Big) + u \frac{|\nabla u|^2}{(1 - u^2)^2} - u |\nabla \phi|^2 + u \big( \lambda_3 - \lambda_1 \sin^2(\phi) \big),
	\end{aligned}
	\right.
\end{equation}
at long  as
$|u| < 1$ on $\R^N$. As shown in the next sections, the hydrodynamical formulation will be essential in the study of solutions of the LL equation.

Although it does not quite have the reputation of e.g.\ the Navier--Stokes equation or the Ricci flow equation, it can be said that the LL equation is among the most intriguing and challenging PDEs. The mathematical appeal relies on the combination of difficulties from nonlinear Schrödinger equations and geometric evolution equations.  The aim of this note is to survey some recent results concerning the different  aspects of the LL equation, as follows.
In Section~\ref{chap:cauchy} we revisit the Cauchy problem for the anisotropic LL equation,
without dissipation. Concerning smooth solutions, the approach follows a methodology
for quasilinear hyperbolic systems based on a priori estimates by using new well-tailored  higher order energies. We also tackle a subtle well-posedness problem in one space dimension in the energy space by invoking the hydrodynamical formulation.

Section~\ref{sec:regime}  examines approximations of the Landau--Lifshitz equation by the Sine--Gordon equation and cubic  Schr\"odinger equations arising in certain singular limits of large easy-plane and easy-axis anisotropy, respectively,  providing  
quantitative convergence results. 

In Section~\ref{chap:stability} we review the orbital and asymptotic stability problems for
sum of solitons and multisolitons for the easy-plane (undamped) LL equation in dimension one.  
Stability problems of this kind are well-established  in the context
of dissipative evolution equations. 
Here the hamiltonian structure plays an essential role, that we exploit in the  hydrodynamical framework. The essential idea is to exploit the variational structure given by the energy and  momentum so that stability is essentially captured
by spectral bounds for the hessian of the combined functional.

Finally, in Section~\ref{chap:self} we consider the (isotropic) dissipative LLG equation. We focus mainly on the one-dimensional  analysis of  self-similar solutions:  expanders and shrinkers evolving from or towards a singular time.
 We survey results concerning their existence and uniqueness by using a moving frame argument that allows us to obtain the asymptotics of the profiles. We also consider the question of stability of expanders that calls for a well-posedness
 result for  solutions with rough initial data.

\section{The Cauchy problem for the LL equation}
\label{chap:cauchy}
Despite some serious efforts to establish a complete Cauchy theory 
for the LL equation, several issues remain unknown. 
In this section  we will focus on the LL equation without damping, for which 
the Cauchy theory is even more delicate to handle. Even in the case where the problem is isotropic, i.e.\ the Schr\"odinger map equation, 
there are several unknown aspects.
Moreover, it is not always possible to adapt results for 
Schr\"odinger map equation to include  anisotropic perturbations. 


The study of  well-posedness in  the presence of a  damping term
is different. Indeed, for the LLG equation,  some techniques related to parabolic equations and for the heat-flow for harmonic maps \eqref{HFHM} can be used. We will discuss this issue in Section~\ref{chap:self}.

\subsection{The Cauchy problem for smooth solutions}
\label{sub:Cauchy-LL}
Let us consider the anisotropic LL equation
\eqref{LL-ani} with $\lambda_1,\lambda_3\geq 0$. Since the associated energy is given by
\bq
\label{energy:ani}
E_{\lambda_1,\lambda_3}(\bmm)= \frac{1}{2} \int_{\R^N} ( |\nabla \bmm|^2
+\lambda_1 m_1^2+\lambda_3 m_3^2),
\eq
the natural functional setting for solving this equation is the energy set
$$\boE_{\lambda_1,\lambda_3}(\R^N) = \big\{ \bv \in L_{\rm loc}^1(\R^N, \R^3) : \abs{\bv}=1\ \textup{a.e.}, \ \nabla \bv \in L^2(\R^N),  \ \lambda_1v_1, \lambda_3 v_3 \in L^2(\R^N) \big\}.$$

In the context of functions taking values on  $\SS^2$, it is standard to use the notation 
$$\boH^\ell(\R^N)=
\big\{ \bv \in L_{\rm loc}^1(\R^N, \R^3) : 
\abs{\bv}=1\ \textup{a.e.},\ \nabla \bv \in H^{\ell-1}(\R^N)  \big\},
$$
for an integer $\ell\geq 1$, where $H^{\ell-1}$ is the classical Sobolev space. Notice that a function $\bv\in \boH^\ell(\R^N)$
does not belong to  $L^2(\R^N, \R^3)$, since this is incompatible with the constraint  $\abs{\bv}=1$. In this manner,  $\boE_{\lambda_1,\lambda_3}(\R^N)$
reduces to $\boH^1(\R^N)$ if $\lambda_1=\lambda_3=0$.

For the sake of simplicity, in this section we drop the subscripts $\lambda_1$ and $\lambda_3$, and denote the energy by $E(\bmm)$ and the space by $\boE(\R^N) $, since the constants $\lambda_1$ and $\lambda_3$ are fixed. 

The first results concerning the existence of weak solutions of \eqref{LL-ani}
in the energy space were obtained by Zhou and Guo in the one-dimensional case $N=1$ \cite{ZhouGuo1}, and by Sulem, Sulem and Bardos~\cite{SulSuBa1} for  $N\geq 1$. The approach followed in 
\cite{ZhouGuo1} was to consider a parabolic regularization by adding 
the  term $\ve \Delta \m$ and letting $\ve\to 0$ (see e.g.\ \cite{GuoDing0}), 
while the strategy in \cite{SulSuBa1} relied on finite difference approximations
and a weak compactness argument. In both cases, no uniqueness was obtained. 
The proof  in \cite{SulSuBa1} can be generalized to include the anisotropic perturbation in \eqref{LL-ani}, leading to the existence of a global 
(weak) solution as follows. 
\begin{thm}[\cite{SulSuBa1}]
	\label{thm:weak-existence}
	For any $\bmm_0 \in \boE(\R^N)$, there exists a global solution of \eqref{LL-ani}
	with $\bmm\in L^\infty(\R^+,\boE(\R^N))$, associated with the initial condition $\bmm_0.$
\end{thm}
The uniqueness  of the solution in Theorem~\ref{thm:weak-existence}
not known. To our knowledge, the well-posedness of the Landau--Lifshitz equation for general initial data in $\boE(\R^N)$ remains an open question.

Let us now discuss some results  about  smooth solutions in $\boH^k(\R^N)$, $k\in \N$, in the isotropic case 
$\lambda_1=\lambda_3=0$. For an initial data in
$\bmm_0 \in \boH^k(\R^N)$, Sulem, Sulem and Bardos~\cite{SulSuBa1}  proved the local existence and uniqueness\footnote{Actually, in \cite{SulSuBa1} they do not study of the difference between two solutions.	It is only asserted that  uniqueness followed from regularity, which it is not clear in this case; see also \cite{JerrSme1}.} of a solution
$\bmm\in L^\infty([0,T),\boH^k(\R^N))$, provided that $k>N/2+2$.
By using a parabolic approximation, Ding and Wang~\cite{DingWan1} proved the local existence in  $L^\infty([0,T),\boH^k(\R^N))$,  provided that $k>N/2+1$. 
They also study the difference between two solutions, obtaining 
uniqueness provided that  the solutions are of class $\boC^3$.
Another approach was used by McMahagan~\cite{McGahag1}, showing the existence 
as the limit of solutions of a perturbed  wave problem, and 
using parallel transport to compare two solutions,  to conclude  local existence and uniqueness in  $L^\infty([0,T),\boH^k(\R^N))$, for $k>N/2+1$. 

When $N=1$,  these results provided the local existence and uniqueness at level $\boH^k(\R^N)$, for $k\geq 2$. Moreover, in this case the solutions are global in time (see \cite{RodRubSta,ChaShUh1}).

Of course, there is a large amount of other works with interesting results about the (local and global) existence and  uniqueness for the LL equation and other related equations, see e.g.\ \cite{BeIoKeT1,GuoDing0,GustSha1,guan,JerrSme1,SongWang} and the references therein. However, it is not straightforward to adapt 
these works to obtain local well-posednes results for smooth solutions to equation \eqref{LL-ani}. For this reason,  in the rest of this section
we  provide an alternative proof for local well-posedness by
introducing high order energy quantities with better symmetrization properties. 

To study the Cauchy problem of smooth solutions, given an integer $k \geq 1$, we introduce the set 
$$\boE^k(\R^N) = \boE(\R^N) \cap \boH^{k}(\R^N),$$
which we endow with the metric structure provided by the norm
$$\| \bv \|_{Z^k} =
\big( \|\grad  \bv \|_{H^{k-1}}^2 + \| \bv_2 \|_{L^\infty}^2 +
\lambda_1 \| \bv_1 \|_{L^2}^2+\lambda_3 \| \bv_3 \|_{L^2}^2 \big)^\frac{1}{2}.
$$
Observe that the energy space $\boE(\R^N)$ identifies with $\boE^1(\R^N)$.
The uniform control on the second component $\bv_2$ in the $Z^k$-norm ensures that  $\| \cdot \|_{Z^k}$ is a norm. 
Of course, this uniform control is not the only possible choice of the metric structure.
%
The main result of this section is the following local well-posedness result.
\begin{thm}[\cite{deLaGra2}]
	\label{thm:LL-Cauchy}
	Let $\lambda_1,\lambda_3\geq 0$  and $k \in \N$, with $k > N/2 + 1$. For any initial condition $\bmm^0 \in \boE^k(\R^N)$, there exist  $T_{\max}>0$ and a unique solution $\bmm : \R^N \times[0, T_{\max}) \to \SS^2$ to the LL equation \eqref{LL-ani}, which satisfies the following statements.
	\begin{enumerate}
		\item 
		The solution $\bmm$ belongs to  $L^\infty([0, T], \boE^k(\R^N))$ and  $\partial_t \bmm \in L^\infty([0, T], \boH^{k - 2}(\R^N))$, for all $T\in (0, T_{\max})$.
		\item  If the maximal time of existence $T_{\max}$ is finite, then
		\begin{equation}
			\label{eq:cond-Tmax-LL}
			\int_0^{T_{\max}} \| \nabla \m(t) \|_{L^\infty}^2 \, dt = \infty.
		\end{equation}
		\item  The flow map $\bmm^0 \mapsto \bmm$ is well-defined and locally Lipschitz continuous from $\boE^k(\R^N)$ to $\boC^0([0, T], \boE^{k - 1}(\R^N))$, for all $T\in (0, T_{\max})$.
%
	\item  The energy \eqref{energy:ani} is conserved along the flow.
	\end{enumerate}
\end{thm}

Theorem~\ref{thm:LL-Cauchy} provides the local well-posedness of the LL equation in the set $\boE^k(\R^N)$.
This kind of statement is standard in the context of hyperbolic systems (see e.g.~\cite[Theorem 1.2]{Taylor03}). The critical regularity for the equation is given by the condition $k = N/2$, so that local well-posedness is expected when $k > N/2 + 1$. This assumption is used to control uniformly the gradient of the solutions by the Sobolev embedding theorem. 

The proof of Theorem~\ref{thm:LL-Cauchy} is based on  energy estimates using well-tailored high order energies. A key observation is that
any smooth function  $\m$ valued into  $\SS^2$, 
satisfies the pointwise identities
$$\langle \m, \partial_i \m \rangle_{\R^3} = \langle \m, \partial_{ii} \m \rangle_{\R^3} + |\partial_i \m|^2 = \langle \m, \partial_{iij} \m \rangle_{\R^3} + 2 \langle \partial_i \m, \partial_{ij} \m \rangle_{\R^3} + \langle \partial_j \m, \partial_{ii} \m \rangle_{\R^3} = 0,$$
for any $1 \leq i, j \leq N$. This allows us to show that a (smooth)
solution to \eqref{LL-ani} satisfies the equation  
\begin{equation}
	\label{eq:second-LL}
	\partial_{tt} \m + \Delta^2 \m - (\lambda_1 + \lambda_3) \big( \Delta m_1 \be_1 + \Delta m_3 \be_3 \big) + \lambda_1 \lambda_3 \big( m_1 \be_1 + m_3 \be_3 \big) = F(\m),
\end{equation}
where we have set
\begin{equation*}
	\begin{aligned}
		 F(\m) =& \sum_{1 \leq i, j \leq N} \Big( \partial_i \big( 2 \langle \partial_i \m, \partial_j \m \rangle_{\R^3} \partial_j \m - |\partial_j \m|^2 \partial_i \m \big) - 2 \partial_{ij} \big( \langle \partial_i \m, \partial_j \m \rangle_{\R^3} \m \big) \Big)\\
		& + \lambda_1  F_{1,3}^+(\m) 
 + \lambda_3  F_{3,1}^-(\m) + \lambda_1 \lambda_3 \big( (m_1^2 + m_3^2) \m + m_1^2 m_3 \be_3 + m_3^2 m_1 \be_1 \big),
	\end{aligned}
\end{equation*}
with
\begin{align*}
	F^{\pm}_{i,j}(\m)=&
	\div \big( (m_j^2 - 2 m_i^2) \nabla \m + (m_1 m_3 \be_3 \pm m_1 \m - m_3^2 \be_1 ) \nabla m_1 + (m_1 m_3 \be_1 \mp m_3 \m - m_1^2 \be_3) \nabla m_3 \big)\\
	& \pm \nabla m_1 \cdot \big( m_1 \nabla \m - \m \nabla m_1 \big) \pm  \nabla m_3 \cdot \big( \m \nabla m_3 - m_3 \nabla \m \big) + m_j |\nabla \m|^2 \be_j\\
	& + \big( m_1 \nabla m_3 - m_3 \nabla m_1 \big) \cdot \big( \nabla m_1 \be_3 - \nabla m_3 \be_1 \big) + \lambda_i m_i^2 \big( m_i \be_i - \m \big).
\end{align*}

In view of \eqref{eq:second-LL}, 
we define the (pseudo)energy of order $k \geq 2$, as
\begin{align*}
	E_{k}(t) =&
	\norm{\partial_t \m}_{\dot H^{k-2}}^2
	+\norm{ \m}_{\dot H^{k}}^2
	+ (\lambda_1 + \lambda_3)
	(\norm{ m_1}_{\dot H^{k-1}}^2+\norm{ m_3}_{\dot H^{k-1}}^2)
 + \lambda_1 \lambda_3 
	(\norm{ m_1}_{\dot H^{k-2}}^2+\norm{ m_3}_{\dot H^{k-2}}^2),
\end{align*}
for any $t \in [0, T]$. This high order energy is an  anisotropic version of the one used in~\cite{SulSuBa1}.

To get good energy estimates, we need to use Moser estimates (also called tame estimates) in Sobolev spaces (see e.g.~\cite{Moser1}). 
Using these estimates and  differentiating $E_k$, we   obtain the following energy estimates.
\begin{prop}
	\label{prop:LL-energy-estimate}
	Let $\lambda_1,\lambda_3\geq 0$ and $k \in \N$, with $k > 1 + N/2$. Assume that $\m$ is a solution to~\eqref{LL-ani} in $\boC^0([0, T], \boE^{k + 2}(\R^N))$, with $\partial_t \m \in \boC^0([0, T], H^k(\R^N))$.
	\begin{enumerate}\item\label{prop:der-energia1}	
		The LL energy is well-defined and conserved along flow on $[0, T]$.
		\item\label{prop:der-energia2}	
		Given any integer $2 \leq \ell \leq k$, the energies $E_\ell$ are of class $\boC^1$ on $[0, T]$, and there exists  $C_k>0$, depending only on $k$, such that their derivatives satisfy
		\begin{equation}
			\label{eq:energy-estimate-LL}
			 E_\ell '(t) \leq C_k \big( 1 + \| m_1(t) \|_{L^\infty}^2 + \| m_3(t) \|_{L^\infty}^2 + \| \nabla \m(t) \|_{L^\infty}^2 \big) \, \Sigma_\ell(t),
		\end{equation}
		for any $t \in [0, T]$. Here, we have set $\Sigma_\ell = \sum_{j = 1}^\ell E_j$.
	\end{enumerate}
\end{prop}

We next discretize the equation by using a finite-difference scheme. The a priori bounds remain available in this discretized setting. We then apply standard weak compactness and local strong compactness results in order to construct local weak solutions, which satisfy statement \ref{prop:der-energia1} in Theorem~\ref{thm:LL-Cauchy}. 
By applying the Gronwall lemma and the condition in~\eqref{eq:cond-Tmax-LL},  inequality \eqref{eq:energy-estimate-LL} 
prevents a possible blow-up.

Finally, we establish uniqueness, as well as continuity with respect to the initial datum, by computing energy estimates for the difference of two solutions. More precisely, we show
\begin{prop}
	\label{prop:LL-diff-control}
	Let $\lambda_1,\lambda_3\geq 0$, and $k \in \N$, with $k > N/2 + 1$. Consider two solutions $\m$ and $\bm{\tilde{m}}$ to~\eqref{LL-ani}, which lie in $\boC^0([0, T], \boE^{k + 1}(\R^N))$, with $\partial_t \m, \partial_t \bm{\tilde{m}} \in \boC^0([0, T], H^{k - 1}(\R^N))$, and set $\bm u = \bm{\tilde{m}} - \m$ and $\bv = (\bm{\tilde{m}} + \m)/2$.
	\begin{enumerate}
		\item  The function
$\gE_0(t) = \norm{ \bu(x, t) - u_2(x,0)  \be_2 }^2_{L^2}$
		is of class $\boC^1$ on $[0, T]$, and there exists $C>0$ such that	for any $t \in [0, T]$,
		\begin{align*}
				\gE_0 '(t) \leq  C \big( 1 +& \| \nabla \tilde{\m} \|_{L^2} + \| \nabla \m(t) \|_{L^2} + \| \tilde{m}_1 \|_{L^2} + \| m_1 \|_{L^2}\\
				& + \| \tilde{m}_3 \|_{L^2} + \| m_3 \|_{L^2} \big) \, \big( \| \bm u - u_2^0  e_2 \|_{L^2}^2 + \| \bm u \|_{L^\infty}^2 + \| \nabla \bm u \|_{L^2}^2 + \| \nabla u_2^0 \|_{L^2}^2 \big).
		\end{align*}
		\item   The function
$\gE_1(t) = \norm{ \nabla\bm  u}^2_{L^2} +
		\norm{\bm u \times \nabla \bm v +\bm  v \times \nabla\bm  u}^2_{L^2}$
		is of class $\boC^1$ on $[0, T]$, and there exists  $C>0$ such that
		\begin{equation*}
			\begin{split}
				 \gE_1'(t) \leq C & \big( 1 + \| \nabla \m \|_{L^\infty}^2 + \| \nabla \tilde{\m} \|_{L^\infty}^2 \big) \, \big( \|\bm  u  \|_{L^\infty}^2 + \| \nabla\bm u  \|_{L^2}^2 \big) \times\\
				& \times \big( 1 + \| \nabla \m \|_{L^\infty} + \| \nabla \tilde{\m} \|_{L^\infty} + \| \nabla \m \|_{H^1} + \| \nabla \tilde{\m} \|_{H^1} \big).
			\end{split}
		\end{equation*}
		
		\item  Let $2 \leq \ell \leq k - 1$, 
		\begin{align*}
			\gE_\ell(t)=&
			\norm{\partial_t \bm u}_{\dot H^{k-2}}^2
			+\norm{ \bm u}_{\dot H^{k}}^2
			+ (\lambda_1 + \lambda_3)
			(\norm{ u_1}_{\dot H^{k-1}}^2+\norm{ u_3}_{\dot H^{k-1}}^2)
			+ \lambda_1 \lambda_3 
			(\norm{ u_1}_{\dot H^{k-2}}^2+\norm{ u_3}_{\dot H^{k-2}}^2),
		\end{align*}
		and $\gS_{\rm LL}^\ell = \sum_{j = 0}^\ell \gE_{\rm LL}^j$. Then $\gE_\ell \in \boC^1([0, T]),$ and there exists  $C_k>0$,  such that
		\begin{equation*}
			\begin{aligned}
				  \gE_\ell'(t) \leq &C_k \Big( 1 + \| \nabla \m \|_{H^\ell}^2 + \| \nabla \tilde{\m} \|_{H^\ell}^2 + \| \nabla \m \|_{L^\infty}^2 + \| \nabla \tilde{\m} \|_{L^\infty}^2\\
				& + \delta_{\ell = 2} \big( \| \tilde{m}_1 \|_{L^2} + \| m_1 \|_{L^2} + \| \tilde{m}_3 \|_{L^2} + \| m_3 \|_{L^2} \big) \Big) \, \big( \gS_{\rm LL}^\ell + \| \bm u  \|_{L^\infty}^2 \big).
			\end{aligned}
		\end{equation*}
		\end{enumerate}
\end{prop}

When $\ell \geq 2$, the quantities $\gE_{\rm LL}^\ell$ in Proposition~\ref{prop:LL-diff-control} are anisotropic versions of the ones used in~\cite{SulSuBa1} for similar purposes. Their explicit form is related to the linear part of the second-order equation in~\eqref{eq:second-LL}. The quantity $\gE_{\rm LL}^0$ is tailored to close off the estimates.

The introduction of the quantity $\gE_{\rm LL}^1$ is of a different nature. The functions $\nabla \bm u$ and $\bm u \times \nabla \bm v + \bm v \times \nabla \bm  u$ in its definition appear as the good variables to perform hyperbolic estimates at an $H^1$-level. They provide a better symmetrization corresponding to a further cancellation of the higher order terms. Without any use of the Hasimoto transform, nor of parallel transport, this makes possible a direct proof of local well-posedness at an $H^k$-level, with $k > N/2 + 1$ instead of $k > N/2 + 2$.
\subsection{Local well-posedness for smooth solutions}
\label{sec:local-cauchy-HLL}

To state a well-posedness result for \eqref{HLL}, we need to introduce a functional setting in which we can legitimate the use of the hydrodynamical framework. 
Under the condition $\abs{\m}<1$, it is natural to work in the Hamiltonian framework in which the solutions $\m$ have finite energy. In the hydrodynamical formulation, the  energy is given by
\begin{equation}
	E_{\rm H}(u, \varphi) = \frac{1}{2} \int_{\R^N} \Big( \frac{|\nabla u|^2}{1 - u^2} + (1 - u^2) |\nabla \varphi|^2 + \lambda_1 (1 - u^2) \sin^2(\varphi) + \lambda_3 u^2 \Big).
\end{equation}
As a consequence, we work in the nonvanishing 
\begin{align*}
	\boN\boV_{\sin}^k(\R^N) &= \big\{ (u, \varphi) \in H^k(\R^N) \times H_{\sin}^k(\R^N) : |u| < 1 \ {\rm on} \ \R^N \big\},
\end{align*}
	where 
\begin{align*}
H_{\sin}^k(\R^N) &= \big\{ v \in L_{\rm loc}^1(\R^N) : \nabla v \in H^{k - 1}(\R^N) \ {\rm and} \ \sin(v) \in L^2(\R^N) \big\}.
\end{align*}
The set $H_{\sin}^k(\R^N)$ is an additive group, which is naturally endowed with the pseudometric distance
\begin{equation*}
	d_{\sin}^k(v_1, v_2) = \| \sin(v_1 - v_2) \|_{L^2}+ \| \nabla v_1 - \nabla v_2 \|_{H^{k - 1}},
\end{equation*}
that vanishes if and only if $v_1-v_2 \in \pi \Z$. This quantity is not a distance on the group $H_{\sin}^k(\R^N)$, but it is on the quotient group $H_{\sin}^k(\R^N)/\pi \Z$. In the sequel, we identify the set $H_{\sin}^1(\R^N)$ with this quotient group when necessary, in particular when a metric structure is required. This identification is not a difficulty as far as we deal with the hydrodynamical form of the LL equation and with the Sine--Gordon equation. Both the equations are indeed left invariant by adding a constant number in $\pi \Z$ to the phase functions. This property is one of the motivations for introducing the pseudometric distance $d_{\sin}^k$. We refer to \cite{deLaGra2} for more details concerning this distance, as well as the set $H_{\sin}^k(\R^N)$.

From Theorem~\ref{thm:LL-Cauchy}, we obtain the following local well-posedness result for \eqref{HLL}.
\begin{cor}[\cite{deLaGra2}]
	\label{cor:HLL-Cauchy}
	Let $\lambda_1,\lambda_3\geq 0$, and $k \in \N$, with $k > N/2 + 1$. Given any  $(u^0, \phi^0) \in \boN\boV_{\sin}^k(\R^N)$, there exist  $T_{\max}>0$ and a unique solution $(u, \phi) : \R^N \times[0, T_{\max}) \to (- 1, 1) \times \R$ to~\eqref{HLL} with initial data $(u^0, \phi^0)$, which satisfies the following statements.
	
	\begin{enumerate}
		\item  The solution $(u, \phi)$ is in $L^\infty([0, T], \boN\boV_{\sin}^k(\R^N))$, while  $(\partial_t u, \partial_t \phi)$ is in $L^\infty([0, T], H^{k - 2}(\R^N)^2)$, for any $T\in(0, T_{\max})$.
	
	\item  If the maximal time of existence $T_{\max}$ is finite, then
	$$\int_0^{T_{\max}} \Big( \Big\| \frac{\nabla u(t)}{(1 - u(t)^2)^\frac{1}{2}} \Big\|_{L^\infty}^2 + \Big\| (1 - u(t)^2)^\frac{1}{2} \nabla \phi(t) \Big\|_{L^\infty}^2 \Big) \, dt = \infty, \quad {\rm or} \quad \lim_{t \to T_{\max}} \| u(t) \|_{L^\infty} = 1.$$
	
	\item The map $(u^0, \phi^0) \mapsto (u, \phi)$ is locally Lipschitz continuous from $\boN\boV_{\sin}^k(\R^N)$ to $\boC^0([0, T], \boN\boV_{\sin}^{k - 1}(\R^N))$ for any  $T\in(0, T_{\max})$, and the energy $E_{\rm H}$ is conserved along the flow.
	\end{enumerate}

\end{cor}


The proof of Corollary~\ref{cor:HLL-Cauchy} is complicated by the metric structure corresponding to the set $H_{\sin}^k(\R^N)$. Establishing the continuity of the flow map with respect to the pseudometric distance $d_{\sin}^k$ is not so immediate, but  this difficulty can be  by-passed by using some trigonometric identities. 


\subsection{Local well-posedness  in the energy space in dimension one}
\label{sub:origin-Cauchy}
We focus now  on the  LL equation with
easy-plane anisotropy in dimension one, i.e. $\lambda_1=0$
and \eqref{LL-ani} reads
\begin{equation}
	\label{LL-1d}
	\partial_t \m + \m \times( \partial_{xx} \m - \lambda_3 m_3 \be_3) = 0.
\end{equation}
As mentioned before, in the isotropic case $\lambda_3=0$, 
 we have the local  well-posedness  for initial data in $\boH^2(\R)$ \cite{ChaShUh1,NaShVeZ1,RodRubSta}.
Theorem~\ref{thm:LL-Cauchy} gives us for instance, the
$\boH^2$-local well-posedness, while  Theorem~\ref{thm:weak-existence}
provides the existence of a solution in $\boH^1(\R)$, i.e.\ in  the energy space for the isotropic equation. The isotropic equation is energy critical in 
$\boH^{1/2}$, so that one could think that local well-posedness at the $\boH^1$-level would be simple to establish. In this direction,  
when the domain is the torus, some progress has been made 
at the $H^{3/2^+}$-level \cite{ChErTz1}, 
and an ill-posedness type result is given in \cite{JerrSme1}
for the $H^{1/2}$-weak topology.

The purpose of this section is to provide a local well-posedness 
theory for \eqref{LL-1d} in the energy space, in the case $\lambda_3\geq 0$.
To this end, we  use the hydrodynamical version of the equation,
considering  hydrodynamical variables $u = m_3$ and 
$w = -\partial_x \varphi$, that is 
\begin{equation}
	\label{HLL-1d}
	\tag{H1d}
	\left\{
	\begin{aligned}
		\partial_t u &= \partial_x \big( (u^2 - 1) w \big),\\
		\displaystyle \partial_t w &= \partial_x \Big( \frac{\partial_{xx} u}{1 - u^2} + u \frac{(\partial_x u)^2}{(1 - u^2)^2} + u \big( w^2 - \lambda_3) \Big).
	\end{aligned}
	\right.
\end{equation}
We introduce the notation $\gu = (u, w)$, that we will refer to as
hydrodynamical pair. Notice that the LL energy is now expressed as
\begin{equation*}
	E(\gu) = \int_\R e(\gu) = \frac{1}{2} \int_\R \Big( \frac{(u')^2}{1 - u^2} + \big( 1 - u^2 \big) w^2 + \lambda_3 u^2 \Big),
\end{equation*}
and the nonvanishing space is
$${\boN\boV}(\R) = \Big\{ \gv = (v, w) \in H^1(\R) \times L^2(\R), \ {\rm s.t.} \ \max_\R |v| < 1 \Big\},$$
endowed  with the metric structure corresponding to the norm
$\| \gv \|_{H^1 \times L^2} = \| v \|_{H^1} + \| w \|_{L^2}.$

Another formally conserved quantity is the momentum $P$, which is  defined by
$P(\gu) = \int_\R u w.$ As we will see in Section~\ref{chap:stability}, the momentum $P$, as well as the energy $E$, play an important role in the construction and the qualitative analysis of the solitons.

Concerning the Cauchy problem for \eqref{HLL-1d}, we have the following local well-posedness result.
\begin{thm}[\cite{deLaGra1}]
	\label{thm:local-hydro-Cauchy}
	Let $\lambda_3\geq 0$	and  $\gu^0 = (u^0, w^0) \in \boN\boV(\R)$. There exist  $T_{\max}>0$ and  $\gu = (u, w) \in \boC^0([0, T_{\max}), \boN\boV(\R))$, such that 
	the following statements hold.
	\begin{enumerate}
		\item 
		The map $\gu$ is the unique solution to \eqref{HLL-1d}, with initial condition $\gu^0$, such that there exist smooth solutions $\gu_n \in \boC^\infty(\R \times [0, T])$ to \eqref{HLL-1d}, which satisfy
			$\gu_n \to \gu$ in $\boC^0([0,T], \boN\boV(\R)),$
		as $n \to  \infty$, for any $T \in (0, T_{\max})$. In addition, 
the energy $E$ and the momentum $P$ are constant on $(0, T_{\max})$.
		\item The maximal time $T_{\max}$ is characterized by the condition
		$$\lim_{t \to T_{\max}} \, \max_{x \in \R} |u(x, t)| = 1, \quad \text { if } T_{\max} <  \infty.$$
		\item When	$\gu_n^0 \to \gu^0$ in $H^1(\R) \times L^2(\R),$
		as $n \to  \infty$, the maximal time of existence $T_n$ of the solution $\gu_n$ to \eqref{HLL-1d}, with initial condition $\gu_n^0$, satisfies
$			T_{\max} \leq \liminf_{n \to  \infty} \, T_n,$
		and
		$\gu_n \to \gu$ in $\boC^0([0, T], H^1(\R) \times L^2(\R)),$
		as $n \to  \infty$, for any $T \in (0, T_{\max})$.
\end{enumerate}
\end{thm}

In other words, Theorem~\ref{thm:local-hydro-Cauchy} provides the existence and uniqueness of a continuous flow for \eqref{HLL-1d} in the energy space $\boN\boV(\R)$.	On the other hand, this does not prevent from the existence of other solutions which could not be approached by smooth solutions. In particular, we do not claim that there exists a unique local solution to \eqref{HLL-1d} in the energy space for a given initial condition. To our knowledge, the question of the global existence in the hydrodynamical framework of the local solution $\gv$ remains open.
Concerning  the   equation  \eqref{LL-1d}, since we 
are in the one-dimensional case, it is possible to endow the energy space 
with the metric structure corresponding to the distance
$$d_\boE(\bu, \bv) = \Big( \big| \check{\bu}(0) - \check{\bv}(0) \big|^2 + \big\| \bu' -\bv' \big\|_{L^2}^2 +\lambda_3 \big\| u_3 - v_3 \big\|_{L^2}^2 \Big)^\frac{1}{2},$$
and to translate Theorem~\ref{thm:local-hydro-Cauchy} into the original framework of the LL equation. This provides the existence of a unique continuous flow for \eqref{LL-1d} in the neighborhood of solutions $\m$, such that the third component $\m_3$ does not reach the value $\pm 1$. The flow is only locally defined due to this restriction. 

The most difficult part 
in Theorem~\ref{thm:local-hydro-Cauchy} is
the continuity with respect to the initial data in the energy space $\boN\boV(\R)$ when 
$\lambda_3>0$. In this case, by performing a change of variables, we can assume that 
$\lambda_3=1$.The proof relies on the strategy developed by Chang, Shatah and Uhlenbeck in~\cite{ChaShUh1} (see also \cite{GustSha1, NaShVeZ1}), by introducing the map
\begin{equation}
	\label{def:Psi}
	\Psi = \frac{1}{2} \Big( \frac{\partial_x u}{(1 - u^2)^\frac{1}{2}} + i (1 - u^2)^\frac{1}{2} w \Big) \exp i \theta,\quad   \text{ with }
	\theta(x, t) = - \int_{- \infty}^x u(y, t) w(y, t) \, dy.
\end{equation}
Then $\Psi$ solves the nonlinear Schr\"odinger equation
\begin{equation}
	\label{eq:Psi}
	i \partial_t \Psi + \partial_{xx} \Psi + 2 |\Psi|^2 \Psi + \frac{1}{2} u^2 \Psi - \Re \Big( \Psi \big( 1 - 2 F(u, \overline{\Psi}) \big) \Big) \big( 1 - 2 F(u, \Psi) \big) = 0,
\end{equation}
with
$	F(u, \Psi)(x, t) = \int_{- \infty}^x u(y, t) \Psi(y, t) \, dy,$
while the function $u$ satisfies
\begin{equation}
		\partial_t u = 2 \partial_x \Im \Big( \Psi \big( 2 F(u, \overline{\Psi}) - 1 \big) \Big),\quad 
		\partial_x u = 2 \Re \Big( \Psi \big( 1 - 2 F(u, \overline{\Psi}) \big) \Big).
\end{equation}
In this setting, deriving the continuous dependence in $\boN\boV(\R)$ of $\gu$ with respect to its initial data reduces to establish it for $u$ and $\Psi$ in $L^2(\R)$. This can be  done  by 
combining an energy method for $u$ and classical Strichartz estimates for $\Psi$.



\section{Asymptotics regimes}
\label{sec:regime}

In this section we  will study the connection between the LL equation
\begin{equation}
	\label{LL-ani:asypm}
	\partial_t \bmm+  \bmm \times (\Delta \bmm-\lambda_1m_1\be_1- \lambda_3m_3 \be_3)=0, 
\end{equation}
with $\lambda_1,\lambda_3\geq 0$,
and the Sine--Gordon and  the NLS equations, for certain types of anisotropies.
More precisely, 
we investigate the cases when $\lambda_1\ll \lambda_3$ and when $1\ll \lambda_1= \lambda_3$. A conjecture in the physical literature  \cite{Sklyani1,FaddTak0} is that  in the former case, the dynamics of  \eqref{LL-ani:asypm} can be described by the Sine--Gordon equation, while in the latter case, can be approximated by the cubic NLS equation.

It is well-known that deriving asymptotic regimes is a powerful tool in order to tackle the analysis of intricate equations. In this direction, we expect that these rigorous derivations will be a useful tool  to describe the dynamical properties of the LL equation, in particular the role played by the solitons in this dynamics. For instance, this kind of strategy has been
useful in order to prove the asymptotic stability of the dark solitons of the Gross-Pitaevskii equation by using its link with the KdV equation
(see \cite{ChirRou2,BeGrSaS3}).

\subsection{The Sine--Gordon regime}
In order to provide a rigorous mathematical statement
for the anisotropic LL equation with  $\lambda_1\ll \lambda_3$,
i.e.\ for a strong easy-plane anisotropy regime, 
we consider a small  parameter $\varepsilon>0$, a fixed constant  $\sigma>0$, and set the anisotropy values
$	\lambda_1 = \sigma \varepsilon$ and $ \lambda_3 = {1}/{\varepsilon}$.

Assuming  that the map $\check{m} = m_1 + i m_2$, associated with a solution $\bmm$ to \eqref{LL-ani:asypm} does not vanish, we write
$\check{m} = (1 - m_3^2)^\frac{1}{2} \big( \sin(\phi) + i\cos(\phi) \big)$,
so that the variables $u = m_3$ and $\phi$
satisfy the system \eqref{HLL}, as long as the nonvanishing condition  holds. To study the behavior of the system as $\ve\to0$, we introduce the rescaled variables $U_\varepsilon$ and $\Phi_\varepsilon$ given by 
$$U_\varepsilon ( x, t)=\frac{ u(x/\sqrt{\ve},t)}{\ve},
\quad {\rm and} \quad \Phi_\varepsilon( x, t)= {\phi(x/\sqrt{\ve},t)},
$$
which satisfy the hydrodynamical system 
\begin{equation}
	\hspace*{-0.1cm}
	\tag{${\rm H}_\varepsilon$}
	\label{HLLeps}
	\left\{
	\begin{aligned}
		\partial_t U_\varepsilon &= \div \big( (1 - \varepsilon^2 U_\varepsilon^2) \nabla \Phi_\varepsilon \big) - \frac{\sigma}{2} (1 - \varepsilon^2 U_\varepsilon^2) \sin(2 \Phi_\varepsilon),\\
		\partial_t \Phi_\varepsilon &= U_\varepsilon \big( 1 - \varepsilon^2 \sigma \sin^2(\Phi_\varepsilon) \big) - \varepsilon^2 \div \Big( \frac{\nabla U_\varepsilon}{1 - \varepsilon^2 U_\varepsilon^2} \Big) +\varepsilon^4 U_\varepsilon \frac{|\nabla U_\varepsilon|^2}{(1 - \varepsilon^2 U_\varepsilon^2)^2} - \varepsilon^2 U_\varepsilon |\nabla \Phi_\varepsilon|^2. 
	\end{aligned}
	\right.
\end{equation}
Therefore, as $\varepsilon \to 0$, we formally see that  the limit system is 
\begin{equation}
	\label{sys:SG}
	\tag{SGS}
\partial_t U = \Delta \Phi - \frac{\sigma}{2} \sin(2 \Phi),\quad 
		\partial_t \Phi = U,
\end{equation}
so that  the limit function $\Phi$ is a solution to the Sine--Gordon equation
\begin{equation}
	\tag{SG}
	\label{SG}
	\partial_{tt} \Phi - \Delta \Phi + \frac{\sigma}{2} \sin(2 \Phi) = 0.
\end{equation}

As seen in Corollary~\ref{cor:HLL-Cauchy},  the hydrodynamical system \eqref{HLLeps}
is locally well-posed in the space $\boN\boV_{\sin}^k(\R^N)$ for $k>N/2+1$.  
However, this result gives us time of existence $T_\ve$ that could 
vanish as $\ve\to 0$. Therefore, we need to find a uniform estimate 
for $T_\ve$  to prevent this phenomenon. As we will recall later,  the Sine--Gordon equation  is also
locally well-posed at the same level of regularity, 
so that we can compare the evolution of the difference in an interval of time independent of $\ve$. A further analysis of 
\eqref{HLLeps} involving good energy estimates, will lead us 
to the following result.

\begin{thm}[\cite{deLaGra2}]
	\label{thm:conv-SG}
	Let $N \geq 1$ and $k \in \N$, with $k > N/2 + 1$, and $ \varepsilon\in (0,1)$. Consider an initial condition $(U_\varepsilon^0, \Phi_\varepsilon^0) \in \boN\boV_{\sin}^{k + 2}(\R^N)$, and set
\bq\label{def:boK} \boK_\varepsilon = \big\| U_\varepsilon^0 \big\|_{H^k} + \varepsilon \big\| \nabla U_\varepsilon^0 \big\|_{H^k} + \big\| \nabla \Phi_\varepsilon^0 \big\|_{H^k} + \big\| \sin(\Phi_\varepsilon^0) \big\|_{H^k}.
\eq
	Consider similarly an initial condition $(U^0, \Phi^0) \in L^2(\R^N) \times H_{\sin}^1(\R^N)$, and denote by $(U, \Phi) \in \boC^0(\R, L^2(\R^N) \times H_{\sin}^1(\R^N))$ the unique corresponding solution to~\eqref{sys:SG}. Then, there exists $C>0$, depending only on $\sigma$, $k$ and $N$, such that, if 
	\begin{equation}
		\label{cond:key}
		C \, \varepsilon \, \boK_\varepsilon \leq 1,
	\end{equation}
	then the following statements hold.
\begin{enumerate}
\item There exists a positive number
$		T_\varepsilon \geq (C \boK_\varepsilon^2)^{-1},$
	such that there is a unique solution $(U_\varepsilon, \Phi_\varepsilon) \in \boC^0([0, T_\varepsilon], \boN\boV_{\sin}^{k + 1}(\R^N))$ to~\eqref{HLLeps} with initial data $(U_\varepsilon^0, \Phi_\varepsilon^0)$.
	
\item 	 If $\Phi_\varepsilon^0 - \Phi^0 \in L^2(\R^N)$, then, for any $0 \leq t \leq T_\varepsilon$, 
	\begin{equation}
		\label{est:0}
			\big\| \Phi_\varepsilon(t) - \Phi(t) \big\|_{L^2} \leq C \, \big( \big\| \Phi_\varepsilon^0 - \Phi^0 \big\|_{L^2} + \big\| U_\varepsilon^0 - U^0 \big\|_{L^2} + \varepsilon^2 \, \boK_\varepsilon \, \big( 1 + \boK_\varepsilon^3 \big) \big) \, e^{C t}.
	\end{equation}
	\item  If $N \geq 2$, or $N = 1$ and $k > N/2 + 2$, then we have, for any $0 \leq t \leq T_\varepsilon$, 
	\begin{equation}
		\big\| U_\varepsilon(t) - U( t) \big\|_{L^2} +
		d_{\sin}^1(\Phi_\varepsilon( t), \Phi( t) ) 
		\leq C \big(  \big\| U_\varepsilon^0 - U^0 \big\|_{L^2} + 
		d_{\sin}^1(
		\Phi_\varepsilon^0,\Phi^0 )
		+ \varepsilon^2 \, \boK_\varepsilon \, \big( 1 + \boK_\varepsilon^3 \big) \big)  e^{C t}.
	\end{equation}
	\item  Let $(U^0, \Phi^0) \in H^k(\R^N) \times H_{\sin}^{k + 1}(\R^N)$ and set
	$\kappa_\varepsilon = \boK_\varepsilon + \big\| U^0 \big\|_{H^k} + \big\| \nabla \Phi^0 \big\|_{H^k} + \big\| \sin(\Phi^0) \big\|_{H^k}.$
	There exists  $A>0$, depending only on $\sigma$, $k$ and $N$, such that the solution $(U, \Phi)$ lies in $\boC^0([0, T_\varepsilon^*], H^k(\R^N) \times H_{\sin}^{k + 1}(\R^N))$, for some $T_\varepsilon^*\in [\frac{1}{A \kappa_\varepsilon^2}, T_\ve]$.
	Moreover, when $k > N/2 + 3$, we have, 	for any $0 \leq t \leq T_\varepsilon^*$, 
	\begin{align*}
&		\big\|  U_\varepsilon(t) - U(t) \big\|_{H^{k - 3}} + \big\| \nabla \Phi_\varepsilon(t) - \nabla \Phi(t) \big\|_{H^{k - 3}} + \big\| \sin(\Phi_\varepsilon(t) - \Phi(t)) \big\|_{H^{k - 3}}\\
&		\leq A \,  e^{A (1 + \kappa_\varepsilon^2) t}
		\big(  \big\| U_\varepsilon^0 - U^0 \big\|_{H^{k - 3}} + \big\| \nabla \Phi_\varepsilon^0 - \nabla \Phi^0 \big\|_{H^{k - 3}} + \big\| \sin(\Phi_\varepsilon^0 - \Phi^0) \big\|_{H^{k - 3}} + 
		\varepsilon^2 \kappa_\varepsilon \big( 1 + \kappa_\varepsilon^3 \big) \big).
	\end{align*}
\end{enumerate}
\end{thm}

In arbitrary dimension, Theorem~\ref{thm:conv-SG} provides a quantified convergence of the LL equation towards the Sine--Gordon equation in the regime of strong easy-plane anisotropy. Three types of convergence are proved depending on the dimension, and the levels of regularity of the solutions. 
This trichotomy is related to the analysis of the Cauchy problems for the LL and Sine--Gordon equations.

In its natural Hamiltonian framework, the Sine--Gordon equation is globally well-posed and 
its Hamiltonian is the Sine--Gordon energy:
\begin{equation}
	\label{def:E-SG}
	E_{\rm SG}(\phi) = \frac{1}{2} \int_{\R^N} \big( (\partial_t \phi)^2 + |\nabla \phi|^2 + \sigma \sin(\phi)^2 \big).
\end{equation}
More precisely, given an initial condition $(\Phi^0, \Phi^1) \in H_{\sin}^1(\R^N) \times L^2(\R^N)$, there exists a unique corresponding solution $\Phi \in \boC^0(\R, H_{\sin}^1(\R^N))$ to~\eqref{SG}, with $\partial_t \Phi \in \boC^0(\R, L^2(\R^N))$. Moreover, the Sine--Gordon equation is locally well-posed in the spaces $H_{\sin}^k(\R^N) \times H^{k - 1}(\R^N)$, when $k > N/2 + 1$. In other words, the solution $\Phi$ remains in $\boC^0([0, T], H_{\sin}^k(\R^N))$, with $\partial_t \Phi \in \boC^0([0, T], H^{k - 1}(\R^N))$, at least locally in time, when $(\Phi^0, \Phi^1) \in H_{\sin}^k(\R^N) \times H^{k - 1}(\R^N)$. 
We refer to \cite{deLaGra2,BuckMil1} for more details about  the  Cauchy problem for \eqref{SG}. 

As seen in Section~\ref{chap:cauchy}, the LL equation is locally well-posed at the same level of high regularity as the Sine--Gordon equation. In the hydrodynamical context, this reads as the existence of a maximal time $T_{\max}$ and a unique solution $(U, \Phi) \in \boC^0([0, T_{\max}), \boN\boV_{\sin}^{k - 1}(\R^N))$ to~\eqref{HLLeps} corresponding to an initial condition $(U^0, \Phi^0) \in \boN\boV_{\sin}^k(\R^N)$, when $k > N/2 + 1$ (see Corollary~\ref{cor:HLL-Cauchy}); note the loss of one derivative here. This loss explains why we take initial conditions $(U_\varepsilon^0, \Phi_\varepsilon^0)$ in $\boN\boV_{\sin}^{k + 2}(\R^N)$, though the quantity $\boK_\varepsilon$ is already well-defined when $(U_\varepsilon^0, \Phi_\varepsilon^0) \in \boN\boV_{\sin}^{k + 1}(\R^N)$.

In view of this local well-posedness result, we restrict our analysis of the Sine--Gordon regime to the solutions $(U_\varepsilon, \Phi_\varepsilon)$ to the rescaled system~\eqref{HLLeps} with sufficient regularity. A further difficulty then lies in the fact that their maximal times of existence possibly depend on  $\varepsilon$.

Statement (i) in Theorem~\ref{thm:conv-SG} provides an explicit control on these maximal times.
Since 	$T_\varepsilon \geq (C \boK_\varepsilon^2)^{-1},$ these maximal times are bounded from below by a positive number depending only on the choice of the initial data $(U_\varepsilon^0, \Phi_\varepsilon^0)$. Notice in particular that if a family of initial data $(U_\varepsilon^0, \Phi_\varepsilon^0)$ converges towards a pair $(U^0, \Phi^0)$ in $H^k(\R^N) \times H_{\sin}^k(\R^N)$, as $\varepsilon \to 0$, then it is possible to find   $T>0$ such that all the corresponding solutions $(U_\varepsilon, \Phi_\varepsilon)$ are well-defined on $[0, T]$. This property is necessary in order to make possible a consistent analysis of the limit $\varepsilon \to 0$.

Statement (i) only holds when the initial data $(U_\varepsilon^0, \Phi_\varepsilon^0)$ satisfy the condition in~\eqref{cond:key}. However, this condition is not a restriction in the limit $\varepsilon \to 0$. It is satisfied by any fixed pair $(U^0, \Phi^0) \in \boN\boV_{\sin}^{k + 1}(\R^N)$ provided that $\varepsilon$ is small enough, so that it is also satisfied by a family of initial data $(U_\varepsilon^0, \Phi_\varepsilon^0)$, which converges towards a pair $(U^0, \Phi^0)$ in $H^k(\R^N) \times H_{\sin}^k(\R^N)$ as $\varepsilon \to 0$.

Statements (ii) and (iii) in Theorem~\ref{thm:conv-SG} provide two estimates between the previous solutions $(U_\varepsilon, \Phi_\varepsilon)$ to~\eqref{HLLeps}, and an arbitrary global solution $(U, \Phi)$ to~\eqref{sys:SG} at the Hamiltonian level. The first one yields an $L^2$-control on the difference $\Phi_\varepsilon - \Phi$, while the second one, an energetic control on the difference $(U_\varepsilon, \Phi_\varepsilon) - (U, \Phi)$. Due to the fact that the difference $\Phi_\varepsilon - \Phi$ is not necessarily in $L^2(\R^N)$, statement (ii) is restricted to initial conditions satisfying this property.

Finally, statement (iv) bounds the difference between the solutions $(U_\varepsilon, \Phi_\varepsilon)$ and $(U, \Phi)$ at the same initial Sobolev level. In this case, we also have to control the maximal time of regularity of the solutions $(U, \Phi)$. This follows from the control from below for $T_\ve^*$, which is of the same order as the one in $T_\ve$.

We then obtain the Sobolev estimate in (iv) of the difference $(U_\varepsilon, \Phi_\varepsilon) - (U, \Phi)$ with a loss of three derivatives. Here, the choice of the Sobolev exponents $k > N/2 + 3$ is tailored 
to gain a uniform control on the functions $U_\varepsilon - U$, $\nabla \Phi_\varepsilon - \nabla \Phi$ and $\sin(\Phi_\varepsilon - \Phi)$, by the Sobolev embedding theorem.

 A loss of derivatives is natural in the context of long-wave regimes; it is related to the terms with first and second-order derivatives in the right-hand side of~\eqref{HLLeps}. This loss is the reason why the energetic estimate in statement (iii) requires an extra derivative in dimension one, that is the condition $k > N/2 + 2$. Using the Sobolev bounds in ~\eqref{borne-W-bis}, we can (partly) recover this loss by a standard interpolation argument, and deduce an estimate in $H^\ell(\R^N) \times H_{\sin}^{\ell + 1}(\R^N)$ for any number $\ell < k$. In this case, the error terms are no more of order $\varepsilon^2$.

As a by-product of the analysis, we can also analyze the wave regime for the LL equation. This regime is obtained by allowing  the parameter $\sigma$ to converge to $0$. Indeed, at least formally, a solution $(U_{\varepsilon, \sigma}, \Phi_{\varepsilon, \sigma})$ to~\eqref{HLLeps} satisfies the free wave system
\begin{equation}
	\tag{FW}
	\label{FW}
		\partial_t U = \Delta \Phi,\qquad 		\partial_t \Phi = U,
\end{equation}
as $\varepsilon \to 0$ and $\sigma \to 0$. In particular, the function $\Phi$ is solution to the wave equation
$\partial_{tt} \Phi - \Delta \Phi = 0.$
The following result provides a rigorous justification for this asymptotic approximation.

\begin{thm}[\cite{deLaGra2}]
	\label{thm:conv-wave}
	Let $N \geq 1$ and $k \in \N$, with $k > N/2 + 1$, and $0 < \varepsilon, \sigma < 1$. Consider an initial condition $(U_{\varepsilon, \sigma}^0, \Phi_{\varepsilon, \sigma}^0) \in \boN\boV_{\sin}^{k + 2}(\R^N)$
	and set
\bqq
\boK_{\varepsilon, \sigma} = \big\| U_{\varepsilon, \sigma}^0 \big\|_{H^k} + \varepsilon \big\| \nabla U_{\varepsilon, \sigma}^0 \big\|_{H^k} + \big\| \nabla \Phi_{\varepsilon, \sigma}^0 \big\|_{H^k} + \sigma^\frac{1}{2} \big\| \sin(\Phi_{\varepsilon, \sigma}^0) \big\|_{L^2}.
\eqq
	Let $m \in \N$, with $0 \leq m \leq k - 2$. Consider similarly an initial condition $(U^0, \Phi^0) \in H^m(\R^N) \times H^{m - 1}(\R^N)$, and denote by $(U, \Phi) \in \boC^0(\R, H^{m - 1}(\R^N) \times H^m(\R^N))$ the unique corresponding solution to~\eqref{FW}. Then, there exists  $C>0$, depending only on $k$ and $N$, such that, if the initial data satisfies the condition
$		C \, \varepsilon \, \boK_{\varepsilon, \sigma}^0 \leq 1,$
	the following statements hold. Then there exists a positive number
	\begin{equation}
		\label{cond:Teps-sigma}
		T_{\varepsilon, \sigma} \geq \frac{1}{C \max ( \varepsilon, \sigma ) (1 + \boK_{\varepsilon, \sigma})^{\max ( 2, k/2 ) }},
	\end{equation}
	such that there is a unique solution $(U_{\varepsilon, \sigma}, \Phi_{\varepsilon, \sigma}) \in \boC^0([0, T_{\varepsilon, \sigma}], \boN\boV_{\sin}^{k + 1}(\R^N))$ to~\eqref{HLLeps} with initial data $(U_{\varepsilon, \sigma}^0, \Phi_{\varepsilon, \sigma}^0)$. Moreover,  if $\Phi_{\varepsilon, \sigma}^0 - \Phi^0 \in H^m(\R^N)$, then we have the estimate,
		for any $0 \leq t \leq T_{\varepsilon, \sigma}$, 
	\begin{align*}
			\big\| U_{\varepsilon, \sigma}(t)   - U(t) \big\|_{H^{m - 1}} &+ \big\| \Phi_{\varepsilon, \sigma}(t) - \Phi(t) \big\|_{H^m} \leq C \big( 1 + t^2 \big) \, \Big( \big\| U_{\varepsilon, \sigma}^0 - U^0 \big\|_{H^{m - 1}}\\
			& + \big\| \Phi_{\varepsilon, \sigma}^0 - \Phi^0 \big\|_{H^m} + \max \big( \varepsilon^2,\sigma^{1/2} \big) \, \boK_{\varepsilon, \sigma}^0 \, \big( 1 + \boK_{\varepsilon, \sigma}^0 \big)^{\max ( 2, m )} \Big).
	\end{align*}
\end{thm}

The wave regime of the LL equation was first derived rigorously by Shatah and Zeng~\cite{ShatZen1}, as a special case of the wave regimes for the Schr\"odinger map equations with values into arbitrary K\"ahler manifolds.
The derivation in~\cite{ShatZen1} relies on energy estimates, which are similar in spirit to
the ones we establish in the sequel, and a compactness argument. Getting rid of this
compactness argument provides the quantified version of the convergence in
Theorem~\ref{thm:conv-wave}. 
This improvement is based on the arguments developed by B\'ethuel, Danchin and Smets~\cite{BetDaSm1} 
in order to quantify the convergence of the Gross--Pitaevskii equation towards the free wave equation in a similar long-wave regime.
Similar arguments were also applied in~\cite{Chiron9} in order to derive rigorously the (modified) KdV and (modified) KP regimes of the LL equation (see also~\cite{GermRou1}).



%

Concerning the proof of Theorem~\ref{thm:conv-SG}, the first step is to provide a control on $T_{\max}$. In view of the conditions in statement $(ii)$ of Corollary~\ref{cor:HLL-Cauchy}, this control can be derived from uniform bounds on the functions $U_\varepsilon$, $\nabla U_\varepsilon$ and $\nabla \Phi_\varepsilon$. Taking into account the Sobolev embedding theorem and the fact that $k > N/2 + 1$, we are left with the computations of energy estimates for the functions $U_\varepsilon$ and $\Phi_\varepsilon$   in the spaces $H^k(\R^N)$ and  $H_{\sin}^k(\R^N)$, respectively. 

In this direction, we recall that the LL energy corresponding to the scaled hydrodynamical system~\eqref{HLLeps} writes as
$$E_\varepsilon(t) = \frac{1}{2} \int_{\R^N} \Big( \varepsilon^2 \frac{|\nabla U_\varepsilon|^2}{1 - \varepsilon^2 U_\varepsilon^2} + U_\varepsilon^2 + (1 - \varepsilon^2 U_\varepsilon^2) |\nabla \Phi_\varepsilon|^2 + \sigma (1 - \varepsilon^2 U_\varepsilon^2) \sin^2(\Phi_\varepsilon) \Big).$$
Inspired in this formula, we proposed to define an energy of order $k\geq 1$ as
\begin{equation*}
		E_\varepsilon^k(t) =
		\frac{1}{2} \sum_{|\alpha| = k - 1} \int_{\R^N} \Big(
 \varepsilon^2 \frac{|\nabla \partial_x^\alpha U_\varepsilon|^2}{1 - \varepsilon^2 U_\varepsilon^2} + |\partial_x^\alpha U_\varepsilon|^2 + (1 - \varepsilon^2 U_\varepsilon^2) |\nabla \partial_x^\alpha \Phi_\varepsilon|^2
 + \sigma (1 - \varepsilon^2 U_\varepsilon^2) |\partial_x^\alpha \sin(\Phi_\varepsilon)|^2 \Big).
\end{equation*}
The factors $1 - \varepsilon^2 U_\varepsilon^2$ in this expression, as well as the non-quadratic term corresponding to the function $\sin(\Phi_\varepsilon)$, are of substantial importance since  they provide a better symmetrization of the energy estimates,
by inducing cancellations in the higher order terms.  More precisely, we have the following key proposition.
\begin{prop}
	\label{prop:estimate}
	Let $\varepsilon>0$  and $k \in \N$, with $k>N/2 + 1$. Consider a solution $(U_\varepsilon,\Phi_\varepsilon)$ to~\eqref{HLLeps}, with $(U_\varepsilon, \Phi_\varepsilon) \in \boC^0([0, T], \boN\boV_{\sin}^{k + 3}(\R^N))$ for some  $T>0$.
	Assume that 
	\begin{equation}
		\label{borne-W}
		\inf_{\R^N \times [0, T]} (1 - \varepsilon^2 U_\varepsilon^2 )\geq 1/2. 
	\end{equation}
	Then there exists $C>0$, depending only on $k$ and $N$, such that
	\begin{equation}
		\label{der-E-j}
		\begin{aligned}
			\big[ E_\varepsilon^\ell & \big]'(t) \leq C \, \max \big( 1, \sigma^{3/2} \big) \, \big( 1  + \varepsilon^4 \big) \, \Big( \| \sin(\Phi_\varepsilon(t)) \|_{L^\infty}^2 + \| U_\varepsilon(t) \|_{L^\infty}^2 + \| \nabla \Phi_\varepsilon(t) \|_{L^\infty}^2 + \| \nabla U_\varepsilon(t) \|_{L^\infty}^2\\
			&  + \| d^2 \Phi_\varepsilon(t) \|_{L^\infty}^2 + \varepsilon^2 \| d^2 U_\varepsilon(t) \|_{L^\infty}^2
			 + \varepsilon \, \| \nabla \Phi_\varepsilon(t) \|_{L^\infty} \, \big( \| \nabla \Phi_\varepsilon(t) \|_{L^\infty}^2 + \| \nabla U_\varepsilon(t) \|_{L^\infty}^2 \big) \Big) \, \Sigma_\varepsilon^{k + 1}(t),
		\end{aligned}
	\end{equation}
	for any $t \in [0, T]$ and any $2 \leq \ell \leq k + 1$. Here, we have set $\Sigma_\varepsilon^{k + 1} = \sum_{j = 1}^{k + 1} E_\varepsilon^j$.
\end{prop}
Thanks to the condition  $k>N/2 + 1$ and the Sobolev embedding, 
we get from \eqref{der-E-j} a differential inequality for $y(t)=\Sigma_\varepsilon^{k}$, of the type
\bq 
\label{toy:diff}
y'(t)\leq A y^2(t),
\eq
at least on the interval where $y$ is well-defined and  $y(t)\leq 2y(0)$.
Here  $A$ is a constant depending on $y(0)$. Integrating \eqref{toy:diff},
we conclude that 
$$y(t)\leq \frac{y(0)}{1-Ay(0)t}\leq 2y(0),$$
provided that $t\leq 1/(2Ay(0))$. Using this argument, we deduce from Proposition~\ref{prop:estimate}, that maximal time $T_{\max}$ is at least of order $1/(\| U_\varepsilon^0 \|_{H^k} + \varepsilon \| \nabla U_\varepsilon^0 \|_{H^k} + \| \nabla \Phi_\varepsilon^0 \|_{H^k} + \| \sin(\Phi_\varepsilon^0) \|_{H^k})^2$, when the initial conditions $(U_\varepsilon^0, \Phi_\varepsilon^0)$ satisfy the inequality in~\eqref{cond:key}. 
In particular, the dependence of $T_{\max}$ on the small parameter $\varepsilon$ only results from the possible dependence of the pair $(U_\varepsilon^0, \Phi_\varepsilon^0)$ on $\varepsilon$.
Choosing suitably these initial conditions, we can assume without loss of generality, that $T_{\max}$ is uniformly bounded from below when $\varepsilon$ tends to $0$, so that analyzing this limit makes sense and we can work in an interval of the form $[0,T_\ve]$. Moreover, we also get 
the energy estimate on $[0,T_\ve]$ in terms og $\boK_{\ve}$ defined in \eqref{def:boK},
	\begin{equation}
		\label{borne-W-bis}
			\big\| U_\varepsilon(t) \big\|_{H^k} + \varepsilon \big\| \nabla U_\varepsilon(t) \big\|_{H^k} +  \big\| \nabla \Phi_\varepsilon(t) \big\|_{H^k} + \big\| \sin(\Phi_\varepsilon(t)) \big\|_{H^k}\leq 
		C \boK_{\ve}.
	\end{equation}

The final  ingredient in the proof of Theorem~\ref{thm:conv-SG} is  the consistency of~\eqref{HLLeps} with the Sine--Gordon system in the limit $\varepsilon \to 0$. Indeed, we can rewrite~\eqref{HLLeps} as
\begin{equation}
	\label{HLLeps-ter}
 \partial_t U_\varepsilon = \Delta \Phi_\varepsilon - \frac{\sigma}{2} \sin(2 \Phi_\varepsilon) + \varepsilon^2 R_\varepsilon^U,\qquad 
		\partial_t \Phi_\varepsilon = U_\varepsilon + \varepsilon^2 R_\varepsilon^\Phi, 
\end{equation}
where 
\begin{align*}
	R_\varepsilon^U &= - \div \big( U_\varepsilon^2 \, \nabla \Phi_\varepsilon \big) + \sigma U_\varepsilon^2 \, \sin(\Phi_\varepsilon)  \cos(\Phi_\varepsilon),\\
		R_\varepsilon^\Phi &= - \sigma U_\varepsilon \, \sin^2(\Phi_\varepsilon) - \div \Big( \frac{\nabla U_\varepsilon}{1 - \varepsilon^2 U_\varepsilon^2} \Big) + \varepsilon^2 U_\varepsilon \, \frac{|\nabla U_\varepsilon|^2}{(1 - \varepsilon^2 U_\varepsilon^2)^2} - U_\varepsilon \, |\nabla \Phi_\varepsilon|^2.
\end{align*}
In view of the Sobolev control in~\eqref{borne-W-bis}, the remainder terms $R_\varepsilon^U$ and $R_\varepsilon^\Phi$ are bounded uniformly with respect to $\varepsilon$ in Sobolev spaces, with a loss of three derivatives. Due to this observation, the differences $u_\varepsilon = U_\varepsilon - U$ and $\varphi_\varepsilon = \Phi_\varepsilon - \Phi$ between a solution $(U_\varepsilon, \Phi_\varepsilon)$ to~\eqref{HLLeps} and a solution $(U, \Phi)$ to~\eqref{sys:SG} are expected to be of order $\varepsilon^2$, if the corresponding initial conditions are close enough.
The proof of this claim would be immediate if the system~\eqref{HLLeps-ter} would not contain the nonlinear term $\sin(2 \Phi_\varepsilon)$. Due to this extra term, we have to apply a Gronwall argument in order to control the differences $u_\varepsilon$ and $\varphi_\varepsilon$. This can be done since 
$v_\varepsilon$ and $\varphi_\varepsilon$ satisfy
\begin{equation*}
 \partial_t v_\varepsilon = \Delta \varphi_\varepsilon - \sigma \sin(\varphi_\varepsilon)  \cos(\Phi_\varepsilon + \Phi) + \varepsilon^2 R_\varepsilon^U,\qquad 
		\partial_t \varphi_\varepsilon = v_\varepsilon + \varepsilon^2 R_\varepsilon^\Phi,
\end{equation*}
so that we can  perform  energy estimates as before.
\subsection{The cubic NLS regime}
We now focus on the cubic Schr\"odinger equation, which is obtained in a regime of strong easy-axis anisotropy of equation \eqref{LL-ani:asypm}. For this purpose, we consider 
a uniaxial material in the direction corresponding to the vector $\be_2$ and we fix the anisotropy parameters as
$\lambda_1 = \lambda_3 = {1}/{\varepsilon}.$
For this choice, 
let us introduce the complex-valued function $\Psi_\varepsilon$ given by
\begin{equation}
	\label{def:Psi-eps}
	\Psi_\varepsilon(x, t)=\varepsilon^{- {1/2}} \check{m}(x, t) e^{i t/\varepsilon}, \quad \text{ with }
	 \check{m} = m_1 + i m_3,
\end{equation}
 associated with  a solution $\m$ of \eqref{LL-ani:asypm}.
This function is of order $1$ in the regime where the map $\check{m}$ is of order $\varepsilon^\frac{1}{2}$. When $\varepsilon$ is small enough, the function $m_2$ does not vanish in this regime, since the solution $\m$ is valued into the sphere $\SS^2$.
Assuming that $m_2$ is everywhere positive, it is given by the formula
$$m_2 = \big( 1 - \varepsilon |\Psi_\varepsilon|^2| \big)^\frac{1}{2},$$
and the function $\Psi_\varepsilon$ is a solution to the nonlinear Schr\"odinger equation
\begin{equation}
	\tag{NLS$_\varepsilon$}
	\label{NLS-eps}
	i \partial_t \Psi_\varepsilon + \big( 1 - \varepsilon |\Psi_\varepsilon|^2 \big)^{1/2} \Delta \Psi_\varepsilon + \frac{|\Psi_\varepsilon|^2}{1 + (1 - \varepsilon |\Psi_\varepsilon|^2)^{1/2}} \Psi_\varepsilon + \varepsilon \div \Big( \frac{\langle \Psi_\varepsilon, \nabla \Psi_\varepsilon \rangle_\C}{(1 - \varepsilon |\Psi_\varepsilon|^2)^{1/2}} \Big) \Psi_\varepsilon = 0,
\end{equation}
where 	$\langle z_1, z_2 \rangle_\C =\Re (z_1 \bar{z}_2)$. As $\varepsilon \to 0$, the formal limit  is therefore the focusing cubic Schr\"odinger equation
\begin{equation}
	\tag{CS}
	\label{CS}
	i \partial_t \Psi + \Delta \Psi + \frac{1}{2} |\Psi|^2 \Psi = 0.
\end{equation}
The goal is to justify rigorously this cubic Schr\"odinger regime of the LL equation.
We  recall that \eqref{CS}  is locally well-posed in $H^k(\R^N)$, 
for $k \in \N$; we refer to~\cite{Cazenav0} for an extended review on this subject.
Going on with our rigorous derivation of the cubic Schr\"odinger regime, we now express the local well-posedness result in Theorem~\ref{thm:LL-Cauchy} in terms of the nonlinear Schr\"odinger equation~\eqref{NLS-eps} satisfied by the rescaled function $\Psi_\varepsilon$.
\begin{cor}[\cite{deLaGra3}]
	\label{cor:Cauchy-Psi-eps}
	Let $\varepsilon>0$, and $k \in \N$, with $k > N/2 + 1$. Consider a function $\Psi_\varepsilon^0 \in H^k(\R^N)$ such that
	\begin{equation}
		\label{cond:petit-Psi-eps}
		\varepsilon^{1/2} \, \big\| \Psi_\varepsilon^0 \big\|_{L^\infty} < 1.
	\end{equation}
	Then there exist $T_\varepsilon>0$ and a unique solution $\Psi_\varepsilon \in L^\infty([0, T], H^k(\R^N))$ to~\eqref{NLS-eps}, for any  $t\in (0, T_\varepsilon)$. Moreover,  
	the flow map $\Psi_\varepsilon^0 \mapsto \Psi_\varepsilon$ is Lipschitz continuous from $H^k(\R^N)$ to $\boC^0([0, T], H^{k - 1}(\R^N))$ for any  $ T\in (0, T_\varepsilon)$ and 
	the nonlinear Schr\"odinger energy $\gE_\varepsilon$ given by
	$$\gE_\varepsilon(\Psi_\varepsilon) = \frac{1}{2} \int_{\R^N} \bigg( |\Psi_\varepsilon|^2 + \varepsilon |\nabla \Psi_\varepsilon|^2 + \frac{\varepsilon^2 \langle \Psi_\varepsilon, \nabla \Psi_\varepsilon \rangle_\C^2}{1 - \varepsilon |\Psi_\varepsilon|^2} \bigg),$$
	is conserved along the flow.
	
\end{cor}



We are now in position to state the main result concerning the rigorous derivation of the cubic Schr\"odinger regime of the LL equation.

\begin{thm}[\cite{deLaGra3}]
	\label{thm:conv-CLS}
	Let $0 < \varepsilon < 1$, and $k \in \N$, with $k > N/2 + 2$. Consider two initial conditions $\Psi^0 \in H^k(\R^N)$ and $\Psi_\varepsilon^0 \in H^{k + 3}(\R^N)$, and set
	$$\boS_\varepsilon = \big\| \Psi^0 \big\|_{H^k} + \big\| \Psi_\varepsilon^0 \big\|_{H^k} + \varepsilon^\frac{1}{2} \big\| \nabla \Psi_\varepsilon^0 \big\|_{\dot{H}^k} + \varepsilon \big\| \Delta \Psi_\varepsilon^0 \big\|_{\dot{H}^k}.$$
There is  $A>0$, depending only on $k$, such that, if the initial data $\Psi^0$ and $\Psi_\varepsilon^0$ satisfy the condition
	\begin{equation}
		\label{cond:key2}
		A \, \varepsilon^\frac{1}{2} \, \boS_\varepsilon \leq 1,
	\end{equation}
	then there exists a time
		$T_\varepsilon \geq \frac{1}{A \boK_\varepsilon^2},$
	such that both the unique solution $\Psi_\varepsilon$ to~\eqref{NLS-eps} with initial data $\Psi_\varepsilon^0$, and the unique solution $\Psi$ to~\eqref{CS} with initial data $\Psi^0$ are well-defined on the time interval $[0, T_\varepsilon]$. Moreover, we have the error estimate, for any $ t \in [0, T_\varepsilon]$, 
	\begin{equation}
		\label{eq:est-error}
		\big\| \Psi_\varepsilon(t) - \Psi(t) \big\|_{H^{k - 2}} \leq \Big( \big\| \Psi_\varepsilon^0 - \Psi^0 \big\|_{H^{k - 2}} + A \varepsilon \boS_\varepsilon \big( 1 + \boS_\varepsilon^3 \big) \Big) \, e^{A \boS_\varepsilon^2 t}.
	\end{equation}
	
\end{thm}

In this manner, Theorem~\ref{thm:conv-CLS} establishes  rigorously   the convergence of the LL equation towards the cubic Schr\"odinger equation in any dimension.
It is certainly possible to show only convergence under weaker assumptions by using compactness arguments as for the derivation of similar asymptotic regimes (see e.g.~\cite{ShatZen1, ChirRou2, GermRou1} concerning Schr\"odinger-like equations).

Observe that smooth solutions for both the LL and the cubic Schr\"odinger equations are known to exist when the integer $k$ satisfies the condition $k > N/2 + 1$. The additional assumption $k > N/2 + 2$ in Theorem~\ref{thm:conv-CLS} is related to the fact that the proof of~\eqref{eq:est-error} requires a uniform control of the difference $\Psi_\varepsilon - \Psi$, which follows from the Sobolev embedding theorem of $H^{k - 2}(\R^N)$ into $L^\infty(\R^N)$.

Finally, the loss of two derivatives in the error estimate~\eqref{eq:est-error} can be partially recovered by combining standard interpolation theory. Under the assumptions of Theorem~\ref{thm:conv-CLS}, the solutions $\Psi_\varepsilon$ converge towards the solution $\Psi$ in $\boC^0([0, T_\varepsilon], H^s(\R^N))$ for any $0 \leq s < k$, when $\Psi_\varepsilon^0$ tends to $\Psi^0$ in $H^{k + 2}(\R^N)$ as $\varepsilon \to 0$, but the error term is not necessarily of order $\varepsilon$ due to the interpolation process.

Note here that condition~\eqref{cond:key2} is not really restrictive in order to analyze such a convergence. At least when $\Psi_\varepsilon^0$ tends to $\Psi^0$ in $H^{k + 2}(\R^N)$ as $\varepsilon \to 0$, the quantity $\boS_\varepsilon$ tends to twice the norm $\| \Psi^0 \|_{H^k}$ in the limit $\varepsilon \to 0$, so that condition~\eqref{cond:key2} is always fulfilled. Moreover, the error estimate~\eqref{eq:est-error} is available on a time interval of order $1/\| \Psi^0 \|_{H^k}^2$, which is similar to the minimal time of existence of the smooth solutions to the cubic Schr\"odinger equation. 

%


The proof of Theorem~\ref{thm:conv-CLS} is similar to the proof of Theorem~\ref{thm:conv-SG}.
It relies on the consistency between the Schr\"odinger equations~\eqref{NLS-eps} and~\eqref{CS} in the limit $\varepsilon \to 0$. Indeed, we can recast~\eqref{NLS-eps} as
\begin{equation}
	\label{eq:consistency}
	i \partial_t \Psi_\varepsilon + \Delta \Psi_\varepsilon + \frac{1}{2} |\Psi_\varepsilon|^2 \Psi_\varepsilon = \varepsilon \boR_\varepsilon,
\end{equation}
where the remainder term $\boR_\varepsilon$ is given by
\begin{equation}
	\label{def:R-eps}
	\boR_\varepsilon = \frac{|\Psi_\varepsilon|^2}{1 + ( 1 - \varepsilon |\Psi_\varepsilon|^2 )^\frac{1}{2}} \Delta \Psi_\varepsilon - \frac{|\Psi_\varepsilon|^4}{2 (1 + (1 - \varepsilon |\Psi_\varepsilon|^2)^\frac{1}{2})^2} \Psi_\varepsilon - \div \Big( \frac{\langle \Psi_\varepsilon, \nabla \Psi_\varepsilon \rangle_\C}{(1 - \varepsilon |\Psi_\varepsilon|^2)^\frac{1}{2}} \Big) \Psi_\varepsilon .
\end{equation}
In order to establish the convergence towards the cubic Schr\"odinger equation, the main goal is to control the remainder term $\boR_\varepsilon$ on a time interval $[0, T_\varepsilon]$ as long as possible. In particular, we have to show that the maximal time $T_\varepsilon$ for this control does not vanish in the limit $\varepsilon \to 0$. The main argument is to perform suitable energy estimates on the solutions $\Psi_\varepsilon$ to~\eqref{NLS-eps}. These estimates provide Sobolev bounds for the remainder term $\boR_\varepsilon$, which are used to control the differences $u_\varepsilon = \Psi_\varepsilon - \Psi$ with respect to the solutions $\Psi$ to~\eqref{CS}. This further control is also derived from energy estimates. 

Concerning the estimates of the solutions $\Psi_\varepsilon$, we rely on the equivalence with the solutions $\m$ to~\eqref{LL-ani:asypm}. However, the estimates given in Section~\ref{chap:cauchy} are not enough in this case. It is crucial to refine the estimate \eqref{eq:energy-estimate-LL}, 
which can be done when $\lambda_1 = \lambda_3$.

\begin{prop}
	\label{prop:LL-energy-estimate2}
	Let $0 < \varepsilon < 1$, and $k \in \N$, with $k > N/2 + 1$. Assume that 
$	\lambda_1 = \lambda_3 = {1/\varepsilon},$
	and that $\m$ is a solution to~\eqref{LL-ani:asypm} in $\boC^0([0, T], \boE^{k + 4}(\R^N))$, with $\partial_t m \in \boC^0([0, T], H^{k + 2}(\R^N))$. Given any integer $2 \leq \ell \leq k + 2$, the energies $E_{\ell}$  are of class $\boC^1$ on $[0, T]$, and there exists  $C_k>0$, depending possibly on $k$, but not on $\varepsilon$, such that their derivatives satisfy
	\begin{equation}
		\label{eq:energy-estimate-LL2}
E'_\ell(t)
 \leq \frac{C_k}{\varepsilon} \Big( \| m_1(t) \|_{L^\infty}^2 + \| m_3(t) \|_{L^\infty}^2 + \| \nabla \bmm(t) \|_{L^\infty}^2 \Big) \, \Big( E_\ell(t) + E_{\ell - 1}(t) \Big),
	\end{equation}
	for any $t \in [0, T]$. Here we have set $E_1(t)=E(\bmm(t))$, the LL energy.
\end{prop}

As for the proof of  Proposition~\ref{prop:LL-energy-estimate}, the estimates in Proposition~\ref{prop:LL-energy-estimate2} rely on the  identity \eqref{eq:second-LL}, that in the case 
$\lambda_1 = \lambda_3=1/\ve$ can be simplified.  In contrast with the estimate \eqref{eq:energy-estimate-LL}, the multiplicative factor in the right-hand side of~\eqref{eq:energy-estimate-LL2} now only depends on the uniform norms of the functions $m_1$, $m_3$ and $\nabla \m$. This property is key in order to use these estimates in the cubic Schr\"odinger regime.

Finally, it is necessary to find a high order energy,
with suitable cancellation properties to obtain good energy estimates. 
The energy proposed in \cite{deLaGra2}, which allows us  to conclude as in the sine--Gordon equation, is 
\begin{equation*}
	\begin{split}
	\gE_\varepsilon^k (t) =&
	 \big\| \Psi_\varepsilon \big\|_{\dot{H}^{k - 2}}^2 + \big\| \varepsilon \partial_t \Psi_\varepsilon - i \Psi_\varepsilon  \big\|_{\dot{H}^{k - 2}}^2 + \varepsilon^2 \big\| \Delta \Psi_\varepsilon  \big\|_{\dot{H}^{k - 2}}^2\\
&		+  \varepsilon \big( \big\| \partial_t (1 - \varepsilon |\Psi_\varepsilon |^2)^\frac{1}{2} \big\|_{\dot{H}^{k - 2}}^2 + \big\| \Delta (1 - \varepsilon |\Psi_\varepsilon |^2)^\frac{1}{2} \big\|_{\dot{H}^{k - 2}}^2 + 2 \big\| \nabla \Psi_\varepsilon  \big\|_{\dot{H}^{k - 2}}^2 \big),
	\end{split}
\end{equation*}
for any $k \geq 2$. We refer to  \cite{deLaGra2} for detailed computations.

\section{Stability of sum of solitons}
\label{chap:stability}

In dimension one, the LL equation is completely integrable by means of the inverse scattering method \cite{FaddTak0} and, using this technique,
explicit solitons and multisolitons solutions can be constructed \cite{BikBoIt1}.
We consider in this section equation \eqref{LL-1d}, i.e.\  the one-dimensional easy-plane LL
equation. By a change of variable, we assume that $\lambda_3=1$.

We say that a soliton for \eqref{LL-1d} is a traveling wave  of the form 
$\m(x, t) = \u(x - c t).$
The nonconstant solitons are explicitly given by
\begin{equation}
	\label{form:uc}
	\u_c(x) =( c\,{\sech \big( \sqrt{1 - c^2} x \big)}, 
	\tanh \big( 
	\sqrt{1 - c^2} x \big), 
	\sqrt{1 - c^2}\,{\sech\big( \sqrt{1 - c^2} x \big)}
	), \quad  \abs{c}<1, 
\end{equation}
up to the invariances of the equation, i.e.\ translations, rotations around the axis $x_3$ and orthogonal symmetries with respect to any line in the plane $x_3 = 0$. Thus a soliton with speed $c$ may be also written as
$$\u_{c, a, \theta, s}(x) = \big( \cos(\theta) [\u_c]_1 - s \sin(\theta) [\u_c]_2, \sin(\theta) [\u_c]_1 + s \cos(\theta) [\u_c]_2, s [\u_c]_3 \big)(x - a),$$
with $a \in \R$, $\theta \in \R$ and $s \in \{ \pm 1 \}$. We refer to~\cite{deLaGra1,deLaGra2,deLaGra3} for more properties of solitons for the LL equation  \eqref{LL-ani}.

In addition, using the integrability of the equation
and by means of the inverse scattering method, for any $M\in \N^*$, it can be also computed explicit solutions to \eqref{LL-1d} that behave like a sum of $M$ decoupled solitons as $t\to \infty$.
These solutions are often called $M$-solitons or simply multisolitons  (see e.g. \cite[Section 10]{BikBoIt1} for their explicit formula).

We can define properly the solitons in the hydrodynamical framework when $c \neq 0$, since the function $\check{u}_c=[{\u}_c]_1+i [{\u}_c]_2$ does not vanish. More precisely, we recall 
that for a function $\u:\R\to \SS^2$ such that $\abs{\u}\neq 0$, we set  $\check{u} = (1 - u_3^2)^{1/2} i \exp(- i \varphi),$ and we define 
the hydrodynamical variables $v =u_3$ and $w =- \partial_x \varphi$.
Thus, equation \eqref{LL-1d} recasts as in
\eqref{HLL-1d},
and the soliton $\u_c$ in the hydrodynamical variables $\gv_c = (v_c, w_c)$ is given by
\begin{equation}
	\label{form:vc}
	v_c(x) = \sqrt{1 - c^2}\,{\sech\big( \sqrt{1 - c^2} x \big)}, \
	{\rm and} \ w_c(x) = \frac{c \, v_c(x)}{1 - v_c(x)^2} = \frac{c \sqrt{1 - c^2} \cosh \big( \sqrt{1 - c^2} x \big)}{\sinh \big( \sqrt{1 - c^2} x \big)^2 + c^2}.
\end{equation}
Therefore,  the only remaining invariances of solitons in this framework are translations and the opposite map $(v, w) \mapsto (- v, - w)$. Any soliton with speed $c$ may be then written as
$\gv_{c, a, s}(x) = s \, \gv_c(x - a) = ( s \, v_c(x - a), s \, w_c(x - a) ),$
with $a \in \R$ and $s \in \{ \pm 1 \}$.

Our goal in this section is to establish the  stability of a single soliton $u_c$ along the LL flow. 
More generally, we will also consider the case of a sum of solitons. In the original framework, defining this sum is not so easy, since 
the sum of unit vectors in $\R^3$ does not necessarily remain in   $\SS^2$. In the hydrodynamical framework, this difficulty does not longer arise. We can define a sum of $M$ solitons $S_{\gc, \ga, \gs}$ as
$$\SSS_{\gc, \ga, \gs} = (V_{\gc, \ga, \gs}, W_{\gc, \ga, \gs}) = \sum_{j = 1}^M \gv_{c_j, a_j, s_j},$$
with $M \in \N^*$, $\gc = (c_1, \ldots, c_M)$, $\ga = (a_1, \ldots, a_M) \in \R^M$, and $\gs = (s_1, \ldots, s_M) \in \{ \pm 1 \}^M$. However, we have to restrict the analysis to speeds $c_j \neq 0$, since the function $\check{u}_0$, associated with the black soliton, vanishes at the origin.

Coming back to the original framework, we can define properly a corresponding sum of solitons $\RRR_{\gc, \ga, \gs}$, when the third component of $\SSS_{\gc, \ga, \gs}$ does not reach the values $\pm 1$. Due to the exponential decay of the functions $v_c$ and $w_c$, this assumption is satisfied at least when the positions $a_j$ are sufficiently separated, i.e. when the solitons are decoupled. In this case, the sum $\RRR_{\gc, \ga, \gs}$ is given, up to a phase factor, by the expression
$$\RRR_{\gc, \ga, \gs} = \Big( (1 - V_{\gc, \ga, \gs}^2)^\frac{1}{2} \cos(\Phi_{\gc, \ga, \gs}),(1 - V_{\gc, \ga, \gs}^2)^\frac{1}{2} \sin(\Phi_{\gc, \ga, \gs}), V_{\gc, \ga, \gs} \Big),\quad \text{ with } \Phi_{\gc, \ga, \gs}(x) = \int_0^x W_{\gc, \ga, \gs}(y) \, dy,$$
for any $x \in \R$. This definition presents the advantage to provide a quantity with values on the sphere $\SS^2$. On the other hand, it is only defined under restrictive assumptions on the speeds $c_j$ and positions $a_j$. Moreover, it does not take into account the geometric invariance with respect to rotations around the axis $x_3$.

\subsection{Orbital stability in the energy space}
\label{sub:result}
In the sequel, our main results are proved in the hydrodynamical framework. We establish that, if the initial positions $a_j^0$ are well-separated and the initial speeds $c_j^0$ are ordered according to the initial positions $a_j^0$, then the solution corresponding to a chain of solitons at initial time, that is a perturbation of a sum of solitons $S_{\gc^0, \ga^0, \gs^0}$, is uniquely defined, and that it remains a chain of solitons for any positive time.

Let us recall that Theorem~\ref{thm:local-hydro-Cauchy} provides the existence and uniqueness of a continuous flow for \eqref{HLL-1d} in the nonvanishing energy space $\boN\boV(\R)$.
To our knowledge, the question of the global existence (in the hydrodynamical framework) of the local solution $\gv$ is open. In the sequel, we by-pass this difficulty using the stability of a well-prepared sum of solitons $\SSS_{\gc, \ga, \gs}$. Since the solitons in such a sum have exponential decay by \eqref{form:vc}, and are sufficiently well-separated, the sum $\SSS_{\gc, \ga, \gs}$ belongs to $\boN\boV(\R)$. Invoking the Sobolev embedding theorem, this remains true for a small perturbation in $H^1(\R) \times L^2(\R)$. As a consequence, the global existence for a well-prepared sum of solitons follows from its stability by applying a continuation argument.

Concerning the stability of  sums of solitons, our main result is

\begin{thm}[\cite{deLaGra1}]
	\label{thm:mult-stab}
	Let $\gs^* \in \{ \pm 1 \}^\M$ and $\gc^* = (c_1^*, \ldots, c_\M^*) \in ((-1, 1) \setminus \{ 0 \})^\M$ such that
$c_1^* <c_2^*  < \cdots < c_\M^*.$
	There exist positive numbers $\alpha^*$, $L^*$, $\nu$ and $A$,  depending only on $\gc^*$ such that, if $\gv^0 \in \boN\boV(\R)$ satisfies the condition
	\begin{equation}
		\label{def:alpha0}
		\alpha := \big\| \gv^0 - \SSS_{\gc^*, \ga^0, \gs^*} \big\|_{H^1 \times L^2} \leq \alpha^*,
	\end{equation}
	for points $\ga^0 = (a_1^0, \ldots, a_\M^0) \in \R^\M$ such that
$
		L^0 := \min \big\{ a_{j + 1}^0 - a_j^0, 1 \leq j \leq \M - 1 \big\} \geq L^*,
$
	then the solution $\gv$ to \eqref{HLL-1d} with initial condition $\gv^0$ is globally well-defined on $\R_+$, and there exists a function $\ba = (a_1, \ldots, a_\M) \in \boC^1(\R_+, \R^\M)$ such that,
		for any $t\geq 0$,
		\begin{equation}
		\label{est:a'}
		\sum_{j = 1}^\M \big| a_j'(t) - c_j^* \big| \leq A \big( \alpha + e^{ - {\nu L^0}} \big),
\ \textup{ 	and }\
		\big\| \gv(\cdot, t) - \SSS_{\gc^*, \ba(t), \gs^*} \big\|_{H^1 \times L^2} \leq A \big( \alpha +
		e^{ - {\nu L^0}} \big).
	\end{equation}

\end{thm}

Theorem~\ref{thm:mult-stab} provides the orbital stability of well-prepared sums of solitons with different, nonzero speeds for positive time. The sums are well-prepared in the sense that their positions at initial time are well-separated and ordered according to their speeds. As a consequence, the solitons are more and more separated along the LL flow (see estimate \eqref{est:a'}) and their interactions become weaker and weaker. The stability of the chain then results from the orbital stability of each single soliton in the chain.

As a matter of fact, the orbital stability of a single soliton appears as a special case of Theorem~\ref{thm:mult-stab} when $\M=1$.
In this case, stability occurs for both positive and negative times due to the time reversibility of the LL equation. Time reversibility also provides the orbital stability of reversely well-prepared chains of solitons for negative time. The analysis of stability for both negative and positive time is more involved. It requires a deep understanding of the possible interactions between the solitons in the chain (see \cite{MartMer7, MartMer8} for such an analysis in the context of the KdV equation). This issue is of particular interest because of the existence of multisolitons.

Special chains of solitons are indeed provided by the exact multisolitons. However, there is a difficulty to define them properly in the hydrodynamical framework. Indeed, multisolitons can reach the values $\pm 1$ at some times. On the other hand, an arbitrary multisoliton becomes well-prepared for large time in the sense that the individual solitons are ordered according to their speeds and well-separated (see e.g.\ \cite[Section 10]{BikBoIt1}).

If we consider  a perturbation of an arbitrary multisoliton at initial time, our theorem does not guarantee that  a perturbation of this  multisoliton remains 
a perturbation of a multisoliton for large time. In fact, this property would follow from the continuity with respect to the initial datum of LL equation in the energy space, which remains, to our knowledge, an open question. We remark that  Theorem~\ref{thm:mult-stab} only shows the orbital stability of the multisolitons, which do not reach the values $\pm 1$ for any positive time. 

To our knowledge, the orbital stability of the soliton $\u_0$ remains an open question. In the context of the Gross--Pitaevskii equation, the orbital stability of the vanishing soliton (often called black soliton) was proved in \cite{BeGrSaS1, GeraZha1}. Part of the analysis in this further context certainly extends to the soliton $\u_0$ of the LL equation.

Let us remark that in case $\lambda_3 = 0$, there is no traveling-wave solution to \eqref{LL-1d} with nonzero speed and finite energy. However, breather-like solutions were found to exist in \cite{LakRuTh1}, and their numerical stability was investigated in \cite{TjonWri1}. In the easy-axis case, there are traveling-wave solutions (see e.g.\ \cite{BishLon1}), but their third coordinate $m_3(x)$ converges to $\pm 1$ as $|x| \to + \infty$. This prevents from invoking the hydrodynamical formulation, and thus from using the strategy developed below in order to prove their orbital stability.
In the rest of this section, we 

We present now the main elements in the proof of Theorem~\ref{thm:mult-stab},
restricting our attention to the analysis of a single soliton.
We  underline that these arguments  do not make use of the inverse scattering transform. Instead, they rely on the Hamiltonian structure of the LL equation, in particular, on the conservation laws for the energy and momentum. As a consequence, these arguments can presumably be extended to nonintegrable equations similar to the hydrodynamical LL equation.

The strategy of the proof of Theorem~\ref{thm:mult-stab} is reminiscent of the one developed to tackle the stability of well-prepared chains of solitons for the generalized KdV equations \cite{MarMeTs1}, the nonlinear Schr\"odinger equations \cite{MarMeTs2}, or the Gross-Pitaevskii equation \cite{BetGrSm1}.
A key ingredient in the proof is the minimizing nature of the soliton $\gv_c$, which can be constructed as the solution of the minimization problem
\begin{equation}
	\label{def:mini-E-P}
	E(\gv_c) = \min \big\{ E(\gv) \mid  \gv \in \boN\boV(\R) \ {\rm s.t.} \ P(\gv) = P(\gv_c) \big\},
\end{equation}
where we recall that the energy and the momentum of  $\gv=(v,w)$, are given by 
\begin{equation*}
	E(\gv) = \frac{1}{2} \int_\R \Big( \frac{(v')^2}{1 - v^2} + \big( 1 - v^2 \big) w^2 + v ^2 \Big),
	\ \textup{ and }
	P(\gv) = \int_\R v w.
\end{equation*}

This characterization results from the compactness of the minimizing sequences for \eqref{def:mini-E-P}, and the classification of solitons in \eqref{form:vc}. The compactness of minimizing sequences can be proved following the arguments developed for a similar problem in the context of the Gross--Pitaevskii equation \cite{BetGrSa2,delaire-mennuni}.

The Euler--Lagrange equation for \eqref{def:mini-E-P} reduces to the identity
$E'(\gv_c) = c P'(\gv_c),$
where the speed $c$ appears as the Lagrange multiplier of the minimization problem. The minimizing energy is equal to
$E(\gv_c) = 2 (1 - c^2)^\frac{1}{2},$
while the momentum of the soliton $\gv_c$ is given by
$	P(\gv_c) = 2 \arctan ({(1 - c^2)^\frac{1}{2}} /{c}),$ for $c \neq 0$. An important consequence is the 
inequality
\begin{equation}
	\label{eq:der-Pc}
	\frac{d}{dc} \Big( P(\gv_c) \Big) = - \frac{2}{(1 - c^2)^\frac{1}{2}} < 0,
\end{equation}
which is related to the Grillakis--Shatah--Strauss condition (see e.g.\ \cite{GriShSt1}) for the orbital stability of a soliton. As a matter of fact, we can use inequality \eqref{eq:der-Pc} to establish the coercivity of the quadratic form
$$Q_c = E''(\gv_c) - c P''(\gv_c),$$
under suitable orthogonality conditions. More precisely, we show
\begin{prop}
	\label{prop:coer-single}
	Let $c \in  (- 1, 1) \setminus \{ 0 \}$. There exists  $\Lambda_c>0$, such that
	\begin{equation}
		\label{eq:coer-Qc}
		Q_c(\beps) \geq \Lambda_c \| \beps \|_{H^1 \times L^2}^2,
	\end{equation}
	for any pair $\beps \in H^1(\R) \times L^2(\R)$ satisfying the two orthogonality conditions
	\begin{equation}
		\label{eq:ortho-Qc}
		\langle \partial_x \gv_c, \beps \rangle_{L^2 \times L^2} = \langle P'(\gv_c), \beps \rangle_{L^2 \times L^2} = 0.
	\end{equation}
	Moreover, the map $c \mapsto \Lambda_c$ is uniformly bounded from below on any compact subset of
	 $(- 1, 1) \setminus \{ 0 \}$.
\end{prop}

The first orthogonality condition in \eqref{eq:ortho-Qc} originates in the invariance with respect to translations of \eqref{HLL-1d}. Due to this invariance, the pair $\partial_x \gv_c$ lies in the kernel of $Q_c$. The quadratic form $Q_c$ also owns a unique negative direction, which is related to the constraint in \eqref{def:mini-E-P}. This direction is controlled by the second orthogonality condition in \eqref{eq:ortho-Qc}.

As a consequence of Proposition~\ref{prop:coer-single}, the functional
$F_c(\gv) = E(\gv) - c P(\gv),$
controls any perturbation $\beps = \gv - \gv_c$ satisfying the two orthogonality conditions in \eqref{eq:ortho-Qc}. More precisely, we derive from the Euler--Lagrange equation and \eqref{eq:coer-Qc} that
\begin{equation}
	\label{eq:control-eps}
	F_c(\gv_c + \beps) - F_c(\gv_c) \geq \Lambda_c \| \beps \|_{H^1 \times L^2}^2 + \boO \big( \| \beps \|_{H^1 \times L^2}^3 \big),
\end{equation}
as $\| \beps \|_{H^1 \times L^2} \to 0$. Since 
the energy $E(\gv)$ and the momentum $P(\gv)$ are conserved along the flow, the left-hand side of \eqref{eq:control-eps} remains small for all time if it was small at the initial time. As a consequence of \eqref{eq:control-eps}, the perturbation $\beps$ remains small for all time, which implies the stability of $\gv_c$. We refer to \cite{deLaGra1}
for more detail about the proof of Theorem~\ref{thm:mult-stab}.

\subsection{Asymptotic stability}

We consider now the long-time asymptotics  of a solution 
to \eqref{LL-1d}, with initial condition a perturbation of a soliton.
We would like to determine conditions such that the solution converges to a (possible different) soliton.
Let us remark that the convergence as $t\to\infty$ cannot 
hold in the energy space. For instance, we could consider 
a solution $\gv$ to \eqref{HLL-1d} with an initial condition 
$\gv^0 \in \boN\boV(\R)$, such that $\gv$  converges to a hydrodynamical soliton  $\gv_c$ in the norm $\| \cdot \|_{H^1 \times L^2}$, as $t\to\infty$. By the  continuity of the energy and the momentum (with respect to this norm), we have
$$E \big( \gv(\cdot, t) \big) \to E(\gv_c) \quad {\rm and} \quad 
P \big( \gv(\cdot, t) \big) \to P(\gv_c),$$
as $t \to  \infty$. Since these quantities are conserved by the flow, 
we conclude that $E(\gv^0)=E(\gv_c)$
and $P(\gv^0)=P(\gv_c)$.
Thus, the variational  characterization of solitons implies 
that $\gv^0$ must be a soliton. Therefore, 
the only solutions that converge (in energy norm) to a soliton 
as $t\to\infty$, are the solitons.

In conclusion, to establish the asymptotic stability, we  need to weaken the notion of  convergence. Indeed, using the weak  convergence in the space $\boN\boV(\R)$, Bahri~\cite{Bahri1}  proved 
the asymptotic stability of solitons in the hydrodynamical framework.

\begin{thm}[\cite{Bahri1}]
	\label{thm:asympt-HLL}
	Let $c \in (-1, 1)\setminus\{0\}$. There is  $\alpha^* > 0$ such that, if  the initial condition $\gv^0 \in \boN\boV(\R)$ satisfies that
$\| \gv^0 - \gv_c \|_{H^1 \times L^2} < \alpha^*,$
	then there exist a unique global associated solution $\gv \in \boC^0(\R, \boN\boV(\R))$ to \eqref{LL-1d},  $c^* \in (-1, 1)\setminus\{0\}$ and 
	$a \in \boC^1(\R, \R)$ such that, as $t \to  \infty$, 
	$$\gv(\cdot + a(t), t) \big) \rightharpoonup \gv_{c^*} \quad \text{ in }  H^1(\R) \times L^2(\R), \quad \text{ and }\quad 
	a'(t) \to c^*.$$
\end{thm}


This theorem provides the weak convergence towards a soliton, but this long-time dynamics 
needs to take into account the geometric invariances of the problem, 
i.e.\ the translations. This is precisely the role of the parameter $a(t)$, whose derivative  converges to the speed of the limit soliton $\gv_{c^*}$. In this fashion, the solution propagates with the same speed as the  limit soliton, as $t$ goes to infinity, as expected.

The weak convergence in Theorem~\ref{thm:asympt-HLL} can probably be improved. Indeed, Martel and Merle~\cite{MartMer4, MartMer6} proved the asymptotic stability of solitons of the KdV equation, establishing a locally (strong) convergence in the energy space. It is possible that a similar result can be shown for the  asymptotic stability of hydrodynamical solitons of the LL equation satisfy a similar, i.e.\  a strong convergence in a norm of the type
${H^1([- R(t), R(t)]) \times L^2([- R(t), R(t)])}$, 
where $R(t)$ is a linear function of time.

The proof of Theorem~\ref{thm:asympt-HLL} is based on an approach developed by Martel and Merle for the KdV  equation \cite{MartMer4, MartMer6}. Their strategy can be decomposed in three steps, that we would explain in our context, i.e.\ in the hydrodynamical setting. First, the orbital stability provided by Theorem~\ref{thm:mult-stab}
guarantees 
that a  solution $\gv$, with initial  condition $\gv^0$
close  enough to a soliton  $\gv_c$, remains in a neighborhood of the orbit of the soliton. In particular, the solution $\gv$ is bounded in the nonvanishing space  $\boN\boV(\R)$ for any $t\geq 0$. It is then possible 
to construct a sequence of times  $(t_n)$, with $t_n\to\infty$, 
and a limit function  $\gv_*^0 \in \boN\boV(\R)$, 
such that, up to a subsequence, 
$$\gv(\cdot, t_n) \rightharpoonup \gv_*^0 \quad {\rm in} \quad H^1(\R) \times L^2(\R),$$
as $n \to  \infty$. In  addition,  $\gv_*^0$ remains close  to the orbit of the  soliton $\gv_c$. Moreover,  the solution $\gv_*$ to \eqref{HLL-1d} with initial condition $\gv_*^0$ is global, and is also close to this orbit.
We point out that is also necessary to introduce  a modulation parameter 
due to the invariance by translation, but we will omit it for the sake of clarity.

We need to prove that the limit profile  $\gv_*^0$, and the associated  solution  $\gv_*$, are indeed solitons. Thus, the second step 
is to study the regularity and decay properties of $\gv_*$. To this end, 
it is useful to establish the weak   continuity of the flow of the hydrodynamical equation with respect to the initial condition, which implies that the solution 
$\gv$ converges to $\gv_*$, i.e.\  for any  $t\in \R$ (fixed),
$$\gv(\cdot, t_n + t) \rightharpoonup \gv_*(\cdot, t) \quad {\rm in} \ H^1(\R) \times L^2(\R), \ 
\text{ as } n \to  \infty.$$
Using also a monotonicity formula for the momentum, from  this convergence it is possible to deduce that   $\gv_*$ is localized in space, uniformly in time, and that 
 $\gv_*$ has an  exponential decay in space, uniformly in time. Thus, using the Kato smoothing effect
that gives regularizing properties of the  Schr\"odinger-type equations, it follows that $\gv_*$ is of class $\boC^\infty$ on $\R \times \R$, and that all its derivatives also decay in space, uniformly in time. 

The third step is to show that in the neighborhood of a soliton, the only solutions
to \eqref{HLL-1d} having this behavior are the solitons. This rigidity property follows from a Liouville type theorem. The proof of this theorem requires 
another monotonicity formula, and it is the most difficult part of the argument.
We refer to \cite{Bahri1} for more details.

By refining the approach described above, Bahri~\cite{Bahri2}  also established 
the asymptotic stability for initial data close to a
sum of solitons, that  are as usual well-prepared
according to their speeds and have sufficiently separated initial positions. 
The proof of this theorem relies on the strategy developed by Martel, Merle and Tsai 
in~\cite{MarMeTs1} for the KdV equation. Let us also remark that the locally strong asymptotic stability result for multisolitons in~\cite{MarMeTs1}  is stronger than the statement in \cite{Bahri2}  with $M=2$. Indeed, the proof
in~\cite{MarMeTs1} is based on a monotonicity argument for the localized energy. 
It is an open problem if  this kind of argument can be adapted to the study of the LL equation, or more generally, if it possible to get  a locally strong asymptotic stability result. 

In the higher dimensional case $N\geq 2$, most of the questions about solitons are still open.
We refer to \cite{deLaire4} and the references therein for more details.

\section{Self-similar solutions for the LLG equation}
\label{chap:self}
In this section we will study the dissipative  LLG equation 
\eqref{LLG-intro}.
We will focus on the existence of self-similar solutions and provide their asymptotics in dimension $N=1$. We also
analyze the qualitative and
quantitative effect of the damping $\alpha$  on the dynamical behavior of these self-similar solutions.

As we will see, these kinds of solutions do not belong to classical Sobolev spaces, and we cannot invoke the Cauchy theory developed in Section~\ref{chap:cauchy} to give a meaning to their stability.
Therefore, we will provide a well-posedness result in a more general framework related to the BMO space to give some stability results.
We point out  that the proof of the  well-posedness result 
uses the parabolic behavior of the equation in presence of damping, 
and cannot be applied for the pure dispersive equation 
(i.e.\ $\alpha =0$) analyzed in previous sections.

\subsection{Self-similar solutions}
A natural question, that has been proven relevant for  understanding the global behavior of solutions and formation of singularities, is whether  there exist solutions which are invariant under scalings of the equation. In the case of equation \eqref{LLG-intro},  it is straightforward to see that it is invariant under the following scaling:  If  $\mm$  is a solution of \eqref{LLG-intro}, then   $\mm_{\lambda}(x,t)=\mm(\lambda x,\lambda^2 t)$  is also a solution, 
for any  $\lambda>0$. Associated with this invariance, a solution $\mm$ of \eqref{LLG-intro} defined for $I=\R^+$ or $I=\R^-$ is called {\it{self-similar}} if it is invariant under rescaling, that is
\bq\label{scaling}
\mm(x,t)= \mm(\lambda x,\lambda^2 t),\quad \forall \lambda>0, \quad \forall x\in \R^N,\quad  \forall t\in I.
\eq
Setting $T\in \R$ and performing a translation in time, this definition leads to  two types of self-similar solutions: 
A forward self-similar solution, or {\em expander}, is a solution of the form
$	\mm(x,t)=\f({x}/{\sqrt{t-T} })$ for $(x,t)\in \R^N\times(T,\infty)$,
and  a backward self-similar solution, or {\em shrinker},  is a solution of the form
$	\mm(x,t)=\f({x}/{\sqrt{T-t} })$ for $(x,t)\in \R^N\times(-\infty, T),$
for certain profile $\f:\mathbb{R}^N\longrightarrow \mathbb{S}^2$. Expanders evolve from a singular value at time $T$, while shrinkers evolve towards a singular value at time $T$. 

Self-similar solutions have brought a lot of attention in the study on nonlinear PDEs
because they can provide some important information about the dynamics of the equation.
While expanders are related to nonuniqueness phenomena, resolution of singularities and 
long time description of solutions, shrinkers 
are often related to phenomena of singularity formation (see e.g.\ \cite{book-giga,eggers}). On the other hand, the construction and understanding of the dynamics and properties of self-similar solutions also provide an idea of which are the natural spaces to develop a well-posedness theory, that captures these very often physically relevant structures.  
Examples of equations for which self-similar solutions have  been considered, and a substantial  work around these types of solutions has been done, include among others the Navier--Stokes equation, semilinear parabolic equations, and geometric flows such as Yang-Mills, mean curvature flow and harmonic map flow. We refer to \cite{jia-tsai,struwe96} and the references therein for more details.


Most of the works in the literature related to the study of self-similar solutions to the  LLG equation are confined to the heat flow for harmonic maps equation, i.e.\ $\alpha=1$. In this setting,  the main works on the subject restrict the analysis to corotational maps
taking values in $\SS^d$, which reduces the analysis of \eqref{HFHM} to the study of a second order real-valued ODE. Then tools such as the maximum principle or the shooting method can be used to  show  the existence of solutions. We refer to \cite{fan,gastel,germain-rupflin,biernat-donninger,bizon-wasserman,biernat-bizon,germain-ghoul} for more details on such  results for 
maps taking values in $\SS^d$, with $d\geq 3.$ Recently, 
Deruelle and Lamm~\cite{deruelle-lamm} have studied the Cauchy problem 
for  the harmonic map heat flow with initial data $\mm^0 : \R^N\to \SS^d$, with $N\geq 3$
and $d\geq 2$, where $\mm^0$ is Lipschitz 0-homogeneous function, homotopic to a constant,
which implies the existence of expanders coming out of $\mm^0$.

When $0<\alpha\leq 1$,  we established  the existence of self-similar {expanders} for the LLG equation   in \cite{gutierrez-delaire2}. This result  is a consequence of a well-possedness theorem for the  LLG equation 
considering an initial data $\mm^0 : \R^N\to \SS^2$ in the space BMO
of functions of bounded mean oscillation.
Notice that this result  includes in particular the case of the  harmonic map heat flow.
We will explain more precisely this result in 
Section~\ref{sec:BMO}.

As seen before, in absence of damping ($\alpha=0$), \eqref{LLG-intro} reduces to  the Schr\"odinger map equation \eqref{SM}, which is reversible in time, so that the notions of expanders and shrinkers coincide.
For this equation, Germain, Shatah and Zeng \cite{germain-shatah-zeng} established the
existence of ($k$-equivariant) self-similar profiles $\f : \R^2\to \SS^2$.

In the one-dimensional  case, when $\alpha=0$, \eqref{SM} is closely related Localized Induction Approximation (LIA), and self-similar profiles $\f: \R\to \SS^2$ 
were obtained and analyzed in \cite{vega-gutierrez,vega-gutierrez1,lakshmanan0}.
In the context of LIA, self-similar solutions constitute a family of smooth solutions that develop a singularity in the shape of a corner in finite time. For further work related to these solutions, including the study of the continuation of these solutions after the blow-up time and their stability, we refer to the reader to \cite{banica-vega-3,banica-vega}. At the level of the Schr\"odinger map equation, these self-similar solutions provide examples of smooth solutions that develop a jump singularity in finite time.  

In this section we explain how to construct 
the family of expanders profiles for  $\alpha\in [0,1]$, and provide their analytical study
and we  discuss the Cauchy problem associated with these solutions and their stability. Finally, we construct and  analyze the family of shrinkers profiles. 

\subsection{Expanders in dimension one}
\label{sec:expanders}
We consider in this section  equation \eqref{LLG-intro} in dimension $N=1$, 
and $\alpha\in[0,1]$, in order to include both the damped and undamped cases.
We seek self-similar solutions of the form
\begin{equation}\label{self-similar}
	\mm(x,t)=\gm ({x}/{\sqrt{t}}),\qquad
	x\in \R, \ t>0,
\end{equation}
and we will say that $\gm$ is the profile of the solution $\mm$.
Observe that if  $\m$ is a smooth  solution to \eqref{LLG-intro}, it can be checked that 
$\gm$ solves the  following system of ODEs
\begin{equation}
	\label{eq:EDO:intro}
	\alpha \gm''+\alpha\abs{\gm'}^2\gm +\beta (\gm\times \gm')'+\frac{x \gm'}2=0,
	\quad\textup{ on }\R,
\end{equation}
due to the fact that $\gm$ takes values in $\SS^2$.
Thus, we can give a weak formulation to this equation in the form 
$-(\boA(\gm) \gm')'=G(x,\gm,\gm')$, with 
$$ \boA(\u)=\begin{pmatrix}
	\alpha & -\beta u_3 & \beta u_2\\
	\beta u_3 & \alpha & -\beta u_1\\
	-\beta u_2 & \beta u_1 & \alpha,
\end{pmatrix}, \ 
G(x,\u,\bm p)=
\begin{pmatrix}
	\alpha u_1 \abs{\bm p}^2-\frac{xp_1}2 \\
	\alpha u_2 \abs{\bm p}^2-\frac{xp_2}2\\
	\alpha u_3 \abs{\bm p}^2-\frac{xp_3}2\\
\end{pmatrix},
$$
where $\u =(u_1,u_2,u_3)$ and $\bm p=(p_1,p_2,p_3)$. 

Therefore, if $\alpha>0$, the system is  uniformly elliptic, 
since $\boA(\u) \bm \xi \cdot \bm \xi= \alpha \abs{\bm\xi}^2$, for all $\bm\xi,\u\in \R^3$,
and we can then  invoke the regularity theory for quasilinear elliptic systems, to verify that the solutions are smooth.

In the limit case $\alpha=0$, we can show directly that the solutions are also smooth.
Most importantly, we have the following theorem that provides a rigidity result concerning the possible solutions to \eqref{eq:EDO:intro}: The modulus of the  gradient of any  solution {\em must} be $ce^{-\alpha x^2/4}$, for some $c\geq 0$.
\begin{thm}[\cite{gutierrez-delaire3}]
	\label{thm:EDO}
	Let $\alpha\in [0,1]$. Assume that $\gm\in H^{1}_{\loc}(\R;\SS^2)$ is a weak solution to \eqref{eq:EDO:intro}. 
	Then $\gm$ belongs to $\boC^\infty(\R;\SS^2)$ and there exists $c\geq 0$ such that
	$\abs{\gm'(x)}=ce^{-\alpha x^2/4}$, for all $x\in \R$. 
\end{thm}

In the limit cases $\alpha=1$ and $\alpha=0$, it is possible to find 
explicit solutions to 
\eqref{eq:EDO:intro}, as we will see later on.
However, this seems unlikely in the case $\alpha\in (0,1)$,
and even the existence of such solutions is not clear.
We proceed now to give a way of establishing the {\em existence} of solutions 
satisfying the condition 	$\abs{\gm'(x)}=ce^{-\alpha x^2/4}$, 
for any $c>0$  and any $\alpha\in[0,1]$ (notice that the case $c=0$  
 corresponds to the trivial constant solution).

The idea is to look for $\gm$ as the tangent vector to a curve in $\R^3$, so we  
first  recall some facts about curves in the space. 
Given   $\gm : \R\to \SS^2$ a smooth function, we can define the curve 
\bq
\label{def:curve}
\bm X_{\gm}(x)=\int_0^x \gm(s)ds,
\eq
so that $\bm X_{\gm}$ is  smooth, parametrized by arclenght,
and its  tangent vector is $\gm$. In addition, if $\abs{\gm'}$ does not vanish on $\R$, 
we can define the normal vector 
${\gn(x)}={\gm'(x)}/{\abs{\gm'(x)}}$ and  
the binormal vector $\gb(x)=\gm(x)\times \gn(x)$.
Moreover, we can define the curvature and torsion of $\bm X_{\gm}$ 
as $k(x)={\abs{\gm'(x)}}$ and $\tau(x)=-{\gb'(x)\cdot \gn(x)}$.
Since $\abs{\gm(x)}^2=1,$ for all $x\in\R$, we have that $\gm(x)\cdot \gn(x)=0$, for all $x\in \R$,
that the vectors $\{\gm,\gn,\gb\}$ are orthonormal and it is standard to check that they satisfy the Serret--Frenet system
\begin{equation}\label{profile-serret}
	\begin{aligned}
		\gm'=k\gn, \quad 		\gn'=-k\gm +\tau \gb,\quad 
		\gb'=-\tau \gn.
	\end{aligned}
\end{equation}

Let us apply this method to find a solution to \eqref{eq:EDO:intro}.
We define $\bm X_{\gm}$ 
as in \eqref{def:curve},
and we remark  that 
equation \eqref{eq:EDO:intro} rewrites in terms of $\{\gm, \gn , \gb   \}$ as 
$$
-\frac{x}{2}k\gn=\beta(k'\gb-\tau k\gn )-\alpha(-k'\gn-k\tau\gb).
$$
Therefore, from the orthogonality of the vectors $\gn$ and $\gb$,  we conclude that 
the curvature and torsion of $\bm X_{\gm}$ are solutions of the equations
$-x k=2\alpha k'-\beta \tau k$ and $\beta k'+\alpha k\tau=0,$
that is
\begin{equation}
	\label{profile-kt}
	k(x)= c e^{-{\alpha x^2}/{4}}\qquad {\hbox{and}}\qquad \tau(x)={\beta x}/{2},
\end{equation}
for some $c\geq 0$. Of course, the fact that  $k(x)= c e^{{-\alpha x^2}/{4}}$
is in agreement with  $\abs{\gm'(x)}=c e^{-\alpha x^2/4}$.

Now, given $\alpha\in[0,1]$ and $c>0$, consider the Serret--Frenet system \eqref{profile-serret} with curvature and torsion function given by \eqref{profile-kt}
and initial conditions
${\gm}(0)=\be_1$, ${\gn}(0)=\be_2$, ${\gb}(0)=\be_3$.
Then, by standard ODE theory, there exists a unique global solution $\{\gm_{c,\alpha}, \gn_{c,\alpha}, \gb_{c,\alpha}\}$ in $(\mathcal{C}^{\infty}(\mathbb{R}; \mathbb{S}^2))^{3}$, and these vectors are orthonormal. Also, it is straightforward to verify that $\gm_{c,\alpha}$
is a solution to  \eqref{eq:EDO:intro} satisfying  
$\abs{\gm_{c,\alpha}'(x)}=c e^{-\alpha x^2/4}$.

Finally, using the uniqueness of the Cauchy--Lipschitz theorem and the Serret--Frenet system, 
it is simple to  show the uniqueness of such  solutions, up to rotations. 
\begin{thm}[\cite{gutierrez-delaire3}]
	\label{prop:unicidad}
	The set of nonconstant solutions to 
	\eqref{eq:EDO:intro} is $\{\boR \gm_{c,\alpha} : c>0,\boR \in SO(3)\}$,
	where $SO(3)$ is the group  of rotations about the origin preserving orientations.
\end{thm}

The above proposition reduces the study of expanders to the understanding  
of the family of expanders associated with the 
profiles $\{\gm_{c,\alpha}\}_{c,\alpha}$. 
The next result summarizes  the properties of these solutions.
\begin{thm}[\cite{gutierrez-delaire1}]
	\label{Theorem1} Let $\alpha\in[0,1]$, $c\geq 0$
	and ${\gm}_{c,\alpha}$ be the solution of the Serret--Frenet system  constructed above.
 Let 
$\mm_{c,\alpha}(x,t)=
	\gm_{c,\alpha}\left({x}/{\sqrt{t} }\right)$, for $(x,t)\in \R\times(0,\infty).$
	Then the following statements hold.
	\begin{enumerate}[label=\textup{(\roman*)},ref=\textup{({\roman*})}]
		\item\label{regular} The function $\mm_{c, \alpha}$ is a  $\mathcal{C}^\infty$-solution of \eqref{LLG-intro}  on $\R\times (0,\infty)$, 
		with 		$	\abs{\partial_{x}\m_{c,\alpha} (x,t)}=\frac{c}{\sqrt t}e^{-{\alpha x^2}/{4t}}.$

		\item\label{converge} There exists a unitary vector $\A^{+}_{c,\alpha}=(A_{j,c,\alpha}^{+})_{j=1}^{3} \in \SS^2$ such that $\mm_{c,\alpha}(\cdot ,t )$ converges pointwise to the initial condition 
		\begin{equation} \label{data}
			\m^0_{c,\alpha}=\A^+_{c,\alpha}\chi_{\R^+}+\A^-_{c,\alpha}\chi_{\R^-},
		\end{equation}
		i.e.\ 
		\begin{equation}\label{convergencia} 
			\lim_{t\to0^+}\mm_{c,\alpha}(x,t)=\A^+_{c,\alpha}, \text{ if }x>0,
		\	\text{ and } \			\lim_{t\to0^+}\mm_{c,\alpha}(x,t)=\A^-_{c,\alpha}, \text{ if }x<0,
		\end{equation}
		where $\A^-_{c,\alpha}=(A_{1,c,\alpha}^+,-A_{2,c,\alpha}^+,-A_{3,c,\alpha}^+)$ and 
				 $\chi_E$ is the characteristic function of the set $E$.
		\item\label{rate} Moreover, there exists a constant $C(c,\alpha,p)$ such that for
		all $t>0$ and all for all $p\in (1,\infty)$,
		\begin{equation}\label{cota-m} \|\mm_{c,\alpha}(\cdot,t)-
				\m^0_{c,\alpha} \|_{L^p(\R)}\leq
			C(c,\alpha,p) t^\frac{1}{2p}. \end{equation} 
\end{enumerate}
\end{thm}

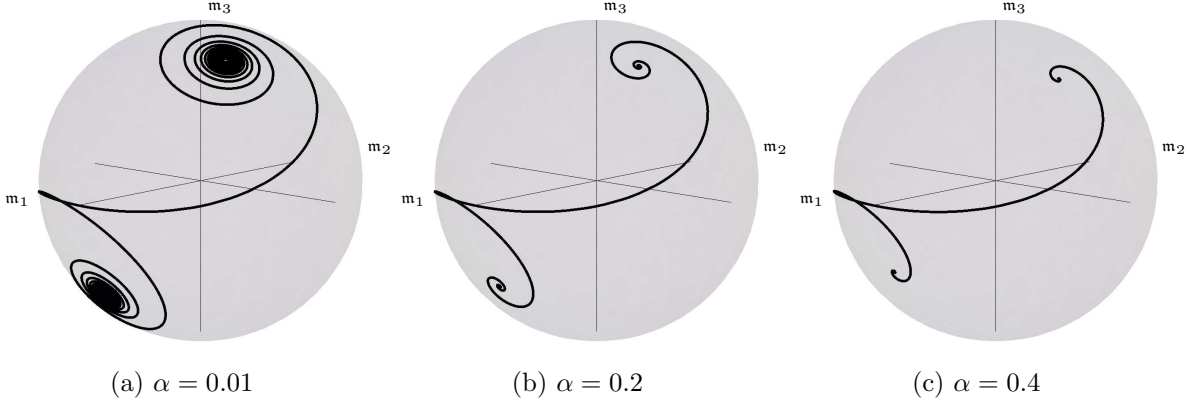
\begin{figure}[ht]
\smallskip
	\begin{subfigure}[c]{0.3\textwidth}
		\centering
		\begin{overpic}[trim=20 45 0 0,clip,scale=0.67
			]{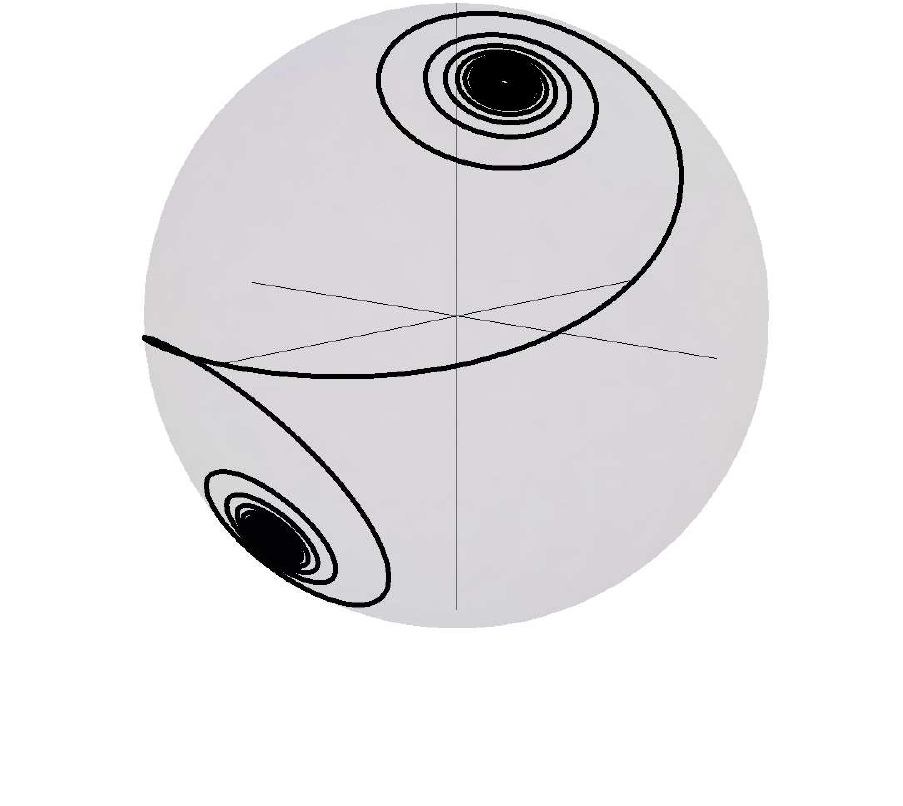}
			\put(1,35){\tiny{$\gm_1$}}
			\put(85,47){\tiny{$\gm_2$}}
			\put(48,80){\tiny{$\gm_3$}}
		\end{overpic}
		\caption{$\alpha=0.01$}
	\end{subfigure}
	\hspace{0.15cm}
	\begin{subfigure}[c]{0.3\textwidth}
		\centering
		\begin{overpic}[trim=20 45 0 0,clip,scale=0.67]{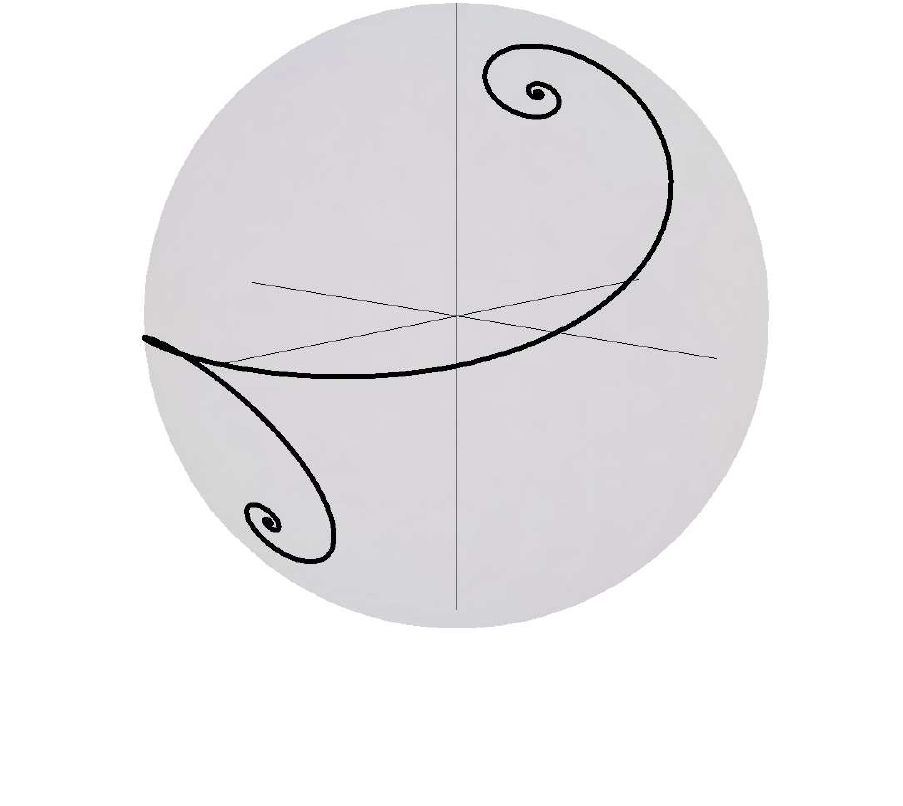}
\put(1,35){\tiny{$\gm_1$}}
\put(85,47){\tiny{$\gm_2$}}
\put(48,80){\tiny{$\gm_3$}}		\end{overpic}
		\caption{$\alpha=0.2$}
	\end{subfigure}
	\hspace{0.19cm}
	\begin{subfigure}[c]{0.3\textwidth}
		\centering
		\begin{overpic}[trim=20 45 0 0,clip,scale=0.67]{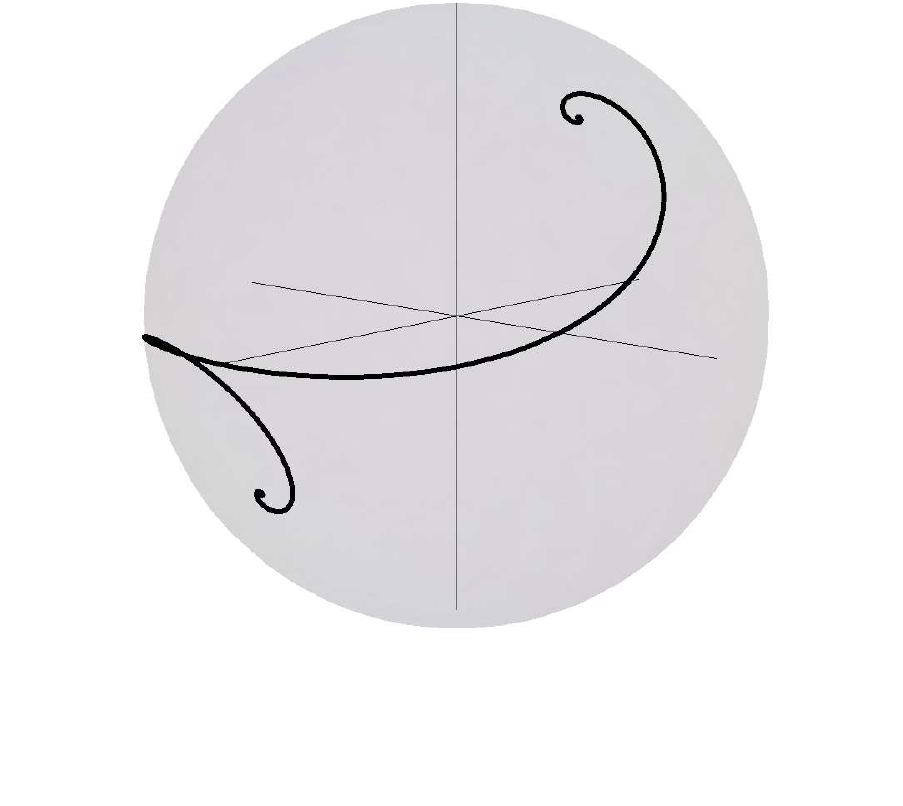}
\put(1,35){\tiny{$\gm_1$}}
\put(85,47){\tiny{$\gm_2$}}
\put(48,80){\tiny{$\gm_3$}}		\end{overpic}
		\caption{$\alpha=0.4$}
	\end{subfigure}
	\caption{The profile $\gm_{c,\alpha}$ for $c=0.8$ and different values of $\alpha$.}
	\label{fig-tan}
\end{figure}

The graphics in Figure~\ref{fig-tan} depict the profile
$\gm_{c,\alpha}$ for fixed $c=0.8$ and the values of $\alpha=0.01$, $\alpha=0.2$,
and $\alpha=0.4$. In particular, it can be observed how the
convergence of $\gm_{c,\alpha}$ to
$\A^\pm_{c,\alpha}$ is accelerated by the diffusion $\alpha$.
%

Notice that the initial condition
$\m^0_{c,\alpha}$
has a jump singularity at the point $x=0$ whenever the vectors
$\A^+_{c,\alpha}$ and $\A^-_{c,\alpha}$ satisfy
$\A^+_{c,\alpha}\neq \A^-_{c,\alpha}.$
In this situation (and we will be able to prove analytically that  this is the case, at
least for certain ranges of the parameters $\alpha$ and $c$, see Proposition~\ref{jump} below),
Theorem~\ref{Theorem1} provides a family of global
smooth solutions of \eqref{LLG-intro} associated with a discontinuous singular  initial
data (jump-singularity).

As  already mentioned, in the absence of damping ($\alpha=0$), singular
self-similar solutions of the Schr\"odinger map equation were previously obtained in \cite{vega-gutierrez,lakshmanan0}.
In this framework,
Theorem~\ref{Theorem1} establishes the persistence of a jump singularity for self-similar solutions in the presence of dissipation.

When $\alpha=0$, the stability of the self-similar solutions was considered in a series of papers by Banica and Vega \cite{banica-vega,Banica-Vega-2,banica-vega-Arch}.
The stability in the case $\alpha>0$ is a natural question that we will
discuss later.

Some further remarks on the results stated in
Theorem~\ref{Theorem1} are in order.
First, the  energy  is given  by
\begin{align*}
	\label{total-energy}
	E_{\text{LLG}}(t)=  \frac{1}{2}\int_{-\infty}^{\infty} |\partial_x\mm_{c,\alpha}(x,t)|^2\, dx=
	\frac{1}{2}\int_{-\infty}^{\infty} \left(  \frac{c}{\sqrt{t}} e^{-\frac{\alpha x^2}{4t}}\right)^2\, dx=
	c^2\sqrt{\frac{\pi}{\alpha t}},\quad  t>0.
\end{align*}
It follows that the  energy at the initial time $t=0$ is infinite, while it 
  becomes finite for all positive times, showing
the dissipation of energy in the system in the presence of damping.

Secondly, it is also important to remark that in the setting of Schr\"odinger
equations,  for fixed $\alpha\in[0,1]$ and $c>0$, the solution
$\mm_{c,\alpha}$  is associated through the Hasimoto
transformation  with the filament function  \cite{hasimoto}, that is 
$	u_{c,\alpha}(x,t)=\frac{c}{\sqrt{t}}e^{(-\alpha+i\beta)\frac{x^2}{4t}},
$
which solves
\begin{equation}
	\label{Schrodinger}
	i\partial_{t}u+(\beta-i\alpha)\partial_{xx}u+\frac{u}{2}\left(\beta\abs{u}^2+2\alpha\int_0^x
	\Im(\bar u \partial_{x}u)-A(t)\right)=0,\quad {\hbox{with}}\quad A(t)=\frac{\beta c^2}{t},
\end{equation}
with initial condition a Dirac delta function since
$\lim_{t\to 0^+} u_{c,\alpha}(x,t)=2c\sqrt{\pi(\alpha+i\beta)}\delta_0$.

Therefore $u_{c,\alpha}$ is very rough at initial time
and  the standard arguments (e.g.\ a Picard iteration scheme based on Strichartz
estimates and Sobolev--Bourgain spaces) cannot be applied, at least not straightforwardly, to study the local
well-posedness of the initial value problem for the Schr\"odinger equation (\ref{Schrodinger}).
The existence of solutions to  equation \eqref{Schrodinger} associated with an initial data proportional
to a Dirac delta opens the question of developing a well-posedness theory for Schr\"odinger equations of the type considered
here to include initial data of infinite energy.
In the case $\alpha=0$,  $A(t)=0$ and when the initial condition is proportional to the Dirac delta,
Kenig, Ponce and Vega \cite{kenig-ponce-vega} proved that the Cauchy problem for \eqref{Schrodinger} is ill-posed due to some oscillations. Moreover, even after removing these oscillations, Banica and Vega \cite{banica-vega} showed that   equation \eqref{Schrodinger} (with $\alpha=0$ and $A(t)=c^2/t$)
is still ill-posed. This question was also addressed by  Vargas and Vega in \cite{vargas-vega} and Gr\"unrock in \cite{grunrock}
for other types of initial data of infinite energy (see also \cite{banica-vega1}),
but we are not aware of any result in this setting when $\alpha>0$ (see \cite{GuoDing0} for related well-posedness results in the case $\alpha>0$
for initial data in Sobolev spaces of positive index).

\subsubsection{Asymptotics for the profile} 
\label{sub:asymp-expanders}

We want now to study  the qualitative and
quantitative effect of the damping $\alpha$ and the
parameter $c$ on the dynamical behavior of  the family 
$(\mm_{c,\alpha} )_{c,\alpha}$  of self-similar solutions
of \eqref{LLG-intro} found in Theorem~\ref{Theorem1}. Precisely, in an
attempt to fully understand the regularization of the solution at
positive times close to the initial time $t=0$, and to understand how
the presence of damping affects the dynamical behavior of these
self-similar solutions, we aim to give answers to the following
questions:
 Can we obtain a more precise behavior of the
	solutions $\m_{c,\alpha}$ at positive times $t$ close to
	zero?
 Can we understand the limiting vectors
	$\A^{\pm}_{c,\alpha}$ in terms of the parameters $c$
	and $\alpha$?

In order to address our first question, we observe that, due to the self-similar nature of these solutions, the behavior of the family of solutions
$\m_{c,\alpha}$ at positive times close to the initial time $t=0$ is directly related to the study of the asymptotics of the associated profile $\gm_{c,\alpha}(x)$ for large values of $\abs{x}$.
In addition, the symmetries of $\gm_{c,\alpha}$  (see Theorem~\ref{thm-conver} below) allow to reduce ourselves to obtain the behavior of the profile as $x\to\infty$. The precise asymptotics of the profile is given in the following theorem.
%
\begin{thm}[\cite{gutierrez-delaire1}] \label{thm-conver}
	Let $\alpha\in[0,1]$, $c>0$.	The components of ${\gm}_{c,\alpha}$ satisfy respectively that
$\gm_{1, c,\alpha}$ is an even function,  and $\gm_{j, c,\alpha}$ is an odd function for $j\in\{2,3\}$. In addition,  for all 	
		$s\geq s_0=4\sqrt{8+c^2}$,
		\begin{align*}
			\gm_{c,\alpha}(s)
			=&\A^+_{c,\alpha}-\frac{2c}{s}\B^+_{c,\alpha}e^{-\alpha s^2/4}(\alpha \sin(\bm{\phi}_{c,\alpha}(s))+\beta \cos(\bm{\phi}_{c,\alpha}(s)))
			\nonumber 			 
			-\frac{2c^2}{s^2} \A^+_{c,\alpha}e^{-\alpha s^2/2} +\boO\big(\frac{e^{-\alpha s^2/4}}{ s^3}\big).
		\end{align*}
		Here, $\sin(\bm{\phi}_{c,\alpha})$ and $\cos(\bm{\phi}_{c,\alpha})$ are understood
		acting on each of the components given by 
			\begin{equation}\label{def-phi} \phi_{j,c,\alpha}(s)=a_{j,\alpha,c}+\beta
			\int_{s_0^2/4}^{s^2/4}\sqrt{1+ c^2\frac{e^{-2\alpha
						\sigma}}{\sigma}}\,d\sigma, \quad j\in\{1,2,3\},
		\end{equation}
		for some constants $a_{1,\alpha,c},a_{2,\alpha,c},a_{1,\alpha,c}\in [0,2\pi)$, and the vector $\B^{+}_{c,\alpha}$ is given in terms of $\A^+_{c,\alpha}$ by
		$
		\B^+_{c,\alpha}=((1-(A_{1,c,\alpha}^+)^2)^{1/2},(1-(A_{2,c,\alpha}^+)^2)^{1/2},(1-(A_{3, c,\alpha}^+)^2)^{1/2}).$
		%
\end{thm}

The convergence and rate of convergence of the solutions $\mm_{c,\alpha}$  to $\mm_{c,\alpha}^0$ established in  Theorem~\ref{Theorem1} are simple  consequences of 
 the  asymptotics  in Theorem~\ref{thm-conver}. Also, similar asymptotics hold for 
 the normal vector $\gn_{c,\alpha}$ and the  binormal vector $\gb_{c,\alpha}$.

With regard to the asymptotics in Theorem~\ref{thm-conver}, it is important to mention
that the error depends only on $c$. More
precisely, we use the notation $\boO(f(s))$ to denote a function
for which there exists a constant  $C(c)>0$ depending on $c$, but
but on  $\alpha$, such that \begin{equation}\label{cota-rmk}\left|
	\boO\left(f(s)\right)\right|\leq C(c)\abs{f(s)}, \quad \text{for
		all } s\geq s_0. \end{equation}
At first glance, one might think that the term
$-2c^2 \A^+_{c,\alpha}e^{-\alpha s^2/2}/s^2$ in the asymptotics   could be included in the error term
$\boO(e^{-\alpha s^2/4}/s^3)$. However, we cannot do this because
 in our notation the big-$\boO$ must be independent of $\alpha$.

%

%
When $\alpha=1$ (so $\beta=0$),  we can solve explicitly
the Serret--Frenet system, to obtain 
\bq\label{explicit-for}
\gm_{c,1}(s)=(\cos(c\Erf(s)), \sin(c\Erf(s)),0),\eq
 for all $s\in \R$, where $\Erf$ is the non-normalized error function
$
\Erf(s)=\int_0^s e^{-\sigma^2/4}\, d\sigma.
$
In particular, the limiting vectors in Theorem~\ref{thm-conver} are given  by
\begin{equation}\label{limit-alpha-1}
	\A^\pm_{c,1}=(\cos(c\sqrt\pi),\pm \sin(c\sqrt\pi),0 ), \quad \B^+_{c,1}=(\abs{\sin(c\sqrt\pi)},\abs{\cos(c\sqrt\pi)},1).
\end{equation}

When $\alpha=0$, the solution of \eqref{profile-serret} can be
solved explicitly in terms of parabolic cylinder functions or
confluent hypergeometric functions (see \cite{GamayunLisovyy}). 
Another analytical approach using Fourier analysis techniques has been taken
in \cite{vega-gutierrez}, leading to the asymptotics
\begin{equation}\label{asymp-alpha-0}
		\gm_{c,0}(s)
		=\A^+_{c,0}-\frac{2c}{s}\B^+_{c,0}\sin( {\psi}_{c})
		+
		\boO\left(1/{s^2}\right),
\quad \text{ with }
	{\psi_{c}(s)}=\frac{s^2}4+c^2\ln(s).
\end{equation}
Moreover, 
$\bm A^+_{c,0}$ can be computed  explicitly. 
On the other hand, when $\alpha=0$, the phase $\bm\phi_{c,\alpha}$ in \eqref{def-phi} can be expanded as
\begin{equation*}
	\phi_{j,c,0}(s)=a_{j,c,\alpha}+\frac{s^2}{4}+c^2\ln(s)+C(c)+\boO\left({1/s^2}\right).
\end{equation*}
Thus the asymptotics in Theorem~\ref{thm-conv} allows us to  recover the logarithmic contribution in the oscillation in \eqref{asymp-alpha-0}.

When $\alpha>0$, $\bm\phi_{c,\alpha}$ behaves like
\begin{equation*} \phi_{j,c,\alpha}(s)=a_{j,c,\alpha}+
	\frac{\beta s^2}{4}+C(\alpha,c)+\boO\big(\frac{e^{-\alpha
			s^2/2}}{\alpha s^2}\big), \end{equation*}
and there is no logarithmic correction in the oscillations in the presence of damping.
Consequently, the phase function $\bm{\phi}_{c,\alpha}$  captures the different nature of the oscillatory character of the solutions in both the absence and the presence of damping.

It can be seen that the terms $\A^+_{c,\alpha}$, $\B^+_{c,\alpha}$, $\B^+_{c,\alpha} \cdot \sin(a_{c,\alpha})$, $\B^+_{c,\alpha}\cdot \cos(a_{c,\alpha})$ and
the error term
depend continuously on $\alpha\in [0,1]$.
Therefore,  the asymptotics in Theorem~\ref{thm-conv} shows how the profile $\gm_{c,\alpha}$ converges to $\gm_{c,0}$ as $\alpha\to 0^+$
and to $\gm_{c,1}$ as $\alpha\to 1^-$. In particular, we recover the asymptotics in  \eqref{asymp-alpha-0}.

%

%
%


Finally, the amplitude of the leading order term controlling the wave-like behavior  of the solution $\gm_{c,\alpha}(s)$ around
$\A^{\pm}_{c,\alpha}$ for values of $s$ sufficiently large is of the order $\displaystyle{c\, e^{-\alpha s^2/4}/{s}}$,
from which one observes how the convergence of the solution to its limiting values $\A^{\pm}_{c,\alpha}$ is accelerated in the presence of damping in the system, as depicted in Figure~\ref{fig-tan}.

Let us discuss now some results answering the second of our questions. Bearing in mind that 
$\A_{c,\alpha}^{-}$ is expressed in terms of the coordinates of $\A_{c,\alpha}^{+}$,  
we only need to focus on 
$\A^+_{c,\alpha}$. When $\alpha=1$ or $\alpha=0$, the vector
$\A^+_{c,\alpha}$ is  explicitly given in terms of the parameter $c$.  
When $\alpha\in(0,1)$, we do not have explicit expressions for these vectors, however
the following result establishes  that the solutions
$\mm_{c,\alpha}$ of the LLG equation found in Theorem~\ref{Theorem1} are indeed associated with a discontinuous initial data at least for certain ranges of $\alpha$ and $c$. 
\begin{thm}[\cite{gutierrez-delaire1}] \label{jump}
\leavevmode
	\begin{enumerate}
		\item Let $\alpha\in(0,1]$. There exists $c^{\ast}>0$ depending on $\alpha$ such that
		$
		\A_{c,\alpha}^{+}\neq \A_{c,\alpha}^{-}$, for all $c\in(0, c^{\ast}).$

		\item Let $c>0$. There exists $\alpha^{\ast}_{0}>0$  such that
		$
		\A_{c,\alpha}^{+}\neq \A_{c,\alpha}^{-}$, for all $\alpha\in(0, \alpha^{\ast}_{0}).$

		\item Let $c>0$, with $c\notin \N \sqrt{\pi}$. There exists $\alpha^{\ast}_{1}\in(0,1)$  such that
		$
		\A_{c,\alpha}^{+}\neq \A_{c,\alpha}^{-}$, for all $\alpha\in(\alpha^{\ast}_{1},1).
		$
	\end{enumerate}
\end{thm}
\begin{remark}\label{remk-jump}
	It can be checked that   $\A_{c,0}^{+}\neq \A_{c,0}^{-}$ for all $c>0$.
	Based on the numerical results in \cite{gutierrez-delaire1}, we conjecture that
	$\A_{c,\alpha}^{+}\neq \A_{c,\alpha}^{-}$ for all $\alpha\in (0,1)$ and $c>0$.
\end{remark}

%

Concerning, the proof of the asymptotics of $\gm_{c,\alpha}$, a key tool is  
a classical change of variables from the differential geometry of
curves that allows us to reduce the nine equations in the Serret--Frenet
system into three complex-valued second order equations (see e.g.\ 
\cite{lamb}). This change of variables is related
to the stereographic projection;  this approach was  used in
\cite{vega-gutierrez}. In our case, the change of variables  
reduces the analysis of the solution
$\{\gm_{c,\alpha},\gn_{c,\alpha},\gb_{c,\alpha}\}$
of the Serret--Frenet system  to the study
of three solutions to the  second order differential equation
\begin{equation}
	\label{eq-f0}
	f_{c,\alpha}''(s)+\frac{s}2(\alpha+i\beta)f_{c,\alpha}'(s)+\frac{c^2}{4} e^{-\alpha s^2/2}f_{c,\alpha}(s)=0,
\end{equation}
associated with three different initial conditions. The analysis of the solutions of
\eqref{eq-f0}  requires  the control of
certain integrals by exploiting their oscillatory character.
This can be achieved by using repeated integration by parts, in the spirit of the
method of stationary phase. We refer to \cite{gutierrez-delaire2} 
for more details of the proof.

\subsection{The Cauchy problem for LLG  in	BMO}
\label{sec:BMO}
A  natural question in the study of the stability properties of the family of solutions $( \m_{c,\alpha})_{c>0}$ is whether  it is possible to develop a well-posedness theory for the Cauchy problem for \eqref{LLG-intro} in a functional framework that allows us to handle initial conditions of the type \eqref{data}. In view of \eqref{data}, such a framework should allow some ``rough'' functions 
(i.e.\ function spaces beyond the ``classical'' energy ones) and step functions.

In the case $\alpha>0$,  global well-posedness results for \eqref{LLG-intro} have been established in $N\geq 2$ by
Melcher~\cite{melcher} and by Lin, Lai and Wang~\cite{lin-lai-wang}
for initial conditions with a smallness condition on the gradient in the $L^N(\R^N)$ and on the Morrey $M^{2,2}(\R^N)$-norm, respectively.
Therefore, these results do not apply to the initial condition $\m^0_{c,\alpha}$.
When $\alpha=1$, global well-posedness results for the heat flow for harmonic maps \eqref{HFHM} have been obtained by Koch and Lamm~\cite{koch-lamm}  for an initial condition $L^\infty$-close to a point and improved to an initial data with small
BMO semi-norm  by Wang~\cite{wang}. The ideas used in \cite{koch-lamm} and \cite{wang} rely on techniques introduced by
Koch and Tataru~\cite{koch-tataru} for the Navier--Stokes equation. Since  $\m^0_{c,\alpha}$
has a small BMO semi-norm if $c$ is small,  the results in \cite{wang} apply to the case $\alpha=1$.

In this subsection we explain the main results in \cite{gutierrez-delaire2}
that allow us to  adapt and extend the techniques developed
in \cite{koch-lamm,koch-tataru,wang} to prove a global well-posedness result for \eqref{LLG-intro}
with $\alpha\in(0,1]$,  for data $\m^0$ in $L^{\infty}(\R^N; \SS^2)$ with small BMO semi-norm.
As an application of these results, we can  establish the stability  of the family of self-similar solutions $(\m_{c,\alpha})_{c>0}$
and derive further properties for these solutions. 
In particular, we can  prove the existence of multiple smooth solutions of \eqref{LLG-intro} associated with the same initial condition, provided that $\alpha$ is close to one.

Our approach to study the Cauchy problem for  \eqref{LLG-intro} consists in analyzing the Cauchy problem
for the associated dissipative quasilinear Schr\"odinger equation through the stereographic projection, and then
``transferring'' the results back to the original equation.
To this end, we use  the stereographic projection from the South Pole defined in \eqref{def:stereo}. As mentioned in the introduction, 
if $\m$ is a smooth solution of \eqref{LLG-intro} with \mbox{$m_3>-1$}, then its stereographic projection $u=\P(\m)$
satisfies the quasilinear dissipative Schr\"odinger equation 
\eqref{DNLS}. At least formally, the Duhamel formula gives the integral equation:
\begin{equation}\label{duhamel}\tag\tag{IDNLS}
	u(x,t)=S_{\alpha}(t) u^0+\int_0^t S_\alpha(t-s)g(u)(s)\,ds,\quad  \text{ with }\ 
g(u)=-2i(\beta-i\alpha)\frac{\bar u (\grad u)^2}{1+\abs{u}^2},
\end{equation}
where $u^0=u(\cdot,0)$   corresponds to the initial condition,
and $S_\alpha(t)$ is the dissipative Schr\"odinger semigroup (also called the complex Ginzburg--Landau semigroup) given by
$S_\alpha(t)\phi=e^{(\alpha+i\beta)t\Delta}\phi$, i.e.
\begin{equation}\label{semigrupo}
	(S_\alpha(t)\phi)(x)=\int_{\R^N}G_\alpha(x-y,t)\phi(y)\,dy, \quad \text{with } \quad G_\alpha(x,t)=\frac{ e^{-\frac{\abs{x}^2}{4(\alpha+i\beta)t} } }{ (4\pi(\alpha+i\beta) t)^{N/2}}.
\end{equation}

One difficulty in studying  \eqref{duhamel} is to handle the term $g(u)$.
We see that $|g(u)|\leq |\grad u|^2,$
so we need to control $|\grad u|^2$. Koch and Taratu dealt with a similar problem when studying the well-posedness for the Navier--Stokes equation in \cite{koch-tataru}.
Their approach was to introduce some new spaces related to BMO  and $\BMO^{-1}$. Later, Koch and Lamm \cite{koch-lamm}, and  Wang \cite{wang} have adapted these spaces to study some geometric flows. Following these ideas, we define the Banach spaces
\begin{align*}
	X(\R^N\times \R^+;F)&=\{v:\R^N\times \R^+\to F \ : \ v,\grad v\in L^1_{\loc}(\R^N\times \R^+), \  \norm{v}_{X}<\infty\}
	\qquad \text{and} \\
	Y(\R^N\times \R^+;F)&=\{v:\R^N\times \R^+\to F \ : v\in L^1_{\loc}(\R^N\times \R^+), \  \norm{v}_{Y}<\infty\},
\end{align*}
where $\norm{v}_{X}=\sup_{t>0}\norm{v}_{L^\infty}+[v]_X$, {with}
\begin{align*}
		[v]_{X}=&\sup_{t>0} \sqrt{t}\norm{\grad v}_{L^\infty}+
	\sup_{\substack{ x\in \R^N\\ r>0}}\left( \frac1{r^N}\int_{Q_r(x)} \abs{\grad v(y,t) }^2\,dt\,dy \right)^\frac12, \qquad\text{and}\\
	\norm{v}_{Y}=&\sup_{t>0} t\norm{v}_{L^\infty}+
	\sup_{\substack{ x\in \R^N\\ r>0}}\frac1{r^N}\int_{Q_r(x)} \abs{ v(y,t) }\,dt\,dy.
\end{align*}
Here $Q_r(x)$ denotes the parabolic ball $Q_r(x)=B_r(x)\times [0,r^2]$ and $F$ is either $\C$ or $\R^3$.
The absolute value stands for the complex absolute value if $F=\C$
and for the euclidean norm if $F=\R^3$. We denote with the same symbol the absolute value
in $F$ and $F^3$.

The spaces $X$ and $Y$ are related to the spaces $\BMO(\R^N)$ and $\BMO^{-1}(\R^N)$ and are well-adapted
to study problems involving the heat semigroup $S_1(t)=e^{t\Delta}$. In order to establish the properties of the semigroup $S_\alpha(t)$
with $\alpha\in (0,1]$, we introduce the spaces
$\BMO_\alpha(\R^N)$ and $\BMO_\alpha^{-1}(\R^N)$ as the space of distributions $f\in S'(\R^N;F)$
such that the semi-norm and norm  given respectively  by
\begin{align*}
	[f]_{BMO_\alpha}=\sup_{\substack{ x\in \R^N\\ r>0}}\Big( \frac1{r^N}\int_{Q_r(x)} \abs{\grad S_\alpha(t)f}^2 \Big)^\frac12,\ \
	\norm{f}_{BMO^{-1}_\alpha}=\sup_{\substack{ x\in \R^N\\ r>0}}\Big( \frac1{r^N}\int_{Q_r(x)}\abs{ S_\alpha(t)f}^2 \Big)^\frac12,
\end{align*}
are finite.

On the one hand, the Carleson measure characterization of BMO functions (see \cite[Chapter~4]{stein} and
\cite[Chapter~10]{lemarie})
yields that for fixed $\alpha\in(0,1]$, ${\BMO_\alpha}(\R^N)$ coincides with the classical $\BMO(\R^N)$ space, 
 that is for all $\alpha\in(0,1]$ there exists a constant  $\Lambda>0$ depending only on $\alpha$ and $N$ such that
$	\Lambda[f]_{BMO}\leq [f]_{BMO_\alpha}\leq \Lambda^{-1}[f]_{BMO}.$
%
On the other hand, Koch and Tataru proved in \cite{koch-tataru} that $\BMO^{-1}$ (or equivalently $\BMO^{-1}_{1}$, using our notation) can be characterized as the space of derivatives of functions in BMO. A straightforward  generalization of their argument shows that the same result holds for $\BMO_\alpha^{-1}$. Hence, using the Carleson measure characterization theorem, we conclude that $\BMO_\alpha^{-1}$
coincides with the space  $\BMO^{-1}$ and that there exists a constant  $\tilde\Lambda>0$, depending only on $\alpha$ and $N$, such that
$	\tilde\Lambda\norm{f}_{BMO^{-1}}\leq  \norm{f}_{BMO_\alpha^{-1}} \leq \tilde\Lambda^{-1}\norm{f}_{BMO^{-1}}.$

\begin{figure}[ht!]
	\begin{center}
		\begin{overpic}[trim=0 3 1 0,clip,scale=0.9,
			]{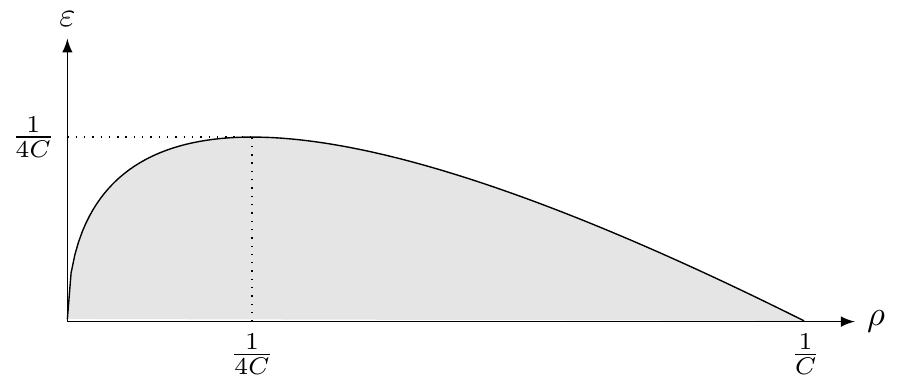}
			\end{overpic}
	\end{center}
	\caption{The shape of the set $\mathcal S(C)$.}
	\label{fig-set}
\end{figure}

The above remarks allow us to use several of the estimates proved in \cite{koch-lamm,koch-tataru,wang} in the case $\alpha=1$,
to study the integral equation \eqref{duhamel} by using a fixed-point approach.
Finally, this leads to the  next result that provides  the global well-posedness of the Cauchy problem for \eqref{duhamel} with small initial data in $\BMO(\R^N)$.
%
\begin{thm}[\cite{gutierrez-delaire2}]\label{thm:cauchy:ss}
	Let $\alpha\in(0,1]$. There exist constants $C,K\geq 1$ such that for every $L\geq 0$, $\ve>0$, and $\rho>0$ satisfying
	\begin{equation}\label{hyp1}
		8C(\rho+\ve)^2\leq \rho,
	\end{equation}
	if $u^0\in L^\infty(\R^N;\C)$,  with
	\begin{equation}\label{cond-CI}
		\norm{u^0}_{L^\infty}\leq L \quad\textup{ and }\quad [u^0]_{BMO}\leq \ve,
	\end{equation}
	then there exists a unique  solution $u\in X(\R^N\times \R^+;\C)$ to \eqref{duhamel} such that
	\begin{equation}\label{cond-small}
		[u]_{X}\leq K(\rho+\ve).
	\end{equation}
	Moreover,  $u\in \boC^\infty(\R^N\times \R^+)$, \eqref{DNLS} holds pointwise, $\sup_{t>0}\norm{u}_{L^\infty}\leq K(\rho+L)$ and 
 $ u(\cdot,t)\to u^0$, as ${t\to0^+}$,  as tempered distributions.

In addition, assume that $u$ and $v$ are respectively
		solutions to \eqref{duhamel} fulfilling \eqref{cond-small} with initial conditions $u^0$ and $v^0$
		satisfying \eqref{cond-CI}. Then	$			\norm{u-v}_{X}\leq 6K \norm{u^0-v^0}_{L^\infty}.$
\end{thm}
Although  condition \eqref{hyp1} appears naturally from the fixed-point used in the proof, it may not be so clear at first glance. To better understand it, let us define for $C>0$
\begin{equation}\label{set-S}
	\mathcal S(C)=\{ (\rho,\ve)\in \R^+\times \R^+ : C(\rho+\ve)^2\leq \rho\}.
\end{equation}
We see that if $(\rho,\ve)\in \mathcal S(C)$, then $\rho,\ve>0$ and
$	\ve\leq \frac{\sqrt{\rho}}{\sqrt{C}}-\rho.$
Therefore, the set $\mathcal S(C)$ is non-empty and bounded. The shape of this set is depicted in Figure~\ref{fig-set}. In particular, we infer that if $(\rho,\ve)\in \mathcal S(C)$, then
$\rho\leq \frac1C$  and $\ve \leq \frac1{4C}.$
In addition, if $\tilde C\geq C$, then $\mathcal S(\tilde C)\subseteq \mathcal S(C).$

Moreover, taking  $\rho=1/(32C)$, Theorem~\ref{thm:cauchy:ss} asserts that for fixed $\alpha\in(0,1]$,
we can take for instance $\ve=1/(32C)$ (that depends  on $\alpha$ and $N$, but not on the $L^{\infty}$-norm of the initial data)
such that for any given initial condition $u^0\in L^{\infty}(\R^N)$ with $[u^0]_{BMO}\leq \ve$,  there exists a global (smooth) solution $u\in X(\R^N\times \R^{+}; \C)$ of \eqref{DNLS}. Notice that $u^0$ is allowed to have a large $L^\infty$-norm as long as $[u^0]_{BMO}$ is sufficiently small;
this is a weaker requirement that asking for  the $L^\infty$-norm of $u^0$
to be sufficiently small, since
$[f]_{BMO}\leq 2\norm{f}_{L^\infty}$, for all $f\in L^\infty(\R^N).$

We remark that the smallness condition in \eqref{cond-small} is necessary for the uniqueness of the solution.	As we will see in Theorem~\ref{thm-non-uniq}, at least in dimension one, it is possible to construct multiple solutions
	of \eqref{duhamel} in $X(\R^N\times \R^+;\C)$, if $\alpha$ is close enough to 1.

By using the inverse of the stereographic projection
$\P^{-1} :\C\to \SS^2\setminus\{0,0,-1\}$, that is explicitly given by
$\m=(m_1, m_2, m_3)=\P^{-1}(u)$, with
\begin{equation}\label{inverse-P}
	m_1=\frac{2\Re u}{1+\abs{u}^2}, \quad m_2=\frac{2\Im u}{1+\abs{u}^2}, \quad m_3=\frac{1-\abs{u}^2}{1+\abs{u}^2},
\end{equation}
we can deduce from Theorem~\ref{thm:cauchy:ss} a
global well-posedness result for \eqref{LLG-intro}.
Moreover, the choice of the South Pole is of course arbitrary. By using the invariance of \eqref{LLG-intro} under rotations,
we have the existence of solutions provided that the essential range of the initial condition  $\m^0$ is far from an arbitrary point $\QQ\in \SS^2$.

\begin{thm}[\cite{gutierrez-delaire2}]\label{thm:cauchy-LLG}
	Let $\alpha\in(0,1]$. There exist constants $C\geq 1$ and $K\geq 4$, such that if $\delta\in (0,2]$,    $\ve_0,\rho>0$, $\delta\in (0,2]$,  $\ve_0>0$ and  $\rho>0$
	satisfy
	\begin{equation}\label{cond-LLG}
		8K^4C \delta^{-4}(\rho+8\delta^{-2}\ve_0)^2\leq \rho,
	\end{equation}
the following holds. Given any   $\m^0=(m_1^0,m_2^0,m_3^0)\in L^\infty(\R^N;\SS^2)$ and any $\QQ\in \SS^2$ satisfying
	\begin{equation*}
		\inf_{\R^N}\abs{\m^0-\QQ}^2\geq 2\delta \quad\textup{ and }\quad[\m^0]_{BMO}\leq \ve_0,
	\end{equation*}
	there exists a unique smooth solution $\m\in X(\R^N\times \R^+;\SS^2)$ of \eqref{LLG-intro}
	with initial condition $\m^0$ such that
	\begin{equation}\label{cond-cor}
		\inf_{\substack{x\in \R^N\\ t>0}} \abs{\m(x,t)-\QQ}^2\geq \frac{4}{1+K^2(\rho+\delta^{-1})^2}
		\quad \textup{ and }\quad [\m]_{X}\leq 4K(\rho+8\delta^{-2}\ve_0).
	\end{equation}
\end{thm}

We point out that the results are valid only for $\alpha>0$. If we let $\alpha\to0^+$, then
the estimates  blow up. Indeed, the proofs rely on 
the exponential decay of the semigroup $e^{(\alpha+i\beta)t\Delta}$, so that these techniques cannot be generalized (at least not straightforwardly) to cover the critical
case $\alpha=0$. In particular, we cannot recover the stability results for the
self-similar solutions in the case of Schr\"odinger maps proved by Banica and Vega in
\cite{banica-vega,Banica-Vega-2,banica-vega-Arch}.

As mentioned before, in \cite{lin-lai-wang} and \cite{melcher} some global well-posedness results for \eqref{LLG-intro} with $\alpha\in (0,1]$ were proved for initial conditions with small gradient in $L^N(\R^N)$ and $M^{2,2}(\R^N)$, respectively.
In view of the embeddings 
$$L^N(\R^N)\subset M^{2,2}(\R^N)\subset BMO^{-1}(\R^N),$$
for $N\geq2$, Theorem~\ref{thm:cauchy-LLG} can be seen as generalization of these results
since it covers the case of less regular initial conditions.
The arguments in \cite{lin-lai-wang,melcher} are based on the method of moving frames that produces a covariant
complex Ginzburg--Landau equation. 

The existence and uniqueness results given by Theorem~\ref{thm:cauchy-LLG} require the initial condition to be small in the 
BMO semi-norm. Without this condition, the solution 
could develop a singularity in finite time. In fact, in dimensions $N=3,4$, Ding and Wang \cite{ding-wang}  proved that
for some smooth  initial conditions with small energy, the associated solutions of \eqref{LLG-intro}  blow up in 
finite time.

Another consequence of Theorem~\ref{thm:cauchy-LLG} is the existence of self-similar solutions 
of expander type in $\R^N$, in any dimension  $N\geq 1$, 
i.e.\ a solution $\m$ of the form 
$\m(x,t)=\f({x}/{\sqrt t}),$
for some profile $\f:\mathbb{R}^N\to \mathbb{S}^2$.
In particular, we have the relation
$\f(y)=\m(y,1)$, for  all $y\in \R^N$.
From the scaling \eqref{scaling}, we see that, at least formally,
a necessary condition for the existence of a self-similar solution is that
initial condition $\m^0$ be homogeneous of degree 0, i.e.\
$\m^{0}(\lambda x)=\m^0(x)$, for all $\lambda>0.$
Since the norm in $X(\R^N\times \R^+;\R^3)$ is invariant under this scaling,
Theorem~\ref{thm:cauchy-LLG} yields the following result concerning the existence
of self-similar solutions.
\begin{corollary}\label{cor-self-similar}
	With the same notations and  hypotheses as in Theorem~\ref{thm:cauchy-LLG}, assume also that $\m^0$
	is homogeneous of degree zero. Then the solution
	$\m$ of \eqref{LLG-intro} provided  by Theorem~\ref{thm:cauchy-LLG} is forward self-similar.
	In particular there exists a smooth profile $\f : \R^N\to \SS^2$ such that
	$\m(x,t)=\f ({x}/{\sqrt t}),$
	for all $x\in\R^N$  and $t>0$.
\end{corollary}

Other authors have considered expanders for the harmonic map flow \eqref{HFHM} in different settings.
Actually, equation \eqref{HFHM} can be generalized for maps $\m : \mathcal{M}\times \R^+\to \mathcal{N}$,
with $\mathcal M$ and $\mathcal N$ Riemannian manifolds. Biernat and Bizo\'n \cite{biernat-bizon}
established results when  $\mathcal M=\mathcal N=\SS^d$ and $3\leq d\leq 6.$
Also, Germain and Rupflin \cite{germain-rupflin} have investigated the case $\mathcal M=\R^d$ and  $\mathcal N=\SS^d$,
in $d\geq 3$. In both works the analysis is done only for equivariant  solutions and does not  cover the  case
$\mathcal M=\R^N$ and $\mathcal{N}=\SS^2$.

\subsection{LLG with a jump initial data}
We want now to apply the well-posedness result 
to the self-similar solutions $\m_{c,\alpha}$
with initial conditions 
$\m^0_{c,\alpha}=\A^+_{c,\alpha}\chi_{\R^+}+\A^-_{c,\alpha}\chi_{\R^-}.$
Let us remark that the first term in the definition of $[\v]_{X}$ allows us to capture a blow-up rate of $1/\sqrt{t}$ for  $ {\| \nabla v(t)\|}_{L^{\infty}}$, as $t\to 0^+$.  This is exactly the blow-up rate for the self-similar solutions $\m_{c,\alpha}$.  The integral term in the semi-norm $[\cdot]_X$ is also well-adapted to these solutions. 	Indeed, for any $\alpha\in (0,1]$ and $c\geq 0$, we have
\begin{equation}
	\label{cond-4}
	[\m^0_{c,\alpha}]_{BMO}\leq  2c{\sqrt{2\pi}}/{\sqrt\alpha}
	\quad  \text{ and }\quad 
	[\m_{c,\alpha}]_{X}\leq {4 c}/{\alpha^{\frac{1}{4}}}.
\end{equation}

Let us start by considering a more general problem:  the LLG equation, in dimension one,  with a jump initial data given by $\m^0_{\A^{\pm}}=\A^+\chi_{\R^+}+\A^-\chi_{\R^-},$
where $\A^{\pm}$ are two given unitary vectors in $\SS^2$.
The smallness condition in the BMO semi-norm of $\m^0_{\A^{\pm}}$ is equivalent to the smallness of the angle between $\A^{+}$ and $\A^{-}$. From Theorem~\ref{thm:cauchy-LLG} we can deduce that the solution associated with $\m^0_{\A^{\pm}}$ is a rotation of a self-similar solution $\m_{c,\alpha}$ for an appropriate value of $c$. Precisely, 
\begin{thm}[\cite{gutierrez-delaire2}]
	\label{thm-small-angle}
	Let $\alpha\in (0,1]$. There exist $L_1, L_2>0$, $\delta^*\in(-1,0)$ and $\vartheta^*>0$  such that the following holds. Let  $\A^+$, $\A^- \in \SS^2$ and let $\vartheta$ be the angle between them.
	If $		0<\vartheta\leq \vartheta^*,$
	then there exists a solution $\m$ of \eqref{LLG-intro}
	with initial condition $\m^0_{\A\pm}$. Moreover, there exists $0<c<\frac{\sqrt{\alpha}}{2\sqrt{\pi}}$, such that
	$\m$ coincides up to a rotation with the self-similar solution  $\m_{c,\alpha}$, i.e.\
	there exists $\boR\in SO(3)$, depending only on $\A^+$, $\A^-$, $\alpha$ and $c$, such that
	$\m=\boR \m_{c,\alpha},$
	and $\m$ is the unique solution satisfying
	\begin{equation}\label{unicidad-self}
		\inf_{\substack{x\in \R\\ t>0}} m_3(x,t)\geq \delta^*
		\quad \textup{ and }\quad [\m]_{X}\leq L_1+L_2 c.
	\end{equation}
\end{thm}

A second consequence of Theorem~\ref{thm:cauchy-LLG} concerns the stability of the self-similar solutions. Indeed, 
from the dependence of the solution with respect to the initial data in this theorem  and the estimates in \eqref{cond-4}, we obtain the following  result: For any given $\m^0\in \SS^2$ close enough to  $\m^0_{\A^{\pm}}$, 
the solution $\m$ of \eqref{LLG-intro} associated with  $\m^0$ given by Theorem~\ref{thm:cauchy-LLG} must remain close to a rotation of a self-similar solution $\m_{c,\alpha}$, for some $c>0$. In particular, $\m$ remains close to a self-similar solution.
The precise statement is provided in the following theorem.

\begin{thm}[\cite{gutierrez-delaire2}]\label{thm:stability}
	Let $\alpha\in (0,1]$. There exist constants $L_1,L_2,L_3>0$, $\delta^*\in (-1,0)$,
	$\vartheta^*>0$ such that the following holds.
	Let $\A^+$,  $\A^- \in \SS^2$ with angle $\vartheta$ between them. If
	$0<\vartheta\leq \vartheta^*,$
	then there is $c>0$ such that for every
	$\m^0$ satisfying
$		\norm{\m^0-\m^0_{\A^\pm} }_{L^\infty}\leq \frac{c\sqrt{\pi}}{2\sqrt{\alpha}},$
	there exists $\boR\in SO(3)$, depending only on $\A^+$, $\A^-$, $\alpha$ and $c$, such that
	there is a unique global smooth solution $\m$ of \eqref{LLG-intro} with initial condition $\m^0$
	that satisfies
	\begin{equation}\label{condition}
		\inf_{\substack{x\in \R\\ t>0}} {(\boR \m)_3(x,t)}\geq \delta^*
		\quad \textup{ and }\quad [\m]_{X}\leq L_1+L_2 c.
	\end{equation}
	Moreover,
$		\norm{ \m-\boR\m_{c,\alpha}}_{X}\leq  L_3\norm{\m^0-\m^0_{\A^\pm} }_{L^\infty}.$
	In particular,
	\begin{equation*}
		\norm{\partial_x \m-\partial_x \boR\m_{c,\alpha}}_{L^\infty}\leq  \frac{L_3}{\sqrt{ t}}\norm{\m^0-\m^0_{\A^\pm} }_{L^\infty}, \quad \text{for all }t>0.
	\end{equation*}
\end{thm}

Let us now discuss the multiplicity of solutions with initial condition $\m^0_{\A^\pm}$.
As seen before, when $\alpha=1$, the self-similar solutions are explicitly given by
\eqref{explicit-for} and limit vectors $\vec \A^\pm_{c,1}$ given in \eqref{limit-alpha-1}.

\begin{figure}[ht!]
	\begin{center}
		\begin{overpic}[trim=0 0 30 0,clip,scale=0.5]{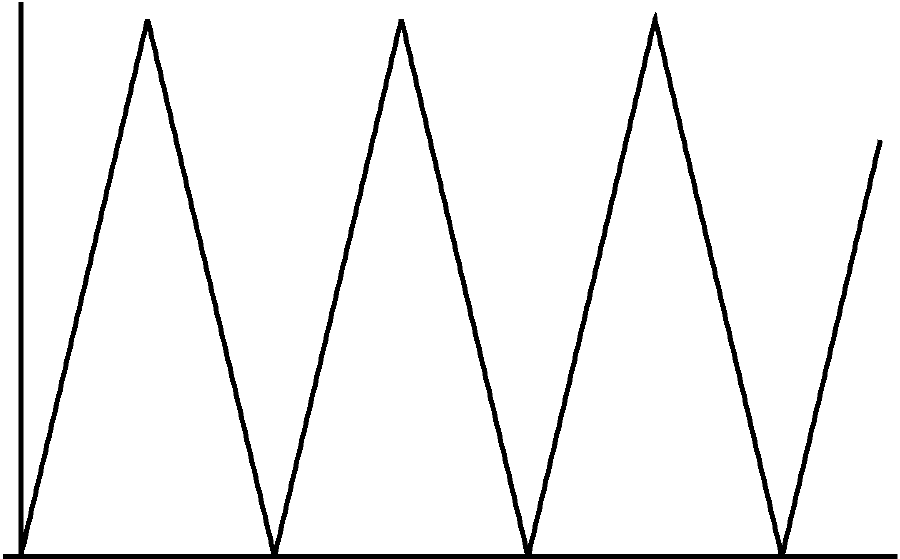}
			\put(4,71){\small{$\vartheta_{c,1}$}}
			\put(-4,64){\small{$\pi$}}
			\put(45,-6){\small{$c$}}
		\end{overpic}
	\end{center}
	\caption{The angle $\vartheta_{c,\alpha}$ as a function of $c$ for $\alpha=1$.}
	\label{fig-angulo}
\end{figure}

Figure~\ref{fig-angulo} shows that there are infinite values of $c$
that allow to reach any angle in $[0,\pi]$.
Therefore, using the  invariance of \eqref{LLG-intro} under rotations, in the case when $\alpha=1$, one can easily prove the existence of multiple solutions associated with a given initial data of the form $\m^{0}_{\A^{\pm}}$ for any given vectors $\A^{\pm}\in \mathbb{S}^2$. In the case that $\alpha$ is close enough to $1$, we can use a continuity argument
to prove that we still have multiple solutions. More precisely, we can establish that  for any given initial data of the form $\m^{0}_{\A^{\pm}}$, with angle between $\A^{+}$ and $\A^{-}$ in the interval $(0,\pi)$, if $\alpha$ is sufficiently close to one, then there exist {\it{at least}} $k$-distinct solutions of \eqref{LLG-intro} associated with the same initial condition, for any $k\in\mathbb{N}$. In other words,  given any angle $\vartheta\in(0,\pi)$ between two  $\A^{+}$ and $\A^{-}$, we can generate any number of distinct solutions by considering values of $\alpha$ sufficiently close to $1$. Precisely, 
\begin{thm}[\cite{gutierrez-delaire2}]\label{thm-non-uniq}
	Let $k\in \N$, $\A^+$, $\A^-\in \SS^2$ and let $\vartheta$ be the angle between $\A^+$ and $\A^-$.
	If $\vartheta\in (0,\pi)$, then there exists $\alpha_k \in (0,1)$ such that for every $\alpha\in [\alpha_k,1]$
	there are at least $k$ distinct smooth self-similar solutions $\{\m_j\}_{j=1}^k$ in $X(\R\times \R^+;\SS^2)$ of \eqref{LLG-intro}
	with initial condition $\m^0_{\A^\pm}$. These solutions are characterized by a  strictly increasing sequence
	of values $\{c_j\}_{j=1}^k$, with $c_k\to\infty$ as $k\to\infty$, such that
	$	\m_j=\boR_j \m_{c_j,\alpha},$
	where $\boR_j \in SO(3)$. In particular
	\begin{equation}\label{derivada-mk}
		\sqrt{t}\norm{\partial_x\m_j(\cdot, t)}_{L^\infty(\R)}=c_j, \quad  \text{for all } t>0.
	\end{equation}
	Furthermore, if $\alpha=1$ and  $\vartheta\in [0,\pi]$, then
	there is an infinite number of distinct smooth self-similar solutions  $\{\m_j\}_{j\geq 1}$ in $X(\R\times \R^+;\SS^2)$ of \eqref{LLG-intro}
	with initial condition $\m^0_{\A^\pm}$.
\end{thm}

It is important to remark that in particular Theorem~\ref{thm-non-uniq} asserts that when $\alpha=1$,
given $\A^+,\A^-\in\mathbb{S}^2$ such that $\A^+=\A^-$,  there exists an infinite number of distinct solutions $\{\m_j\}_{j\geq 1}$ in $X(\mathbb{R}\times\mathbb{R}^{+}; \mathbb{S}^2)$ of \eqref{LLG-intro} with initial condition $\m^0_{\A^\pm}$ such that $[\m^0_{\A^{\pm}}]_{BMO}=0$. This particular case shows that a condition on the size of $X$-norm of the solution  in Theorem~\ref{thm:cauchy-LLG} 
is necessary for the uniqueness of solution. 
We recall that for finite energy solutions of \eqref{HFHM}, there are several nonuniqueness results  based on Coron's technique \cite{coron}
in dimension $N=3$.
Alouges and Soyeur~\cite{alouges-soyeur} successfully adapted this idea to prove the existence of multiple solutions of \eqref{LLG-intro}, with $\alpha>0$,  for maps $ {\boldsymbol{m}}:\Omega\longrightarrow \mathbb{S}^2$, with $\Omega$ a bounded regular domain of $\mathbb{R}^3$. In our case, since  $\{c_j\}_{j=1}^k$ is strictly increasing, we have at least $k$ different \emph{smooth}
solutions. Notice also that  the identity \eqref{derivada-mk} implies that the $X$-norm of the solution is large as $j\rightarrow\infty$.

\subsection{Shrinkers}
\label{sec:shrinkers}
We end this note by discussing the backward self-similar solutions to \eqref{LLG-intro}, i.e.\ the shrinker solutions of  the form 
$	\m(x,t)=\f \left(\frac{x}{\sqrt{T-t}}\right)$, for $	x\in \R$ and $t\in (-\infty,T)$.
As in Section~\ref{self-similar}, we can reduce our problem 
to the study of the ODE
\begin{equation}
	\label{eq:EDO:intro-2}
	\alpha \f''+\alpha\abs{\f'}^2\f +\beta (\f\times \f')'-\frac{x \f'}2=0,
	\quad\textup{ on }\R,
\end{equation}
which is the same equation that we obtained for the expanders, 	except for the minus sign in the last term.
Following similar arguments, we get
\begin{thm}[\cite{gutierrez-delaire3}]
	\label{thm:EDO-2}
	Let $\alpha\in (0,1]$. Assume that $\f\in H^{1}_{\loc}(\R;\SS^2)$ is a weak solution to \eqref{eq:EDO:intro-2}. 
	Then $\f$ belongs to $\boC^\infty(\R;\SS^2)$ and there exists $c\geq 0$ such that
	$\abs{\f'(x)}=ce^{\alpha x^2/4}$, for all $x\in \R$. 
	Moreover, the set of nonconstant solutions to 
	\eqref{eq:EDO:intro-2} is $\{\boR \f_{c,\alpha} : c>0,\boR \in SO(3)\}$,
	where is  $\f_{c,\alpha}$ is given by the solution $\{\f,\g,\h\}$ of the Serret--Frenet system
	with curvature 		$k(x)= c e^{{\alpha x^2}/{4}}$ and torsion
$ \tau(x)=-{\beta x}/{2},$
	and initial conditions
	$	{\f}(0)=\be_1$,	${\g}(0)=\be_2$, and $ {\h}(0)=\be_3.$
\end{thm}

As done for the expanders, we provide now some 
properties of these solutions, that are obtained by studying the Serret--Frenet system.
\begin{thm}[\cite{gutierrez-delaire3}]
	\label{thm-conv}
	Let $\alpha\in(0,1]$, $c>0$, $T\in \R$
	and ${\f}_{c,\alpha}$ as above.
	Set
$		\tilde\m_{c, \alpha}(x,t) =\f_{c,\alpha}\left( \frac{x}{\sqrt{T-t}}  \right)$, for $x \in\R$, $t<T$. 	Then $\tilde\m_{c, \alpha}$ belongs to  $\mathcal{C}^\infty( \R\times (-\infty,T);\mathbb{S}^2)$,  solves  \eqref{LLG-intro} for $t\in (-\infty,T)$,
		and
$\abs{\partial_{x} \tilde\m_{c, \alpha}(x,t)}=\frac{c}{\sqrt{T-t}}e^{\frac{\alpha x^2}{4(T-t)}},$ 
	for all $(x,t)\in \R\times(-\infty,T).$
		Moreover, the following properties hold.
		\begin{enumerate}
		\item\label{parity}
		The first component of  $\f_{c,\alpha}$ is even,  while the others are  odd.
		\item\label{asymptotics}
		There exist constants $\rho_{j,c,\alpha}\in[0,1],$
		$B_{j,c,\alpha}\in [-1,1],$
		and 
		$\phi_{j,c,\alpha} \in [0,2\pi),$
		for $j\in \{1,2,3\}$,
		such that
		we have the  following asymptotics for the profile   $\f_{c,\alpha}$:
		\begin{equation}
			\label{eq:asymp}
			\begin{aligned}
				f_{j,c,\alpha}(x)=&\rho_{j,c,\alpha}\cos(c\Phi_\alpha(x)-\phi_{j,c,\alpha})
				-\frac{\beta B_{j,c,\alpha}}{2c} xe^{-\alpha x^2/4}\\
				&+
				\frac{\beta^2 \rho_{j,c,\alpha}}{8c}\sin(c\Phi_\alpha(x)-\phi_{j,c,\alpha})
				\int_x^\infty s^2 e^{-\alpha s^2/4}ds
				+\frac{\beta}{\alpha^5c^2}\boO(x^2e^{-\alpha x^2/2}),
			\end{aligned}
		\end{equation}
		for all $x\geq 1$, where $	\Phi_\alpha(x)=\int_{0}^{x} e^{\frac{\alpha s^2}{4}}\, ds.$
%
		\item\label{convergence}
		The solution $\tilde\m_{c,\alpha}=(\tilde m_{1,c,\alpha},\tilde m_{2,c,\alpha},\tilde m_{3,c,\alpha})$ satisfies the following pointwise convergences
		\begin{equation}
			\label{convergencia:shrin}
			\begin{aligned}
				\lim_{t\to T^-}(\tilde m_{j,c,\alpha}(x,t)-\rho_{j,c,\alpha}\cos\big(  c\, \Phi_{\alpha} \big( \frac{x}{\sqrt{T-t}}\big) -\phi_{j,c,\alpha}  \big)=0, \text{ if } x>0,\\
				\lim_{t\to T^-}(\tilde  m_{j,c,\alpha}(x,t)-\rho_{j,c,\alpha}^{-}\cos\big(  c\, \Phi_{\alpha} \big( \frac{-x}{\sqrt{T-t}}\big) -\phi_{j,c,\alpha}  \big)=0, \text{ if }x<0,
			\end{aligned}
		\end{equation}
		for $j\in \{1,2,3\}$, where $\rho_{1,c,\alpha}^-=\rho_{1,c,\alpha}$, 
		$\rho_{2,c,\alpha}^-=-\rho_{2,c,\alpha}$ and $\rho_{3,c,\alpha}^-=-\rho_{3,c,\alpha}$.
		\item\label{IVP} $\tilde  \m_{c,\alpha}(\cdot, t)\to 0$ as $t\to T^-$, as a tempered  distribution.
	\end{enumerate}
\end{thm}

\medskip

%

As for the expanders,   the big-$\boO$ in the asymptotics does not depend on $\alpha \in[0,1]$.
In this manner, the constants multiplying the big-$\boO$ are meaningful and in particular,  big-$\boO$  
vanishes when $\beta=0$. 
Let us remark that the behavior of the profile for $x\leq  -1$ follows 
from the symmetries of the profile established in part \ref{parity}.

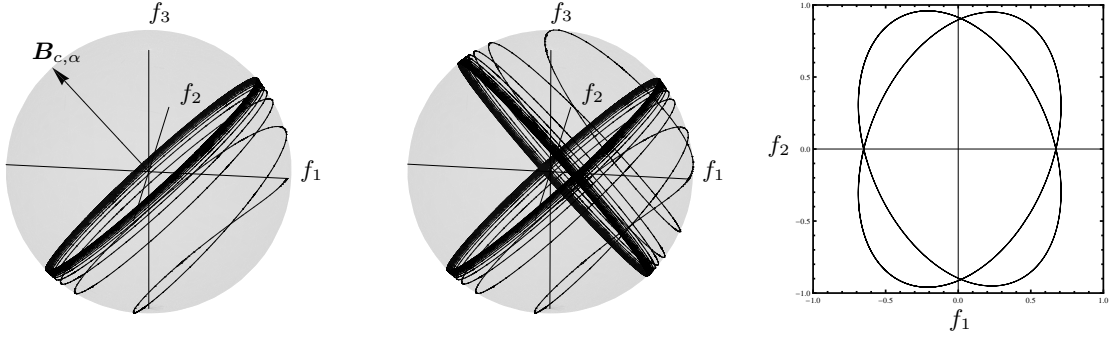
\begin{figure}[h]
	\begin{subfigure}[c]{0.33\textwidth}
		\begin{overpic}[scale=0.5]{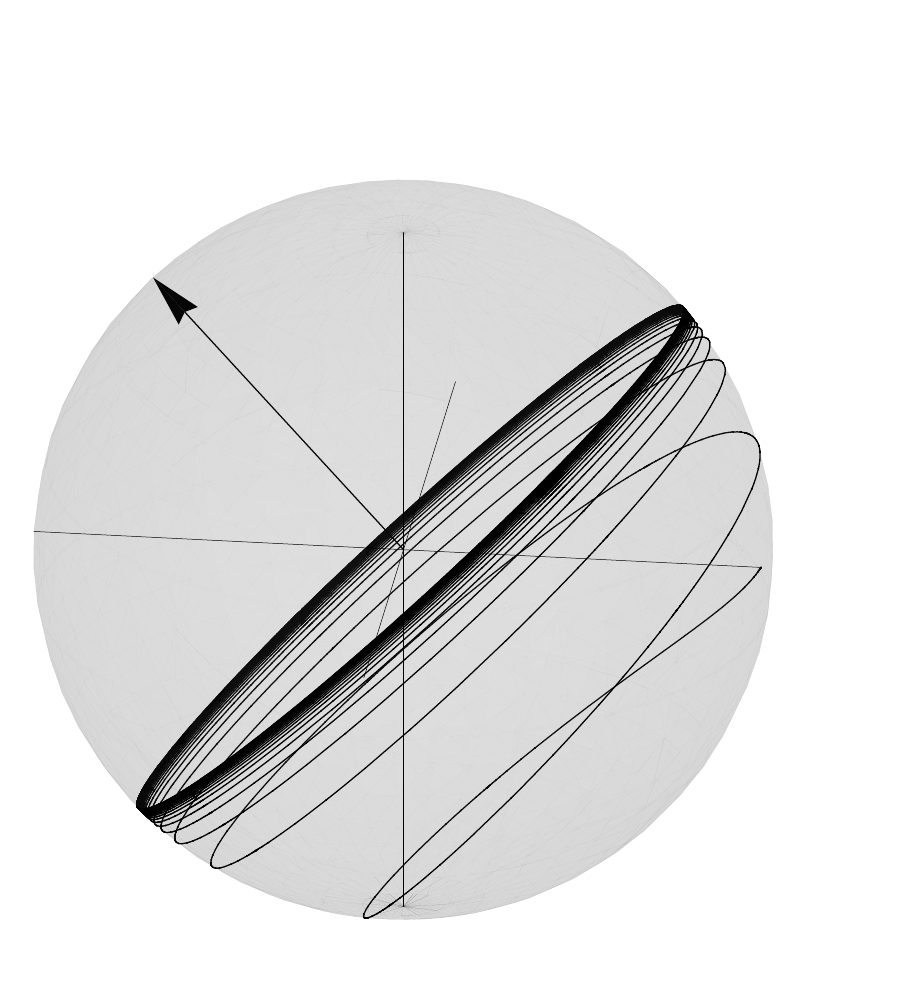}
			\put(79,44){\footnotesize{$f_1$}}
			\put(48,64){\footnotesize{$f_2$}}
			\put(40,85){\footnotesize{$f_3$}}
			\put(10,74.5){\footnotesize{$\bm B_{c,\alpha}$}}
		\end{overpic}
	\end{subfigure}%
	\begin{subfigure}[c]{0.33\textwidth}
		\begin{overpic}[scale=0.5]{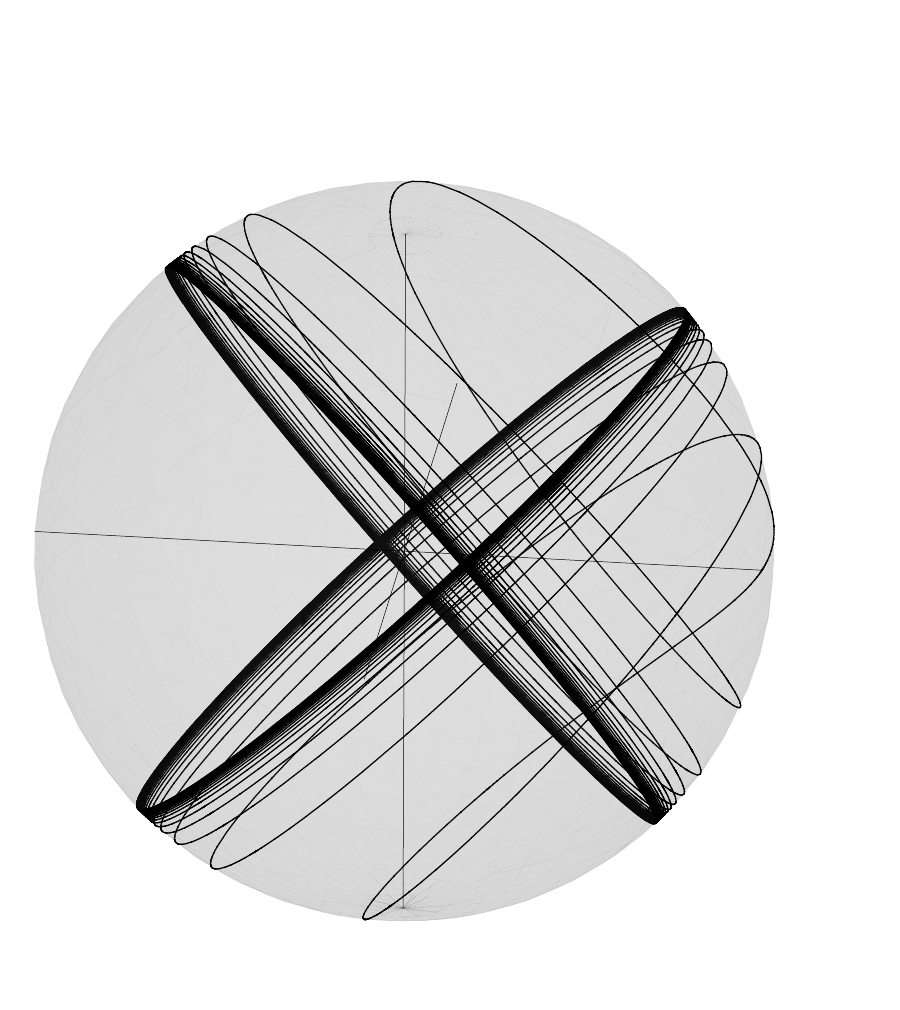}
			\put(79,44){\footnotesize{$f_1$}}
			\put(48,64){\footnotesize{$f_2$}}
			\put(40,85){\footnotesize{$f_3$}}
		\end{overpic}
	\end{subfigure}%
	\begin{subfigure}[c]{0.33\textwidth}
		\begin{overpic}[scale=0.45]{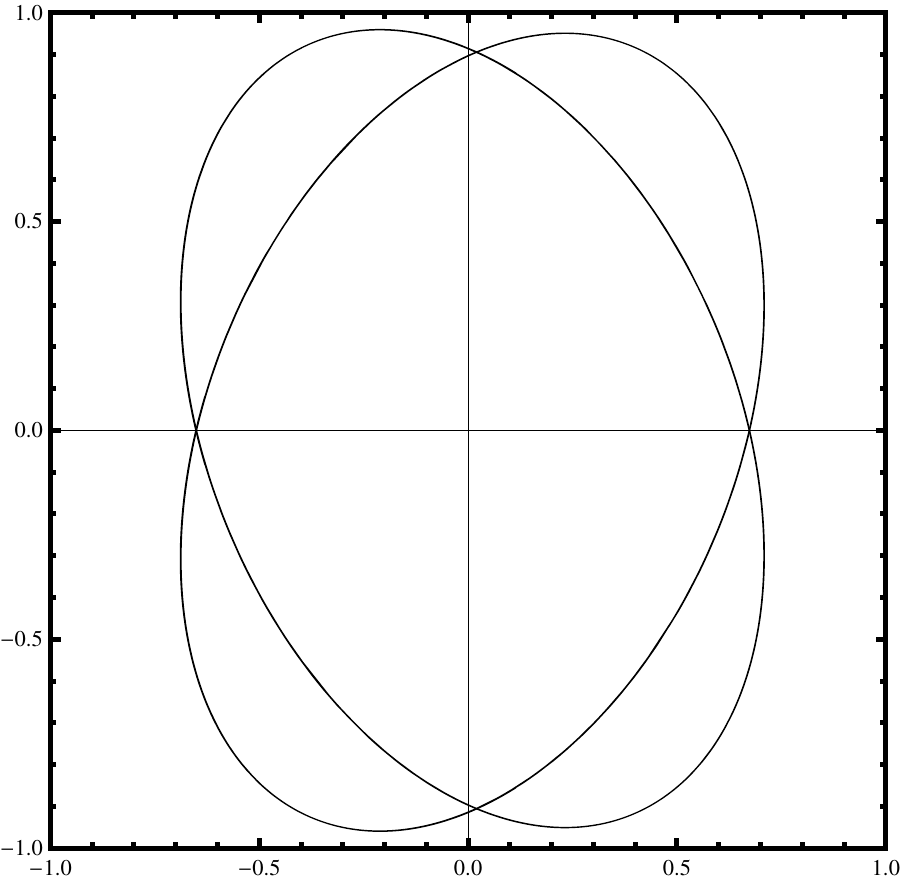}
			\put(49,-7){\footnotesize{$f_1$}}
			\put(-9,49){\footnotesize{$f_2$}}
		\end{overpic}
	\end{subfigure}%
	\caption{Profile $\f_{c,\alpha}$ for $c=0.5$ and $\alpha=0.5$. The figure on the left depicts profile for $x\in \R^+$ and the normal vector $\bm B_{c,\alpha}\approx (-0.72, -0.3, 0.63)$. The figure on the center  shows  the profile for $x\in\R$; the angle between the circles $\boC_{c,\alpha}^\pm$ is  $\vartheta_{c,\alpha}\approx 1.5951$.
		At the right, the projection of limit cycles $\boC_{c,\alpha}^\pm$ on the plane $\R^2$.}
	\label{fig-curvas}
\end{figure}

In Figure~\ref{fig-curvas}, we have depicted the profile $\tilde\m_{c, \alpha}$ 
for $\alpha=0.5$ and $c=0.5$, where we can see the oscillating behavior.
Moreover, the plots in Figure~\ref{fig-curvas} suggest  that the 
limit sets of the trajectories are great circles on the sphere $\SS^2$ when $x\to\pm\infty$. 
This is indeed the case. The next result 
establishes  
analytically that $\tilde \m_{c,\alpha}$ 
oscillates in a plane passing through the origin whose normal vector is given by 
$\bm B_{c,\alpha}^\pm=(B_{1,c,\alpha},B_{2,c,\alpha},B_{3,c,\alpha})$, as  $x\to \pm \infty$, respectively.

\begin{thm}[\cite{gutierrez-delaire2}]
	\label{plane}
Let  $\boP_{c,\alpha}^{\pm}$ be 
	the planes passing through the origin with normal vectors $\bm B_{c,\alpha}^\pm$,
	respectively.
	Let
	$\boC_{c,\alpha}^\pm$ be the circles in $\R^3$ given by
	$\boC_{c,\alpha}^\pm=\boP_{c,\alpha}^{\pm}\cap \SS^2.$
	Then for all $\abs{x}\geq 1$,
		\begin{equation}
			\label{dist1}
			\dist(\tilde \m_{c,\alpha}(x),\boC_{c,\alpha}^{\pm}))\leq \frac{15\sqrt{2}\beta}{c \alpha^2 }\abs{x}e^{-\alpha x^2/4}.
		\end{equation}
		In particular
		\begin{equation*}
			\begin{aligned}
				\lim_{t\to T^-}\dist(\tilde\m_{c,\alpha}(x,t),\boC_{c,\alpha}^+))=0, \text{ if }x>0,
				\text{ and }
				\lim_{t\to T^-}\dist(\tilde\m_{c,\alpha}(x,t),\boC_{c,\alpha}^-))=0,  \text{ if }x<0.
			\end{aligned}
		\end{equation*}
\end{thm}
Theorem~\ref{plane}   establishes the  convergence of the profile  $\f_{c,\alpha}$ to the great circles $\boC_{c,\alpha}^\pm$ as shown in Figure~\ref{fig-curvas}. Moreover, \eqref{dist1}
gives us an exponential rate for this convergence.
In terms of the solution $\tilde\mm_{c,\alpha}$ to the LLG equation, this provides a more precise 
geometric information about the way that the solution blows up at time $T$. The existence of limit cycles for related ferromagnetic models have been 
investigated for instance in \cite{waldner,broggi-meier}, but to the best of our knowledge 
this is  the first time that this type  of phenomenon  has been observed  for the LLG equation.
In Figure~\ref{fig-curvas} one can see that   $\vartheta_{c,\alpha}\approx 1.5951$ for $\alpha=0.5$
and $c=0.5$, where we have chosen the value of $c$ such that the angle is close to $\pi/2$.

In the case  $\alpha=1$, the torsion vanishes, and it easy to deduce that
the profile is explicitly given by the plane curve 
$\f_{c,1}(x)=(\cos(c\Phi_1(x)),\sin(c\Phi_1(x)), 0)$.
In particular, we see  that the asymptotics in Theorem~\ref{thm-conv} are satisfied with
$\rho_{1,c,1}=1$, $\rho_{2,c,1}=1$, $\rho_{3,c,1}=0,$
$\phi_{1,c,1}=0$, $\phi_{2,c,1}=3\pi/2$, $\phi_{3,c,1}\in [0,2\pi).$

In the  case $\alpha=0$,  $\f_{c,0}$ is equal to $\gm_{c,0}$ in  \eqref{asymp-alpha-0}, so that   $\f_{c,0}$ converges to the point $\A^+_{c,0}$,  as $x\to\infty$.
Hence, there is a drastic change  in the behavior of the profile in the
cases $\alpha=0$ and $\alpha>0$: 
In the first case $\f_{c,0}$ converges to a point at infinity, 
while in the second case \eqref{dist1} tells us that $\f_{c,\alpha}$ converges to a great circle. In this sense,  there is a discontinuity in the behavior of
$\tilde \m_{c,\alpha}$ at $\alpha=0$.


\noindent {\bf{Acknowledgements. }} 
A.~de Laire was partially supported by the Labex CEMPI
(ANR-11-LABX-0007-01), the ANR project ``Dispersive and random waves'' (ANR-18-CE40-0020-01), and the MATH-AmSud project EEQUADD-II.
 \bibliographystyle{abbrv}

\end{document}